\documentclass[tbtags,reqno]{amsart}

\usepackage{amssymb,enumitem,cancel}
\usepackage{geometry}
\usepackage[foot]{amsaddr}
 
\usepackage{mathtools,bbm}
\usepackage{mathrsfs}
\usepackage{xfrac}
\usepackage{marginnote}
\usepackage{dsfont}
\usepackage{changepage}
\usepackage{xcolor}
\usepackage{soul}
\usepackage{extdash}
\usepackage{tikz-cd}
\usepackage{resmes}
\usepackage{xltabular}
\usepackage{scalerel,stackengine}
\stackMath
\newcommand\wwidehat[1]{%
\savestack{\tmpbox}{\stretchto{%
  \scaleto{%
    \scalerel*[\widthof{\ensuremath{#1}}]{\kern-.6pt\bigwedge\kern-.6pt}%
    {\rule[-\textheight/2]{1ex}{\textheight}}
  }{\textheight}%
}{0.5ex}}%
\stackon[1pt]{#1}{\tmpbox}%
}
\usepackage{hyperref}
\usepackage{cleveref}
\numberwithin{equation}{section}

\newtheorem{thm}{Theorem}[section]
\newtheorem{prop}[thm]{Proposition}
\newtheorem{lem}[thm]{Lemma}
\newtheorem{cor}[thm]{Corollary}
\crefname{thm}{Theorem}{Theorems}
\crefname{prop}{Proposition}{Propositions}
\crefname{lem}{Lemma}{Lemmas}
\crefname{cor}{Corollary}{Corollaries}

\crefname{claim}{Claim}{Claims}
\crefname{step}{Step}{Steps}

\theoremstyle{definition}
\newtheorem{defn}[thm]{Definition}
\newtheorem{exa}[thm]{Example}
\newtheorem{asss}[thm]{Assumptions}
\crefname{asss}{Assumptions}{Assumptions}
\crefname{defn}{Definition}{Definitions}
\crefname{exa}{Example}{Examples}

\theoremstyle{remark}
\newtheorem{rmk}[thm]{Remark}
\crefname{rmk}{Remark}{Remarks}

\AddToHook{env/prop/begin}{\crefalias{thm}{prop}}
\AddToHook{env/lem/begin}{\crefalias{thm}{lem}}
\AddToHook{env/cor/begin}{\crefalias{thm}{cor}}
\AddToHook{env/claim/begin}{\crefalias{thm}{claim}}
\AddToHook{env/step/begin}{\crefalias{thm}{step}}
\AddToHook{env/defn/begin}{\crefalias{thm}{defn}}
\AddToHook{env/exa/begin}{\crefalias{thm}{exa}}
\AddToHook{env/rmk/begin}{\crefalias{thm}{rmk}}
\AddToHook{env/asss/begin}{\crefalias{thm}{asss}}

\makeatletter 
    \def\@makefnmark{\@textsuperscript{\normalfont\texttt{\@thefnmark}}} 
\makeatother

\newcommand{\norm}[1]{\ensuremath\lVert#1\rVert}
\newcommand{\compl}{^\mathsf{c}}
\newcommand{\defeq}{\vcentcolon=}

\renewcommand{\leq}{\leqslant}
\renewcommand{\geq}{\geqslant}

\let\temp\phi
\let\phi\varphi
\let\varphi\temp
\let\temp\varepsilon
\let\varepsilon\epsilon
\let\epsilon\temp

\newcommand{\tildee}[1]{\widetilde{#1}}

\DeclareMathOperator{\Lip}{Lip}

\DeclareMathOperator{\tr}{tr}

\newcommand{\N}{\mathbb{N}}

\newcommand{\de}{\partial}
\newcommand{\di}{\mathrm{d}}
\newcommand{\R}{{\mathbb{R}}}

\newcommand{\call}[1]{\ensuremath\mathcal{#1}}
\newcommand{\frk}[1]{\ensuremath\mathfrak{#1}}
\newcommand{\scr}[1]{\ensuremath\mathscr{#1}}
\newcommand{\bb}[1]{\ensuremath\mathbb{#1}}

\newcommand{\dto}{\ensuremath{\hspace{0.8pt}\downarrow\hspace{0.8pt}}}
\newcommand{\uto}{\ensuremath{\hspace{0.8pt}\uparrow\hspace{0.8pt}}}
\newcommand{\var}{\mspace{1.8mu}\cdot\mspace{1.8mu}}
\newcommand\loc{\mathrm{loc}}
\newcommand{\Sym}{\scr S}

\newcommand{\1}{\mathbbm{1}}

\newcommand{\GLN}{\mathrm{GL}_N}
\newcommand{\TV}{\mathrm{TV}}
\newcommand{\M}{\call R}
\newcommand{\Ran}{\mathrm{im}}

\newcommand\mres\llcorner

\tikzset{%
  symbol/.style={
    draw=none,
    every to/.append style={
      edge node={node [sloped, allow upside down, auto=false]{$#1$}}
    },
  },
}

\def\be{\begin{equation}}
\def\ee{\end{equation}}
\def\dive{{\rm div}}
\def\rife#1{\eqref{#1}}

\def\cL{{\mathcal L}}

\numberwithin{equation}{section}

\begin{document}

\date{\today}

\title[Hypocoercivity estimates via maximum principle in short and long time scales]{Hypocoercivity estimates via maximum principle \\ 
in short and long time scales}
\author[A.\ Porretta]{Alessio Porretta\textsuperscript{(1)}}
\address{\textsuperscript{(1)}Department of Mathematics \\ University of Rome Tor Vergata \\ Via della Ricerca Scientifica 1 \\ 00133 Roma, Italy.}
\author[D.\ F.\ Redaelli]{Davide Francesco Redaelli\textsuperscript{(1)}}
\email{porretta@mat.uniroma2.it, redaelli@mat.uniroma2.it}

\keywords{hypoelliptic diffusions, hypocoercivity, H\"ormander conditions, fractional diffusions, kinetic equations, long-time decay, smoothing estimates, coupling method, doubling variables, maximum principle \\ 
{\bf AI Statement:} No AI tool, nor any LLM model, have been used, in any form, in the conception, elaboration and writing of the present paper. 
\\
{\bf Funding:} We acknowledge the support of the Excellence Project MatMod@TOV
of the Department of Mathematics of the University of Rome Tor Vergata and  of  Italian
(EU Next Gen) PRIN project 2022W58BJ5 (\emph{PDEs and optimal control methods in mean field games, population dynamics and multi-agent models}, CUP E53D23005910006). The authors are members of INdAM (Istituto Nazionale di Alta Matematica) and received support of annual GNAMPA research projects.}


\begin{abstract} We prove decay estimates for a large class of linear drift-diffusion operators with hypoelliptic structure. In the main example, the diffusion consists of (possibly superposing) local or nonlocal operators  acting on a  restricted set of variables, while diffusivity is transported to the whole space thanks to H\"ormander-type conditions on the drift term.
We prove both short-time and long-time oscillation estimates on the solutions, including bounded-to-Lipschitz smoothing in the  sharp time scale and  exponential decay, as time goes to infinity, whenever a strict supersolution, with the role of Lyapunov function, is supposed to exist. The estimates are proved with a unifying approach relying on 
doubling variables and coupling methods, using  maximum principle in  suitably time scaled moving frames. 
Applications range from nonlocal kinetic Fokker--Planck equations to heat conduction models in a chain  of coupled oscillators. 
\end{abstract}

\maketitle

\setcounter{tocdepth}{2}
\tableofcontents

\section{Introduction}
This paper is concerned with both short-time and long-time oscillation estimates for linear drift-diffusion equations such as
\be\label{TE}
\partial_t u = \cL u + H(t,x) \cdot Du,\qquad u=u(t,x), \quad t>0,\ x\in \R^N,
\ee
in which the diffusion operator $\cL$ is supposed to act only on a subspace of $\R^N$, say on a given subset of variables, while the drift field $H$ satisfies some structure condition, such as \emph{H\"ormander-type conditions}, that spread diffusivity to the whole set of variables. 
As a byproduct of these oscillation estimates, we prove the exponential convergence, as time goes to infinity, to stationary measures, for solutions of the adjoint conservative equation
\be\label{FPE}
\partial_t m = \cL^* m - \dive( H(t,x) m),\qquad  t>0,\ x\in \R^N\,.
\ee
Typical applications of our results include -- under suitable confining conditions -- the long-time stabilization for kinetic Fokker--Planck equations such as
\[
\partial_t m + v\cdot D_\xi m +  (-\Delta_v)^{\frac\sigma2} m + \dive_v( \Phi(\xi,v)m)=0
\]
where $m= m(t,\xi,v)$ and $ \sigma\in (0,2]$ denotes the degree of (possibly fractional) diffusion. However, much more general examples are included in our framework, like for instance if $x=(x^0,x^1, x^2, x^3, x^4) \in (\R^d)^5$ and 
\[
\cL u + H(x) \cdot Du = \Delta_{x^0} - (-\Delta_{x^1})^{\frac\sigma2} + x^1D_{x^2} + x^0D_{x^3} + x^3 D_{x^4},
\]
in which case the diffusivity in $x^0$ is transported recursively  to $x^3$ and $x^4$, while the diffusivity in $x^1$  is transported to $x^2$. By the way, the transport drift $H$ does not need to be linear in $x$, but this is just the simplest model.

\medskip

As it may be clear from the above two examples, our purpose is to present a unifying approach to both short-time smoothing and long-time contractivity estimates that is meant to be sufficiently  general to be applied to possibly different diffusive operators, local and/or nonlocal, while including hypoelliptic structures under H\"ormander-type conditions.
 
By short-time smoothing, we mean an estimate such as 
\[
|u(t,x) - u(t,y)| \leq  C_T 
\bigl(\phi(t,x)+\phi(t,y)\bigr) \bigl\lvert \M(t)^{-1} (x-y) \bigr\rvert \qquad \forall\,t \in (0,T], \ x,y \in \R^N,
\]
for a suitable matrix $\M(t)$ (which keeps trace of time scaling) and some weight function $\phi$. Here the constant $ C_T$  may depend on time horizon $T$ as well as on the $L^\infty$-norm of $\phi^{-1}u$.  
The reader may just take $\phi=1$ at first glance, though possibly general Lyapunov-type functions $\phi$ are relevant to take care of the growth at infinity.  

By long-time contractivity, we mean, for instance, an exponentially decaying estimate such as 
\[
[u(t,\var)]_{\phi,\theta} \leq K e^{-\varpi t} [u(0,\var)]_{\phi,\theta} \qquad \forall\,t\in (0,\infty),
\]
where now $K$ and $\varpi$ are independent of $T$  and 
\begin{equation} \label{eq_defoscintro}
[u(t,\var)]_{\phi,\theta} \defeq \sup_{x,y \in \R^N} \frac{|u(t,x) - u(t,y)|}{\bigl(\phi(x)+\phi(y)\bigr)(|x-y|^\theta \wedge 1)}, \quad \theta \in (0,1],
\end{equation}
is a weighted oscillation semi-norm of $u$.

\medskip

There are a number of motivations, as well as connections to previous works, that lead to our study. 

On one hand, there is the long-time stabilization of Fokker--Planck equations, which has been widely investigated with both analytic and probabilistic methods, linked to the ergodic theory of stochastic processes. Typical approaches rely on Harris' type theorems~\cite{Canizo+al,Eberle+al2, HM, Yoldas}, spectral gap methods \cite{Gualdani+al,MM} or variational tools based on entropy/energy decay \cite{Bakry+al, MaVi}; for kinetic equations or hypoelliptic diffusions, the concept of hypocoercivity  is nowadays a milestone \cite{Desvillettes+V, DMS, H-Nier, herau, herau-nier, Vil, Vil26}. 

It is known that, when the operator $\cL$ in \rife{FPE} is diffusive, local and/or nonlocal, the long-time convergence to stationary measures can be reduced to a single estimate  establishing the  decay of oscillations of solutions of \rife{TE}. This viewpoint,  somehow  highlighted by Hairer and Mattingly~\cite{HM},  becomes effective by using purely probabilistic coupling methods~\cite{Eberle+al2} or analytic doubling variables arguments based on maximum principle \cite{PorAdv}. This approach was extended to the kinetic Fokker--Planck equation in \cite{Eberle+al} and then generalized in \cite{ForPor}, but still limited to what is called H\"ormander condition \emph{of first order}.

One of our goals is to extend  this approach and the results in \cite{ForPor} to the case of general H\"ormander conditions, including possibly (a mix of local and) nonlocal diffusions. When doing so, some form of hypocoercivity must clearly play a role; to this extent, our contribution lies in pointwise (oscillation, gradient) hypocoercivity estimates for the advection problem \rife{TE} which can be interesting in their own. Indeed, so far the classical hypocoercivity approach of Villani was developed essentially in Hilbertian settings, say through $L^2$ or $H^s$ norms, while we provide natural counterparts in $L^\infty$/oscillation and Lipschitz (semi-)norms.

On the other hand, there is the short-time smoothing property which shows how the hypoelliptic H\"ormander-type conditions guarantee that bounded initial data become Lipschitz at positive time. In this direction, along the viewpoint of stochastic processes, a lot of attention was devoted to cascade systems such as
\begin{equation} \label{intro-chosc} \begin{dcases}
\di X^0_t = H^0(t,X^0_t,\dots,X^r_t)\di t + \di Z_t
\\
\di X_t^1 = H^1(t,X^0_t,\dots,X^r_t)\di t\\
\di X_t^2 = H^2(t,X^1_t,\dots,X^r_t)\di t \\[2pt]
 \ \ \smash{\vdots} \\
\di X_t^r = H^r(t,X^{r-1}_t,X^r_t)\di t,
\end{dcases} \end{equation}
where the diffusivity appearing on a group of variables is transmitted to the other components through iterated non-degeneracy conditions on the drifts ($D_{x^{i-1}}H^i >0$). Models of this type appear in heat conduction problems, describing systems of coupled oscillators~\cite{Eckmann,Talay} as well as for option pricing problems in finance (see e.g. \cite{Del-Menozzi,Difra-Pasc, pigato} when $Z_t$ is a Brownian motion, or  \cite{Huang-menozzi} for $Z_t$ a jump process). In this latter context, diffusivity is often exploited through Malliavin's calculus and lower/upper  bounds for the density function, or for the kernel of the associated Markov operator. 

A second goal of our paper is to prove the bounded-to-Lipschitz smoothing property for similar dynamics, with a rather elementary PDE proof,  which may be interesting even for possibly more general transmission models of this type. 
Note that H\"ormander conditions appear in this kind of systems through the upper triangular form of the Jacobian matrix $DH$; this is in fact the main example of structure conditions on the drift $H$ that we have in mind in our results.  

\medskip

Even if we will give later more general --- though more abstract --- statements of our results, let us provide here two simple versions to serve as model statements of reference. For this purpose, it is convenient to introduce the following class of operators.

\begin{defn}[Standard diffusive operators on a subspace] \label{def_opintro}
Given a subspace $V \neq \{0\}$ of $\R^N$, we say that an operator $\call L \neq 0$ of the form 
\[
\call Lu = \tr(AD^2u) + \call I[u] 
\]
is a \emph{standard diffusive operator acting on $V$} if
\[
A \in \Sym_{\geq 0}^N, \quad \Ran\,A \in \{0,V\}
\]
and $\call I$ is given by
\[
\call I[u](t,x) = \int_{V} \Bigl( u(t,x+z) - u(t,x) - Du(t,x) \cdot z \1_{|z| \leq 1} \Bigr) \, \di\nu(z), 
\]
where either $\nu = 0$ or $\nu$ is a L\'evy measure supported on $V$ that satisfies 
\begin{equation} \label{nulb_intro}
\di\nu(z) \gtrsim \frac{\di z}{|z|^{d+\sigma}}, \qquad d = \dim V, \ \ \sigma \in (0,2).
\end{equation}
To such an operator $\call L$, we associate the two values
\[
\underline\sigma_{\cL} \defeq \inf\{ \sigma \in (0,2]:\, \text{\eqref{nulb_intro} holds} \}, \qquad \bar\sigma_\cL \defeq \sup\{ \sigma \in (0,2]:\,  \text{\eqref{nulb_intro} holds} \},
\]
with the convention that \eqref{nulb_intro} holds for $\sigma = 2$ if $\Ran\,A = V$.\footnote{Note that such a convention is formally motivated by the fact that $|\var|^{-d-\sigma}$ is the kernel of a fractional Laplacian of order $\sigma$ and by the well-known identity $\lim_{\sigma \uto 2} (-\Delta)^{\frac\sigma2} = -\Delta$.}
\end{defn} 

We point out that $\call L$ may contain, besides the local diffusion, more than one nonlocal diffusion, as, for instance, \eqref{nulb_intro} is verified by sums of fractional Laplacians of different orders $\sigma_j$ (also allowing $\sigma_j = 2$); in this case, we see that $\underline\sigma_{\cL} = \min_j \sigma_j$ and $\bar\sigma_{\cL} = \max_j\sigma_j$. 

In the particular case that $\R^N \simeq \R^{d_0} \times \cdots \times \R^{d_r}$ has coordinates $(x^0,\dots,x^{r})$, where $x^i \in \R^{d_i}$, and $V = \{ x \in \R^N :\, x^{-i} = 0 \} \simeq \R^{d_i}$ (where $x^{-i} \in \R^{N-d_i}$ denotes the vector $x$ deprived of the $i$-th coordinate), we will simply say that $\cL$ \emph{acts on the variable $x^i$}.
Then,   our first model statement reads as follows.

\begin{thm}
\label{intro-thm1}
Let $\R^N \simeq \R^{d_0} \times \cdots \times \R^{d_r}$ have coordinates $x=(x^0,\dots,x^r)$, where $x^i \in \R^{d_i}$. Let $\cL$ be a standard diffusive operator acting on the variable $x^0$; assume that the drift $H = (H^0,\dots,H^r)$ satisfies $DH \in L^\infty(\R^N)$ and, for each $i \in \{1,\dots,r\}$, 
\begin{equation} \label{eq_HhormC1intro}
D_{x^j} H^i = 0 \quad \text{for} \ j < i-1, \qquad \exists\, J_i \in \R^{d_{i-1} \times d_i}:\  \det(D_{x^{i-1}} H^i J_i) > 0 \ \text{uniformly on $\R^N$}.
\end{equation} 
Let $u$ be a bounded continuous (viscosity) solution to \eqref{TE}. Then:
\begin{enumerate}[leftmargin=2em]
\item given a weight function  $\phi$ satisfying conditions \eqref{lyapt1}--\eqref{eq_pshar}, there is a constant $C_T$, independent of $u$, such that 
\[
|u(t,x) - u(t,y)| \leq C_T \norm{\phi^{-1}u(0,\var)}_\infty \bigl(\phi(t,x)+\phi(t,y)\bigr)
 \bigl\lvert \call S(t)^{-1} (x-y) \bigr\rvert \qquad \forall\,t \in (0,T], \ x,y \in \R^N,
\]
where
\begin{equation} \label{eq_scalintro}
\call S(t) = t^{\bar\sigma_\cL^{-1}} \mathrm{diag}(t^{i}I_{d_i})_{0\leq i\leq r};
\end{equation}
\item if there exists a stationary Lyapunov function $\phi$ satisfying \eqref{eq_pshar} and \eqref{eq_condLyapdec}, for any $\theta \in (0,1]$ there are $\varpi>0$ and $K$, independent of $u$ and $t$, such that
\be\label{osc-decay-intro}
[u(t,\var)]_{\phi,\theta} \leq K e^{-\varpi t} [u(0,\var)]_{\phi,\theta} \qquad \forall\,t>0
\ee
where the oscillation semi-norm is defined as in \eqref{eq_defoscintro}.
\end{enumerate}
\end{thm}

Note that the second condition in \eqref{eq_HhormC1intro} is a uniform non-degeneracy assumption on $D_{x^{i-1}}H^i$, which ensures the transmission of diffusivity from the $x^{i-1}$ to the $x^i$ variable. In the particular case when $d_i = d$ for all $i$ (which is particularly common in cascade systems), this is equivalent to $D_{x^{i-1}}H^i$ being a square matrix with uniformly bounded inverse. On the other hand, a basic setting with different dimensions naturally occurs when considering general H\"ormander's conditions, or chains of heated oscillators (see  \Cref{model_LP}, or \Cref{chains}).

The key point  behind \Cref{intro-thm1} is to show how the mix of time scaling and coupling methods (possibly component-wise coupling) leads to both the short-time {\it bounded-to-Lipschitz} smoothing and the long-time decay of the oscillation.  As in \cite{PoPr}, coupling methods are here developed in the analytical version of doubling variables, only relying on maximum principle, where time scaling may act more efficiently. The use of maximum principle explains why the result is formulated for so-called viscosity solutions, which represent the natural class of weak solutions where the maximum principle continues to hold (in a relaxed form); see  \cite{CIL}. Such a formulation, together with natural stability properties, also guarantee that the estimates are stable in the small limit of diffusive regimes (so called vanishing viscosity approximation for the deterministic variables). Similarly, we point out that the boundedness assumption on $u$, in \Cref{intro-thm1}, can be relaxed by a stability argument; cf.~\Cref{rmk_onuLip}.

The main ideas of our approach are presented in \Cref{model_LP} for the simplest case of linear drifts with general H\"ormander conditions. First we show that, in the case of second order operators, the gradient decay can be equally proved through the classical Bernstein method, which is  instructive to highlight the spirit of our estimates.
 However, it would be unclear to extend the same approach to general nonlocal operators; this is why we opt for a coupling method, which  looks more flexible to mix local and nonlocal diffusions (cf.~\cite{PorAdv}). The way we employ the coupling approach is through splitting; the idea is that the operator should be reduced (eventually after twisting the variables) to the sum of a diffusive part, acting on given subspaces, plus a dissipative transport part acting on the non-diffused variables. This allows us to tune reflection or sinchronous couplings alternatively on two mutually orthogonal sets. In addition, time scaling is crucially used to enhance this structure at arbitrarily small scales.

Let us also observe that the oscillation decay \rife{osc-decay-intro} is just another form of hypocoercivity enjoyed by the drift-diffusion operator. This can take the more familiar form  of gradient decay, if we combine the short and long time scales, whenever one has additionally that $D_{x^j} H^i = 0$ for all $j\neq i-1$.

\begin{cor}\label{sharp-Kolm}
Given the assumptions of \Cref{intro-thm1}, suppose in addition that $D_{x^j} H^i = 0$ if $j \neq i-1$. Then there is a constant $K_0$ (independent of $u$) such that, for any $T>0$ and all $t \in (0,T)$,
\be\label{grad_t_decay}
\norm{D_{x^j}u(t,\var)}_\infty \leq \frac{K_0\norm{u}_\infty}{t^j(t^{\underline\sigma_\cL^{-1}} + t^{\bar\sigma_\cL^{-1}})}.
\ee
\end{cor}

The above corollary includes the polynomial decay for generalizations of the classical Kolmogorov operator.  In particular, when the diffusion is mixing local and nonlocal terms, the estimate \rife{grad_t_decay} shows how the possibly fat tails of fractional kernels affect the long-time scale, whereas the short-time scaling is governed by the fastest diffusion.

\Cref{intro-thm1} admits several possible extensions and we will provide abstract forms of our method later (see \Cref{thm_mainsmnl,prop_oscdec}). The general recipe requires first to split the variables according to the diffusive operators (and their most useful scaling laws), then to use suitable  transformations of coordinates and time scaling. In this regard, here is the second model statement that we propose, which in fact generalizes \Cref{intro-thm1}.

\begin{thm}
\label{intro-thm1gen}
Suppose there exists a decomposition $\R^N \simeq \R^{d_{1,0}} \times \cdots \times \R^{d_{1,r_1}} \times \cdots \times \R^{d_{q,0}} \times \cdots \times \R^{d_{q,r_q}}$ with coordinates $(x^{1,0},\dots,x^{1,r_1},\dots,x^{q,0},\dots,x^{q,r_q})$, where $x^{k,j} \in \R^{d_{k,j}}$, such that we have
\[
\call L = \sum_{k=1}^q \call L_k, \qquad H = (H^{1,0}, \dots, H^{1,r_1},\dots,H^{q,0},\dots,H^{q,r_q}),
\] 
satisfying the following, for each $k \in \{1,\dots,q\}$:
\begin{itemize}[leftmargin=2em]
\item $\call L_k$ is a standard diffusive operator acting at least on $x^{k,0}$ -- i.e., $\call L_k$ acts on $V_k \supseteq \{x \in \R^N :\, x^{-(k,0)} = 0 \}$;
\item $H^k = (H^{k,0},\dots,H^{k,r_k})$ is independent of $x^{-k}$, and verifies the conditions of \Cref{intro-thm1} --- i.e., $DH^k \in L^\infty(\R^N)$ and \eqref{eq_HhormC1intro} holds with $D_{x^{k,h}}H^{k,i}$ in place of $D_{x^h} H^i$. 
\end{itemize}
Let $u$ be a  bounded continuous (viscosity) solution to \eqref{TE}. Then both (1) and (2) from \Cref{intro-thm1} hold, with
\[
\call S(t) = \bigoplus_{k=1}^q t^{\bar\sigma_{\cL_k}^{-1}} \mathrm{diag}(t^iI_{d_{k,i}})_{0 \leq i \leq r_k}.
\]
\end{thm}

In this setting, we allow $\cL$ to be the sum of standard diffusive operators acting on subspaces that may also intersect with one another, as long as we can find chains of coordinates along which the diffusion acting on the first propagates to the others thanks to the structure of the drift highlighted in \Cref{intro-thm1} (where only one chain is present).
The proofs of both \Cref{intro-thm1,intro-thm1gen} are given in \Cref{sec:ex}, alongside with some explicit examples; in particular, we also show how heat conduction models through chains of oscillators, such as those considered e.g.\ in \cite{angeliki}, are included in our framework.

Eventually, as a consequence of \Cref{intro-thm1gen}, the exponential contractivity of the Fokker--Planck equation is deduced, provided that the existence of suitable Lyapunov functions is guaranteed.

\begin{thm}
\label{intro-thm2}
Suppose that the conditions of \Cref{intro-thm1gen} hold, and there exists a stationary Lyapunov function $\phi$ satisfying \eqref{eq_pshar} and \eqref{eq_condLyapdec}. If $m^i$ ($i=1,2$) are solutions to \rife{FPE} such that $\norm{\phi m^i(0)}_\TV < +\infty$ for $i=1,2$, we have 
\begin{equation*}
\norm{\phi(m^1(t) - m^2(t))}_\TV \leq Ke^{-\varpi t} \norm{\phi(m^1(0) - m^2(0))}_\TV \qquad \forall \, t >0,
\end{equation*}
where $K$ and $\varpi$ are positive constants independent of $m^i$ and $t$. 
\end{thm}

\subsection{Relation to previous literature} 
Hypocoercivity estimates, used to prove long-time convergence of Fokker--Planck or kinetic equations, were extensively developed since \cite{Desvillettes+V,H-Nier,herau,herau-nier}; a systematic presentation can be found in the celebrated notes of Villani \cite{Vil} (cf.~also \cite{Vil26}). In those results, the estimates are typically obtained in Hilbertian settings or, alternatively,  for relative entropy functionals. Along a similar direction, refined $L^2$-convergence methods were introduced in \cite{DMS}, later developed and  employed  for a broad class of kinetic models in \cite{Bao-Wang2,Bouin+al,Duan} and \cite{Bouin-frac} which includes the case of fractional diffusion in the kinetic Fokker--Planck equation. A key idea in the so-called DMS approach (after Dolbeault--Mouhot--Schmeiser~\cite{DMS})  is  played by splitting operators techniques (cf.\ also \cite{Gualdani+al}), exploiting the role of the diffusive and transport operators in suitably projected variables. In this respect, our approach shares some remote similarity in the idea of splitting diffusion and transport in the coupling argument (see the presentation in \Cref{model_LP}), though we have a completely different method that relies on pointwise estimates and the maximum principle, definitively far from Hilbertian methods.

New updates of $L^2$-hypocoercivity methods have appeared more recently, using intrinsic metric forms \cite{Armstrong+al} or refined  Poincar\'e inequalities \cite{Arnold2,Brigati-Stoltz}. Otherwise, still in the range of variational methods, the decay of entropy  functionals was also used  to prove the convergence to equilibrium of kinetic or hypoelliptic equations \cite{Arnold-Erb,Conforti+al}. 

Compared to the aforementioned  PDE results, our approach is very different; we do not use Hilbertian norms, nor we rely on Poincar\'e inequalities or properties of the stationary measures, often needed for the relative entropy decay method or for weighted  functional inequalities.
One advantage --- and motivation --- is to include fairly general nonlocal diffusions, hence our results apply to a broader class of operators. In addition, we go beyond pure  kinetic models including  H\"ormander conditions of any order. By contrast, our proofs do not provide in general with sharp rates of long-time convergence, for which variational methods or spectral gap estimates can be more accurate. 

However, in the special case of second order diffusions and linear drifts, we show how we can recover the sharp decay estimates of Arnold and Erb \cite{Arnold-Erb} with a simple Bernstein-type gradient estimate (see \Cref{model_LP}); to be specific, for this model case, \cite{Arnold-Erb} produces sharp decay of the relative entropy, whereas our method  gives the same decay for the total variation. The elementary proof that we give may hopefully be interesting for more general models.

Other different PDE methods are used for the convergence to equilibrium of kinetic equations; we refer, e.g., to \cite{Canizo+al,Yoldas} for an approach based on Harris' ergodicity theorem (relying on Harnack inequalities), whereas spectral gap analysis and abstract resolvent estimates are used in \cite{Gualdani+al,MM}, and \cite{Fathi+al}; see  this latter paper for a broader literature overview in this direction. Most of the above results seem however addressing only the kinetic case, rather than more general H\"ormander-type conditions.

There are results obtained with probabilistic methods which look closer to ours, standing on either gradient bounds or coupling techniques: hypocoercivity estimates were obtained via $\Gamma$-calculus and the Bakry--\'Emery approach \cite{Bau}, whereas some gradient estimate addressing hypoellipticity or hypocoercivity of kinetic equations were obtained via Bismut formula or coupling arguments \cite{Baudoin+al,Guillin-Wang}; coupling methods were used in \cite{Bao-Wang,Eberle+al} to get exponential decay of kinetic Fokker--Planck equations.

Beyond the interest of a purely PDE approach, it seems that our long-time oscillation estimates apply to a more general setting, both in terms of operators and in terms of H\"ormander structure conditions. In addition, the bounded-to-Lipschitz smoothing allows us to upgrade the decay to total variation norm, whereas probabilistic couplings naturally lead to the decay of Wasserstein-like  distances.

Another novelty of our approach is the combination of  short and long time scales, gathered in the same perspective; but if one is only interested in the local-in-time regularity, of course there are plenty of other results in the literature. Local smoothing estimates for hypoelliptic PDEs or stochastic models of chains of oscillators exist  through Malliavin calculus \cite{Del-Menozzi,Huang-menozzi, pigato} or through kernel estimates \cite{Difra-Pasc}; otherwise, the local hypoelliptic regularity theory was deeply developed recently even for rough coefficients, see e.g.\ the survey \cite{mouhot}. This is however a different issue which goes beyond --- and outside --- the spirit of our work. 

\subsection{Notations} We collect in the table below the main symbols we use throughout the paper, with their respective meaning.
\begin{small}
\begin{center}
\renewcommand{\arraystretch}{1} 
\begin{xltabular}{.96\textwidth}{@{} p{.25\textwidth}|X @{}}
$|\var|$ & The Euclidean norm. \\
$\norm\var$ & The operator norm, relative to the Euclidean norm; i.e., if $A$ is a matrix,
\begin{center}$\displaystyle\norm A \defeq \max_{|x|=1}|Ax|$.\vspace{7pt}\end{center} \\
$\norm\var_{L^\infty(\mu)}$, $\norm\var_\infty$ & The $L^\infty$-norm with respect to the measure $\mu$; if $\mu$ is the Lebesgue measure, we use $\norm\var_\infty$. If $M(t)$ is a matrix, $\norm{M}_\infty$ is the $L^\infty$-norm (in time) of the operator norm.\\
$\norm\var_{\mathrm{TV}}$ & The total variation norm; i.e., if $\mu$ is a signed measure on $\R^N$,
\begin{center}$\displaystyle \norm{\mu}_{\mathrm{TV}} \defeq \sup_{\norm{f}_\infty \leq 1} \int_{\R^N} f \,\di\mu$.\vspace{10pt}\end{center} \\
$\R_+$ & The positive real numbers. \\
$\R^{m\times n}$ & The space of $m \times n$ real matrices. \\
$\GLN$ & The space of invertible $N\times N$ real matrices, endowed with the operator norm. \\
$\Sym^N$, $\Sym_{\geq 0}^N$ & The spaces of symmetric and symmetric positive semidefinite $N \times N$ matrices, respectively, endowed with the operator norm. \\
$\Ran$, $\ker$ & The range (or image) and the kernel of a matrix, respectively. \\
$\mathrm{rk}$ & The rank of a matrix; i.e., $\mathrm{rk} = \mathrm{dim}\,\Ran$. \\
$\tr$, $\det$ & The trace and the determinant of a matrix, respectively. \\
$L^1(\call X,\mu;\call Y)$ & The space of $\call Y$-valued $\mu$-integrable functions on $\call X$; we omit $\mu$ if it is the Lebesgue measure, and we omit $\call Y = \R$. \\
$L^\infty(\call X,\mu; \call Y)$, $L^\infty_\loc(\call X,\mu; \call Y)$ & The space of $\call Y$-valued $\mu$-essentially bounded, resp.\ locally bounded, functions on $\call X$; we omit $\mu$ if it is the Lebesgue measure, and we omit $\call Y = \R$. Note that if 
$\phi \colon \call X \to \R$ is positive and measurable, $\phi L^\infty(\call X)$ is then the space of measurable functions $u$ such that $\phi^{-1}u$ is bounded almost everywhere with respect to the Lebesgue measure. 
\\
$C^k(\call X; \call Y)$, $C^k_b(\call X;\call Y)$ &   The space of $\call Y$-valued $k$ times continuously differentiable functions on $\call X$, with the subscript $b$ meaning that all derivatives from order $0$ up to $k$ are bounded. 
We omit $\call Y=\R$. \\
$\Lip(\R^N)$ &   The space of real Lipschitz continuous functions on $\R^N$; i.e.,
\begin{center} $\displaystyle f \in \Lip(\R^N) \quad\iff\quad \sup_{x\neq y} \frac{|f(x)-f(y)|}{|x-y|} < \infty$.\vspace{10pt}\end{center} \\
$\scr M(\call X)$, $\scr P(\call X)$, $\scr M_L(\call X)$ & The space of finite signed Borel, probability, and L\'evy measures on $\call X$, respectively; we recall that if $\call X$ is an Euclidean space, then L\'evy measures are nonnegative Borel measures on $\call X \setminus \{0\}$ satisfying
\begin{center} $\displaystyle\int_{\call X \setminus \{0\}} (1 \wedge |\var|^2)\,\di \nu < \infty$.\vspace{10pt}\end{center} \\
$\scr H^d$ & The $d$-dimensional Hausdorff measure. \\
$\var^*$ & The transpose (of a matrix), or the adjoint (of an operator). \\
$\var^\perp$ & The orthogonal complement of a vector subspace. \\
$\widehat{\var}$ & The normalization of a vector; i.e., $\widehat v \defeq v/|v|$. \\
$\mathrm{diag}$ & The map $\R^{m_0 \times n_0} \times \cdots \times \R^{m_k \times n_k} \to \R^{M \times N}$, with $\smash{M = \sum_{j=0}^k m_j}$ and $\smash{N = \sum_{j=0}^k n_j}$, such that, if $v = (v^0,\dots,v^k)$, then $\mathrm{diag}\,v$ is block diagonal with $(\mathrm{diag}\,v)_{jj} = v^j$. \\
$\otimes$ & The tensor product of vectors, or the Kronecker product of matrices; we also use $\xi^{\otimes2} \defeq \xi \otimes \xi$. \\
$\oplus$ & The direct sum. For vector spaces, $V \oplus W$ is $V+W$ when $V \cap W = \{0\}$; for linear operators $A \oplus B = \mathrm{diag}(A,B)$. We also use $\bigoplus_{j \in \{1,\dots,k\}} V_j \defeq V_1 \oplus \cdots \oplus V_k$. \\ 
$\vee$, $\wedge$ & The maximum and minimum of two real numbers, respectively; i.e.,
\begin{center} $a \vee b \defeq \max\{a,b\}, \qquad a \wedge b \defeq \min\{a,b\}$. \end{center}
Also, the logical conjunction and disjunction, respectively. \\
$\resmes$ & The measure restriction; i.e., if $m$ is a measure on $\call Y$, and $\call X \subseteq \call Y$ and $E \subseteq \call Y$ are measurable, $(m \resmes \call X)(E) = m(E \cap \call X)$. \\
$\var_\sharp$ & The pushforward (of a measure); i.e., if $m$ is a measure on $\call X$ and $f \colon \call X \to \call Z$, $f_\sharp m = m \circ f^{-1}$, where $f^{-1}$ is the pre-image via $f$; equivalently, 
\begin{center}$\displaystyle
\int_{\call Z} g\, \di f_\sharp m = \int_{\call X} g \circ f \,\di m,
$\end{center}
for any non-negative and measurable $g \colon \call Z \to \R$. \\
$\circ$ & The composition of maps; i.e. $(f \circ g)(x) = f(g(x))$. \\
$\lesssim$, $\lesssim_\lambda$ & Less than or equal to, up to a constant; i.e., $a \lesssim b$ means $a \leq Cb$ for some $C>0$. To stress that $C$ may depend on $\lambda$, we write $a \lesssim_\lambda b$. We refer to $C$ as an implied constant. Also used in the opposite direction $\gtrsim$, $\gtrsim_\lambda$. \\
$\asymp$, $\asymp_\lambda$ & Of the same order; i.e., $a \asymp b$ means $a \lesssim b \lesssim a$. To stress that the implied constants may depend on $\lambda$, we write $\asymp_\lambda$. \\
$\de_t$, $D_x$, $\Delta_x$ & The derivative with respect to $t$ (time derivative), the derivative (gradient) with respect to the (space) variable $x$, and the Laplacian with respect $x$, respectively. We omit $x$ if it is the only space variable (equivalently, if it includes all space variables). \\
$\1_{\var}$ & A generalized indicator function: whenever $\mathsf P$ is some property,
\begin{center}$\displaystyle
\1_{\mathsf P} \defeq 
\begin{dcases}
1 & \text{$\mathsf P$ is true} \\
0 & \text{$\mathsf P$ is false};
\end{dcases}$
\end{center}
for example, if $r$ is a real number, $\1_{r \geq 1} = \chi_{[1,+\infty)}(r)$, where $\chi_A$ denotes the usual indicator function of the set $A$. \\
$\call P_V$ & The orthogonal projection of $\R^N$ onto $V \subseteq \R^N$. We also write $\call P_V^\perp \defeq \call P_{V^\perp}$. \\
$[\var]_{\phi,\theta}$ & The oscillation semi-norm of order $\theta \in [0,1]$ weighted by $\phi$ --- cf.~\Cref{def_oscsmn}. \\
$u|^x_y$ & The difference of $u$ at $x$ and $y$; i.e., $u|^x_y \defeq u(x) - u(y)$. 
\end{xltabular}
\end{center}
\end{small}

 \section{A model case: linear drifts under H\"ormander's condition}\label{model_LP}

In this section, we explain the main ideas of our approach, starting with the model case of second-order operators and linear drifts satisfying H\"ormander's condition \cite{Ho}. Indeed, the simplest example of \rife{TE} is given by the second-order equation
\be\label{Hor-eq}
\partial_t u - \tr(AD^2 u) - Bx \cdot Du = 0 \qquad \text{in} \ \ (0,T) \times \R^{N},
\ee
where $A \in \Sym_{\geq 0}^N$ and $B \in \R^{N\times N}$ satisfy \emph{H\"ormander's condition of order $r$} -- that is, letting
\begin{equation} \label{eq_Vkspace}
V_k \defeq \sum_{j=0}^k \Ran(B^jA), \qquad k \in \N,
\end{equation}
we have
\be\label{Hor-r}
\min\{ k \in \N :\ V_k = \R^N \} = r.
\ee
As is well-known, the existence of $r \in \N$ such that \rife{Hor-r} holds is just one reformulation of the hypoellipticity condition introduced by H\"ormander; defining the vector fields $X_0 \defeq Bx \cdot D$ and $X = (X^1,\dots,X^N) \defeq AD$, then \rife{Hor-r} is equivalent to the fact that $\mathrm{Lie}(X^0,\dots,X^N)(x) = T_x \R^N \simeq \R^N$ for all $x \in \R^N$. This follows from the fact that the space generated by the iterated brackets of order $k$ between those vector fields can be identified with $\Ran(B^kA)$. It is also immediate to see the equivalence with 
\emph{Kalman's rank condition} for controllability; namely, defining the block matrices $C_k \defeq (A \mid BA \mid \cdots \mid B^k A)$ for $k \in \N$, \rife{Hor-r} can be rephrased by requiring that $\min\{k \in \N :\ \mathrm{rk}\,C_k = N\} = r$.  We refer, e.g., to \cite[Appendix]{LaPo} for the proof  of  the equivalence of several different formulations of the original H\"ormander's condition.  From that paper, we also recall the following canonical form to which any $(A,B)$ satisfying the hypoellipticity condition can be reduced.

\begin{rmk}[Canonical forms of the matrices] \label{rmk_cfm}
With respect to a fan orthonormal basis for the spaces $V_0, \dots, V_r$, the matrices $A$ and $B$ assume the following structure: there exist natural numbers $d_0 \geq d_1 \geq \cdots \geq d_r \geq 1$ such that
\[
A = \begin{pmatrix}
A_0	&0	& \dots 	& 0& 0 \\
0 	& 0	& \dots 	& 0& 0 \\
0		& 0 	& \dots	& 0& 0 \\
\vdots 	& \vdots	& \ddots	& \vdots	 & \vdots \\
0		& 0		& \dots	& 0	& 0
\end{pmatrix},
\qquad 
B = B_0 + E = \begin{pmatrix}
E_{00} 	& E_{01}	& \dots 	& E_{0,r-1}& E_{0r} \\
B_1 	& E_{11}	& \dots 	& E_{1,r-1}& E_{1r} \\
0		& B_2 	& \dots	& E_{2,r-1}& E_{2r} \\
\vdots 	& \vdots	& \ddots	& \vdots	 & \vdots \\
0		& 0		& \dots	& B_r	& E_{rr}
\end{pmatrix},
\]
where $A_0 \in \Sym^{d_0}$, $E \in \R^{N\times N}$ is the block matrix with entries $E_{ij}$ and $0$, while $B_0 \in \R^{N \times N}$ consists in the sub-diagonal $(B_1,\dots,B_r)$ and $0$ elsewhere, with $B_i \in \R^{d_i \times d_{i-1}}$ having $\mathrm{rk}\,B_i = d_{i}$ for all $i \in \{1,\dots,r\}$. More precisely, one has $d_0 = \mathrm{rk}\, A = \mathrm{dim}\,V_0$ and $d_i = \mathrm{dim}\,V_i - \mathrm{dim}\,V_{i-1}$ for all $i \in \{1,\dots,r\}$. Therefore, up to an orthonormal change of coordinates, we can suppose that $A$ and $B$ have the above form, and that $A_0 \in \Sym^{d_0}_{>0}$ -- i.e., there exists $\gamma > 0$ such that $A_0 \geq \gamma I_{d_0}$.
\end{rmk}

\begin{rmk}[Rank factorization for the canonical form] \label{rmk_frf}
As each $B_i$ in the block decomposition of $B_0$ has maximal rank and $d_i \leq d_{i-1}$, for each $i \in \{1,\dots,r\}$ there exists a $J_i \in \R^{d_{i-1} \times d_i}$ of rank $d_i$ such that $B_i J_i x^i \cdot x^i \geq \gamma |x^i|^2$ for some $\gamma > 0$, for all $x^i \in \R^{d_i}$. In fact, there exists a rank factorization $B_i = C_iF_i$ with $C_i \in \Sym^{d_i}_{>0}$ and $F_i \in \R^{d_{i} \times d_{i-1}}$ of maximal rank $d_i$, hence it suffices to take $J_i$ as the right inverse of $F_i$. 
\end{rmk}

Due to \Cref{rmk_cfm}, we can --- and will --- henceforth assume that equation \rife{Hor-eq} is written in the canonical form, where $x \in \R^{N}$ is split as  $x = (x^0,\dots,x^r)$ with $x^i \in \R^{d_i}$;  \Cref{rmk_frf} in turn tells that \rife{Hor-eq} can in fact be naturally reduced to a drift-diffusion equation having the structure required in \Cref{intro-thm1}.

\subsection{Gradient decay via Bernstein's method.}  \label{sec_BernM}
The following reformulation of the hypoellipticity property is particularly relevant to explain where our approach comes from: in fact (see \cite[Proposition~A.1]{LaPo}), H\"ormander's condition is also equivalent to the positivity property
\be\label{C(t)}
\call C(t)>0 \quad \forall\, t>0\,, \quad \text{where}\quad \call C(t):= \int_0^t e^{-sB}A e^{-sB^*}\di s\,.
\ee
This property has a significant role in the evolution of  $u$, solution of \rife{Hor-eq}. By applying the classical  Bernstein's method, we observe that $u\mapsto \call C(t) Du\cdot Du$ defines the correct metric to observe the decay of the Lipschitz semi-norm. We just need to observe that $\call C(t)= \int_0^t e^{-(t-\tau)B}A e^{-(t-\tau)B^*}d\tau$, hence $\call C$ is the unique solution to the ODE
\be\label{ode_C}
\begin{dcases}
\dot{\call C} + \call C B^* + B \call C - A = 0 \\
\call C(0) = 0.
\end{dcases}
\ee

\begin{prop}[Gradient decay via Bernstein] \label{eq_bernstest}
Let $u$ be a bounded and smooth solution to \eqref{Hor-eq} and let $\call C \colon \R \to \Sym^N$ be defined by \rife{C(t)}. If $A$ and $B$ satisfy H\"ormander's condition (equivalently, \eqref{C(t)} holds), then we have
\begin{equation} \label{eq_esthormbern}
\norm{\sqrt{2\call C(t)} Du(t,\var)}_\infty \leq \norm{u_0}_\infty \qquad \forall\, t >0\,.
\end{equation}
\end{prop}

\begin{proof}
 If $L(u):= \partial_t u-\tr(AD^2 u) - Bx \cdot Du$, then by direct computations we have
\[
L\bigl(\tfrac12 |u|^2 + \call C D u \cdot D u\bigr) = -2\tr(A D^2u\,  \call C D^2 u) + \bigl( \dot{\call C} + 2\call CB^* - A \bigr) Du \cdot Du \leq 0,
\]
where we use that $\call C \geq 0$ and the last quadratic form vanishes due to \eqref{ode_C}. Therefore, $w \defeq \tfrac12 |u|^2 + \call C D u \cdot D u$ is a positive subsolution of \eqref{Hor-eq} such that $w(0,\var) = \frac12|u|^2$. By maximum principle we get  $\|w(t)\|_\infty\leq \|w(0)\|_\infty$, hence \eqref{eq_esthormbern}.
\end{proof}

The estimate \rife{eq_esthormbern} not only contains the bounded-to-Lipschitz effect but also gives it in the right short-time scale; this can be better observed by using the asymptotic expansion of $\call C(t)$ as $t\dto 0$. Indeed, by \cite[Proposition~2.3 and Lemma~3.2]{LaPo}, there is a matrix $C_0 > 0$ such that 
\be\label{LP-asy}
\call C(t) = \call S(t) (C_0 + O(t)) \call S(t) \quad \hbox{as $t \dto 0$,}
\ee 
 where $\call S(t)$ is the diagonal matrix
\begin{equation} \label{eq_defscalmat}
\call S(t) = \mathrm{diag}\bigl( t^{k+\frac12}I_{d_k} \bigr)_{0 \leq k \leq r};
\end{equation}
therefore, \rife{eq_esthormbern} implies that there is $\tau>0$ depending only on $A$ and $B$ such that  
\be\label{scaled-R}
\norm{\call S(t) Du(t,\var)}_\infty \leq K {\norm{u_0}_\infty} \qquad \forall\, t \in (0,\tau),
\ee
where $K$ only depends on the least eigenvalue of $C_0$.

Note that \eqref{eq_defscalmat} gives the same scaling as that provided by \rife{eq_scalintro}. In fact, while we proved \rife{eq_esthormbern} using the explicit role of $\call C(t)$, resulting from Kalman's rank condition on $(A,B)$ and the equation \rife{ode_C}, estimate \rife{scaled-R} can hold in a more general setting (including nonlocal diffusions), by only relying on the time scaling of the diffusivity transmission rate. 
This means of course to replace Bernstein's method as well, which is what we explain in the next subsection.

Let us point out that, as estimate \rife{eq_esthormbern} holds for every $t>0$, it also provides long-time decays in special cases. For example, if $A,B$ satisfy the structure in \Cref{rmk_cfm} with $B = B_0$ (that is, $E=0$ and the canonical form of $B$ coincides with its so-called \emph{principal part} $B_0$),  then $\call C(t)$ diverges component-wise as $t\to \infty$, at  the same time scale, namely $C(t) \sim \mathrm{diag}\bigl( t^{k+\frac12}I_{d_k} \bigr)_{0 \leq k \leq r}$ --- cf.~\cite[eqs.~(2.17), (3.6)]{LaPo}. This implies the long-time polynomial decay of the gradient, namely estimate \rife{grad_t_decay}. 

Of course, when it comes at long-time convergence, some extra condition beyond the H\"ormander's structure is needed. As is well-known, in order to have ergodicity and the convergence of the associated Fokker--Planck equation towards a stationary state, one needs the matrix $B$ to be \emph{asymptotically stable}, meaning that every eigenvalue of $B$ has negative real part; under this condition, for instance, sharp decay of the entropy, for the corresponding Fokker--Planck equation, is proved in \cite{Arnold-Erb}. Correspondingly, we observe that a similar decay can be proved using again Bernstein's method.

One option for that is to replace $\call C(t)$ in \rife{eq_esthormbern} with $\call E(t):= e^{-Bt} e^{-B^*t}$, which is the solution of the pure transport problem
\[
\begin{dcases}
\dot{\call E} + \call E B^* + B \call E=0\\
 \call E(0)=I;
 \end{dcases}
\]
the same Bernstein's method used above would then apply to $\call E Du \cdot Du$ and provides with the estimate 
\[
\norm{e^{-tB^*} Du(t,\var)}_\infty^2 \lesssim \|Du(1,\var)\|_\infty^2 \lesssim \|u_0\|_\infty^2 \qquad \forall\, t\geq 1,
\]
which implies the exponential decay of $\norm{Du(t,\var)}_\infty$. 

Interestingly, the asymptotic stability of $B$ also allows us to obtain the decay in a {\it time-independent metric frame}. In fact, one knows from \cite[Lemma 4.3]{Arnold-Erb} that, since $B$ is stable, there exists a positive definite matrix $P$ such that
\be\label{arnerb}
BP+ PB^* \leq -2\mu P 
\ee
where $\mu = \min_{\lambda \in \mathrm{Spec\,B}} |\Re\lambda|$ (i.e.,  $-\mu$ is the maximum real part of the eigenvalues of $B$), or a strictly smaller value in case such minimum is attained at some defective eigenvalue.  
Then, as in the proof of \Cref{eq_bernstest}, we have
\[
L\bigl( P D u \cdot D u\bigr) \leq \bigl( BP+PB^*  \bigr) Du \cdot D u \leq -2\mu P D u \cdot D u;
\]
this immediately yields that $(t,x) \mapsto e^{2\mu t} P Du(t,x) \cdot Du(t,x)$ is a subsolution of \eqref{Hor-eq}, and by maximum principle we deduce
\[
P D u(t,\var) \cdot D u(t,\var) \leq e^{-2\mu (t-1)} P D u(1,\var) \cdot D u(1,\var)\qquad \forall\, t>1.
\]
Recalling the short-time bounded-to-Lipschitz smoothing and the invertibility of $P$,  we conclude that
\be\label{graddec}
\|Du(t,\var)\|_\infty \lesssim e^{-\mu t}\, \|u_0\|_\infty\qquad \forall t\geq 1,
\ee
which shows the exact exponential decay rate as in \cite{Arnold-Erb}. Note that the function $\phi(x) = Px\cdot x$ is in fact a stationary Lyapunov function for equation \rife{Hor-eq}, according to \cref{ass_lyap}.

Let us also recall that, by duality, the gradient decay of $u$ leads to a corresponding decay in total variation norm for the Fokker--Planck equation; indeed, by mean value theorem, \rife{graddec} implies
\[
\sup_{x,y} \frac{|u(t,x)-u(t,y)|}{1+ |x|+ |y|} \lesssim e^{-\mu t}\, \|u_0\|_\infty,
\]
hence one deduces the same decay in the total variation norm for the solution of the Fokker--Planck equation (in terms of the first moments bounds of the initial data);   cf.~\cite[Section~4]{PorAdv}, or \Cref{prop_decafp}.

\subsection{Oscillation decay via coupling arguments} 
To start with, let us consider the matrix $\call C(t)$ from \rife{C(t)}, once again. In view of estimate \rife{eq_esthormbern}, a natural guess is to rephrase the equation \rife{Hor-eq} in a time-scaled reference frame. Let us define
\[
v(t,x) \defeq u(t,\call R(t)x), \qquad \call R(t) \defeq \sqrt{\call C(t)}
\]
where $\sqrt{\call C(t)}$ denotes the positive semidefinite square root of $\call C(t)$; by direct computations, $v$ solves
\be\label{tildeq}
\de_t  v - \tr(\call  A(t) D^2 v) - \call B(t) x\cdot D v=0 \qquad \text{in} \ (0,T) \times \R^N ,
\ee 
where
\[
\call A(t) \defeq \call R(t)^{-1} A \call R(t)^{-1}, \qquad \call  B(t) \defeq \call R(t)^{-1} \bigl( \dot{\call R}(t) + B\call R(t) \bigr) .
\]
We first note that $\call B$ satisfies
\be\label{skew}
{\call  B} + {\call  B}^*  =  \call  A  \quad \text{on}\ (0,T);
\ee
indeed,  using that $\call C(t)= \call R(t)^2$ in the equation \rife{ode_C} we see that $\dot{\call R} \call R^{-1} + \call R^{-1} \dot{\call R}  + \call R B^* \call R^{-1} + \call R^{-1} B \call R =  \call A$, which implies \rife{skew}. Thus, at the expense of a time-dependent change of frame, the transport term can be made skew-symmetric -- hence, in particular, non-expansive -- on the kernel of the diffusion matrix.

We now apply to equation \rife{tildeq} a coupling argument which is similar, in spirit,  to the construction  in \cite{ForPor}, except that it is here performed at the natural embedded time scale. We define $\bb I(t) \defeq \Ran\,\call A(t) = \call R(t)^{-1} (\R^{d_0} \times \{0\})$ and we let $\Pi(t) \colon \R^N \to \bb I(t)$ be the orthogonal projection onto $\bb I(t)$; we also let $\Pi(t)^\perp$ be the orthogonal projection onto $\bb J(t) \defeq \bb I(t)^\perp$. Recall that since $\Pi^2 = \Pi$, we have $\Pi \dot \Pi + \dot \Pi \Pi = \dot \Pi$, hence $\dot\Pi(t)\bb I(t) \subset \bb J(t)$ and $\dot\Pi(t)\bb J(t) \subset \bb I(t)$. Heuristically, our coupling idea can be described as follows: for $x, x'\in \R^N$, we use the distance
\[
\rho= \rho(t,x,x') \defeq |\Pi(t)(x-x')|  + |\Pi(t)^\perp(x-x')|
\]
and the \emph{coupling matrix}
\[
\frk C(t,x,x') \defeq
\begin{pmatrix}
\call A(t) & \call A(t) - 2\gamma(t) \widehat{\Pi(t)\Delta x}{}^{\otimes2} \1_{\Delta x \notin \ker\call A(t)}\\
\call A(t)  -  2\gamma(t) \widehat{\Pi(t)\Delta x}{}^{\otimes2} \1_{\Delta x \notin \ker\call A(t)} &  \call A(t) \end{pmatrix},
\]
where $\Delta x= x-x'$, $\gamma(t)$ is the least positive eigenvalue of $\call A(t)$, and $\widehat{\Pi(t)\Delta x}$ denotes the normalized vector $|\Pi(t)\Delta x|^{-1} \Pi(t)\Delta x$.
Note that this choice guarantees that $\frk C(t,x,x')\geq 0$, and corresponds in fact to using a synchronous coupling whenever $\Delta x \in \ker\call A(t)$, and a coupling by reflection otherwise. 
Defining now
\[
\Psi (t,x, x') \defeq K_0\psi(\rho ), 
\]
for a suitable increasing Lipschitz function $\psi$, we wish to show by maximum principle that 
\be\label{goal_formal}
v(t,x) - v(t,x') \le  \Psi (t,x, x'),
\ee
and we know that in the PDE version of coupling arguments (see \cite{PoPr}), it will be enough to show that there is a choice of $\psi$ such that $\Psi$ satisfies
\be\label{soprasol_cou}
\partial_t \Psi - \tr\bigl(\frk C(t,x, x')D^2_{(x, x')}\Psi\bigr)- \bigl(\call Bx \cdot D_x \Psi +\call B x' \cdot D_{x'} \Psi  \bigr)\geq 0 \qquad \forall\, x\neq  x';
\ee
in fact, $v(t,x) - v(t, x')$ is a subsolution of the same equation, which will lead to establish the comparison \rife{goal_formal}.

We   sketch  a bit formally  the computations showing \rife{soprasol_cou}, in order to get heuristically the main ingredients we will use later.  First, using that $\rho$ is an even function of $x-x'$, we rewrite \rife{soprasol_cou} as
\be\label{verifica_cou}
\psi'(\rho) \bigl( \partial_t \rho  - \call B(t) \Delta x\cdot D_x\rho \bigr) -  4 \gamma(t) \Bigl( \psi'(\rho) D^2_{xx}\rho\, \widehat{\Pi(t) \Delta x}\cdot \widehat{\Pi(t) \Delta x} + \psi''(\rho)|\widehat{\Pi(t) \Delta x}\cdot D_x\rho|^2 \Bigr) \1_{\Delta x \not \in \ker \call A(t)} \geq 0,
\ee
where we have
\be\label{derrho}
\begin{gathered}
\de_t \rho  = {\dot\Pi   \Delta x \cdot \widehat{\Pi\Delta x}} - {\dot\Pi \Delta x  \cdot \widehat{\Pi^\perp\Delta x}} 
\,,\qquad 
D_x \rho  = \widehat{\Pi\Delta x} + \widehat{\Pi^\perp\Delta x},
\\
D^2_{xx} \rho = \frac{1}{|\Pi\Delta x|} \bigl( \Pi - \widehat{\Pi \Delta x}{}^{\otimes 2} \,\bigr) + \frac{1}{|\Pi^\perp\Delta x| } \bigl( \Pi^\perp - \widehat{\Pi^\perp \Delta x}{}^{\otimes 2} \,\bigr) \,.
\end{gathered}
\ee
Now we distinguish two different, mutually exclusive, regimes. If $\Pi(t)\Delta x=0$, then $\Delta x \in \ker\call A(t)$ by construction, and \rife{verifica_cou} reduces to
\begin{equation} \label{heurintro_case1}
\psi'(\rho) \bigl( \partial_t \rho  - \call B(t) \Delta x\cdot D_x\rho\bigr)\geq 0,
\end{equation}
however, on account of \rife{derrho} (and noting that in this regime $\Delta x= \Pi^\perp\Delta x$), we have
\[
\partial_t \rho  - \call B \Delta x\cdot D_x\rho= - \dot\Pi  \Pi^\perp\Delta x  \cdot \widehat{\Pi^\perp\Delta x} - \call B \Pi^\perp\Delta x\cdot \widehat{\Pi^\perp\Delta x} ,
\]
where both  terms above are zero, because $\dot\Pi \Pi^\perp\Delta x \in \Ran\call A = (\ker \call A)^\perp$ and $\call B$ is skew-symmetric on $\ker\call A$. Hence \rife{verifica_cou} trivially holds in this case with equality. 

In the alternative regime, we have $\Pi(t)\Delta x\neq 0$ and we can use the coupling by reflection. Here we first observe that 
\[
D^2_{xx}\rho\,\widehat{\Pi(t) \Delta x}\cdot \widehat{\Pi(t) \Delta x}=0,\qquad  \widehat{\Pi(t) \Delta x}\cdot D_x\rho=1,
\]
hence \rife{verifica_cou} reduces to
\[
\psi'(\rho) \bigl( \partial_t \rho  - \call B(t) \Delta x\cdot D_x\rho \bigr) - 4\psi''(\rho) \gamma(t)  \geq 0\,.
\]
Now, we need to note that the following time scalings hold, as $t \dto 0$, as a consequence of $A$ and $B$ satisfying H\"ormander's condition:\footnote{Recalling the asymptotic formula \eqref{LP-asy} and the canonical structure of \Cref{rmk_cfm}, the min-max theorem gives
\[
\gamma(t) =  \max_{\dim V = d_0} \min_{x \in V \setminus \{0\}} \frac{\call A(t)x \cdot x}{|x|^2} =  \max_{\dim V = d_0} \min_{y \in V \setminus \{0\}} \frac{ Ay \cdot y}{\call C(t)y \cdot y } \asymp \max_{\dim V = d_0} \min_{z \in V \setminus \{0\}} \frac{ \call S(t)^{-1}A\call S(t)^{-1} z \cdot z}{C_0 z \cdot z} = \frac\lambda t,
\]
where $\lambda$ is the least positive eigenvalue of $\sqrt{C_0}^{-1}A\sqrt{C_0}^{-1}$.
Similarly, letting $B_{\call S}(t) \defeq \call S(t)^{-1} B \call S(t)$, we estimate
\[ 
\norm{\call R(t)^{-1} B \call R(t)}^2 = \max_{x \neq 0} \frac{\call R(t)^{-1} B\call C(t) B^* \call R(t)^{-1} x \cdot x}{|x|^2} \asymp \max_{z \neq 0} \frac{B_{\call S}(t) (C_0 + O(t)) B_{\call S}(t)^* z \cdot z}{C_0 z \cdot z} 
\lesssim \frac1{t^2},
\] 
because by the canonical structure of $B$ we see that $tB_{\call S}(t)$ has bounded entries as $t\dto0$. As for the remaining part of $\call B$, namely $\call R^{-1}\dot{\call R}$, we start by observing that $(\call R \call S^{-1})^*(\call R \call S^{-1}) = C_0 + O(t)$, so taking the traces, the Frobenius norm of $\call R\call S^{-1}$ is bounded, which implies that each block $\call R_{ij}$ ($i,j \in \{0,\dots,r\}$, according to the canonical decomposition) has entries of order at most $t^{\max\{i,j\} + 1/2}$; that is, 
\[
\norm{\call R(t)_{ij}} \lesssim t^{\max\{i,j\} + \frac12} \quad \forall\, i,j = 0,\dots,r. 
\]
At this point, note that Rellich's theorem~\cite[Chapter~I, Section~1]{Rellich} ensures that there is an analytic orthonormal basis of eigenvectors of $\call C(t)$; that is, $\call C(t) = U(t)D(t)U(t)^*$, where $U(t)$ is orthonormal and $D(t)$ is the diagonal of the eigenvalues of $\call C(t)$, both are analytic in a neighborhood of $t=0$, and we know in addition that $D(t) = \mathrm{diag}\bigl(t^{2k+1}(I_{d_k}+O(t))\bigr){}_{0 \leq k \leq r}$, due to \eqref{LP-asy}. Therefore, we see that $t^{-1/2} \call R(t) = U(t)\bigl(t^{-1/2}\sqrt{D(t)}\,\bigr)U(t)^*$ is analytic; combining this piece of information with the previous bound, we deduce that $\call R(t)_{ij} = t^{\max\{i,j\} + 1/2} R(t)_{ij}$, with $R(t)$ analytic, so each entry of $\dot{\call R}(t)_{ij}$ is $O(t^{\max\{i,j\} - 1/2})$. We can then estimate as above  
\[
\norm{\call R(t)^{-1}\dot{\call R}(t)}^2 = \max_{x \neq 0} \frac{\call R(t)^{-1} \dot{\call R}(t)^2 \call R(t)^{-1} x \cdot x}{|x|^2} \asymp \max_{z \neq 0} \frac{|\dot{\call R}(t) \call S(t)^{-1} z|^2}{C_0 z \cdot z} \lesssim t^{2(\max\{i,j\} - j - 1)} \leq \frac1{t^2}.
\]
This concludes the proof that $\norm{\call B(t)} \lesssim \frac1t$. Finally, to estimate $\norm{\dot\Pi(t)}$, note that by an analogous argument starting from $(\call R^{-1}\call S)^*(\call R^{-1}\call S) = (C_0 + O(t))^{-1}$, one also deduces that $({\call R}(t)^{-1})_{ij} = t^{-\min\{i,j\} - 1/2} \tilde R(t)_{ij}$, with $\tilde R(t)$ analytic; this implies that $\call A(t) = \frac1t \mathsf{A}(t)$, with $\mathsf{A}(t)$ analytic, so $\Pi(t)$ is the projection onto the image of $\mathsf{A}(t)$. Recalling that we can write $\Pi = \mathsf{A}^\dagger \mathsf{A}$, where $\mathsf{A}^\dagger$ is the Moore-Penrose pseudoinverse of $\mathsf{A}$, standard computations lead to $\dot\Pi = (\mathsf{A}^\dagger \dot{\mathsf{A}} \Pi^\perp) + (\mathsf{A}^\dagger \dot{\mathsf{A}} \Pi^\perp)^*$, with $\mathsf{A}^\dagger \dot{\mathsf{A}} \Pi^\perp = \Pi \mathsf{A}^\dagger \dot{\mathsf{A}} \Pi^\perp$, whence $\norm{\dot\Pi} = \norm{\mathsf{A}^\dagger \dot{\mathsf{A}} \Pi^\perp} \leq \norm{\mathsf{A}^\dagger}\norm{\dot{\mathsf{A}}} \lesssim 1$, where the last bound is due to the regularity of $\mathsf A(t)$ and the first estimate on $\gamma(t)$  (which yields $\norm{\mathsf{A}^\dagger}$ bounded).
}
\be\label{scaling_dance}
\gamma(t) \asymp \frac1t,  \qquad  \partial_t \rho\geq -\rho \norm{\dot\Pi (t)} \gtrsim - \rho,
\qquad
\call B(t) \Delta x\cdot D_x\rho \leq \norm{\call B(t)} \rho \lesssim \frac \rho t ; 
\ee
hence, for small $t$, we have
\[
\psi'(\rho) \bigl( \partial_t \rho  - \call B(t) \Delta x\cdot D_x\rho \bigr) - 4\psi''(\rho) \gamma(t) \gtrsim - \frac1 t (\lambda \psi'' + \psi'\, \rho ),
\]
for some $\lambda > 0$, so, by choosing $\psi$ conveniently, the last term can be made positive, at least for bounded $\rho$. This shows that \rife{verifica_cou} holds in the second regime, as well. 

Hence $\Psi$ satisfies \rife{soprasol_cou}, which implies that $\call U (t,x,x') \defeq v(t,x) - v(t,x') - \Psi (t,x,x')$ cannot have a positive maximum in $ (0,T) \times \R^N \times \R^N$, at least for bounded $\rho$, e.g.\ for $\rho\leq 1$ from the above heuristic argument; on the other hand, by choosing $K_0= 2\|u\|_\infty/ \psi(1)$, the inequality \rife{goal_formal} would obviously hold for $\rho\geq 1$ as well. In addition, note that $\call U (0,x,  x')\leq 0$. 
Such a comparison would therefore lead to
\[
v(t,x) - v(t,x') \leq K_0\psi(\rho) \leq K_0 \, \|\psi'\|_\infty \rho \lesssim |x-\bar x|;
\]
by the definition of $v$, this means
\[
u(t,x)- u(t,x') \lesssim \bigl\lvert \call R(t)^{-1}(x-\bar x) \bigr\rvert, \qquad \call R(t) = \sqrt{\call C(t)},
\]
which is exactly another way to prove estimate \rife{eq_esthormbern}.

\medskip

To resume our heuristic presentation, the key point of the strategy is to apply a convenient (time-dependent) change of coordinates, represented  above by the matrix $\call R(t)$. Two main goals were realized after such a change of frame. Firstly, we used that the transport term was itself dissipative wherever the diffusion was not active. Secondly, the coercivity induced by the diffusion, exploited through reflection coupling, was effective {\it at any small time scale}.
     
Let us now make some comments about the above approach, as a guidance to what we will actually develop in the next sections.

\begin{enumerate}[leftmargin=*]

\item In terms of short-time scaling, we mostly used the asymptotic expansion of the covariance matrix $\call C(t)$ in order to verify \rife{scaling_dance}. To this purpose, simpler choices would still work for $\call R(t)$; for instance, due to \rife{LP-asy}, a natural guess would be to replace $\call R(t)$ with $\call S(t)\sqrt{C_0}$, where the diagonal matrix $\call S(t)$, given in \rife{eq_defscalmat}, would bring alone the essential scaling of each variable. The great advantage, to our purposes, is that   a similar time scaling is common to a much broader class of models, including nonlocal diffusions and chain of oscillators, without need to use the exact H\"ormander's condition \rife{C(t)}. 
In the abstract version of our result, we will only require the existence of a  time-scaling reference frame which is consistent with the coercivity of the diffusive part and for which the transport part is well-behaved at small scales; see our condition \rife{eq_condnormSH}.

\item  After the change of variables, the new transport term happened to be skew-symmetric on the kernel of the diffusive operator, due to the precise equation satisfied by $\call R(t)$. However, our argument only needed that the transport be dissipative whenever the diffusion vanishes --- that is, inequality is enough in \rife{skew}; hence we only need $\call R(t)$ to provide with a relaxed condition, such as $\call B + \call B^* \leq 0$ on the kernel of $\call A$. On one hand, this will give us a wider choice of $\call R(t)$; on another hand, a strict dissipativity condition may even help to better exploit the exponential decay in long time (cf. inequality \rife{arnerb}, which is available whenever $B$ is stable).

It is not a priori trivial that one can find a different change of frame that relaxes \eqref{skew}, while preserving the optimal scaling provided by $\call R(t)$.;
in particular, the na\"ive choice of $\call S(t)\sqrt{C_0}$ does not really work, as it satisfies exactly \eqref{skew} only if $\call L$ coincides with its principal part (i.e, $E=0$ in the canonical form of \Cref{rmk_cfm}), while in general it may make $\call B$ lose the necessary coercivity along some directions in $\ker\call A$.\footnote{Since $\mathsf R \defeq \call S\sqrt{C_0}$ is such that $\mathsf R(t) \mathsf R(t)^* = \call C(t) + \call S(t) O(t) \call S(t)$ due to \eqref{LP-asy}--\eqref{eq_defscalmat}, the equation~\eqref{ode_C} implies that the resulting diffusion $\mathsf A = \mathsf R^{-1} A (\mathsf R^*)^{-1}$ and drift $\mathsf B = \mathsf R^{-1}(\dot{\mathsf R} + B\mathsf R)$ satisfy $\mathsf B + \mathsf B^* = \mathsf A + O(1)$, where the $O(1)$ matrix may be not negative semidefinite on $\ker\mathsf A$. On the other hand, note that we do not have such an error term if $\call C(t) = \call S(t)C_0\call S(t)$, which happens when $B=B_0$ (see \cite[eq.~(2.17)]{LaPo}).} Nevertheless, a change of variables of the form $\call S(t)\call T$ is indeed desirable, as it has the advantage of resulting in the projection $\Pi \colon \R^N \to \Ran\,\call A$ being independent of $t$,\footnote{In this case, recalling the canonical form of $A$, we have $\call A = \call T^{-1}\call S(t^{-1}) A \call S(t^{-1}) (\call T^*)^{-1} = t^{-1} \call T^{-1}A(\call T^*)^{-1}$.} and there actually is a choice of $\call T$ such that $\call B + \call B^* < 0$ on $\ker\call A$ (see the proof of  \Cref{intro-thm1} for an explicit construction).

We also point out that, like the abstract counterpart of the scaling balance \eqref{scaling_dance} will be given by \eqref{eq_condnormSH}, a second requirement \eqref{eq_condcallN''} will formalize such a need for a strict dissipativity when the reflection coupling is turned off and we cannot get help from any diffusions. By the way, a strict sign in \rife{skew} will also allow us to compensate some error terms produced by unavoidable regularizations needed in the coupling argument.

\item A similar coupling estimate as above could be obtained using a stationary change of frame instead of $\call R(t)$. The output can be a gradient bound depending on the Lipschitz seminorm of the initial datum, and whenever a Lyapunov function $\phi$ exists, it can also provide the decay of weighted oscillation seminorms. 
This strategy was actually used in \cite{ForPor} for the purely kinetic equation, where the comparison \rife{goal_formal} is replaced by the stronger estimate 
\[
v(t,x) - v(t,x') \leq  K_0 e^{-\omega t} ( \phi(x)+ \phi(x') ) \psi(\rho).
\]
One of the improvements of  the present paper is to extend this idea, in order to get both the short-time smoothing (with a time-dependent change of variables) and the long-time decay (with a possibly stationary one).

\item To make the above coupling argument rigorous, a few more technicalities will be needed.  First of all, weak form of maximum principles are helpful, for which we opt to use the general framework of viscosity solutions. Then, as said, some regularization will be needed to handle the $\ell^1$ distance $\rho$ suggested above, as well as some classical perturbation argument to make our comparison effective up to $t=0$, and for $|x|$ and $|x'|$ up to infinity. At the expense of these technicalities, the main advantage is to reduce the idea explained before to a truly pointwise argument. Note that this approach avoids any coupling construction for the underlying SDEs trajectories, looking quite different from the probabilistic methods used, e.g., in \cite{Eberle+al}.

\end{enumerate}

\section{General setting, assumptions and statement of the main results} \label{sec:genset}

Here we introduce the general setting for our results and, in particular, we extend the framework of \Cref{model_LP} to include non-local diffusions and L\'evy jump operators, as well as a nonlinear transport vector field.

\subsection{Mixing local and nonlocal diffusions}

The background space is $\R^N$, we denote by $x=(x^0,\dots,x^{N-1})$ its generic variable. We consider on $(0,T) \times \R^N$ the equation
\begin{equation} \label{eq_nleq}
\de_t u = \call Lu + H\cdot Du, \qquad \call Lu \defeq \tr(AD^2u) + \call I[u],
\end{equation}
where we assume that $H \colon (0,T) \times \R^N \to \R^N$ is a continuous vector field (called the \emph{drift} term), and  
\begin{equation} \label{eq_condAt}
\text{$A \colon (0,T) \to \Sym_{\geq0}(N)$ is continuous and such that } A(t) \geq \bar A \in \Sym_{\geq0}(N) \ \  \forall\, t \in (0,T)\,.
\end{equation}  
When convenient, we will also use the more compact notation $\call L^H \defeq \call L + H \cdot D$.
Note that $\bar A$ is in general degenerate and $\ker\bar A$ is the subspace of $\R^N$ made of the variables that are not \lq\lq locally diffusive\rq\rq. The term $\call I[u]$ denotes the \emph{jump part} of the operator, and we assume that
\begin{equation} \label{eq_defcallI}
\call I[u](t,x) \defeq \int_{\R^N} \Bigl( u(t,x+z) - u(t,x) - Du(t,x) \cdot z \1_{|z| \leq 1} \Bigr) \, \nu(t)(\di z), \qquad \nu = \sum_{1 \leq j \leq \ell} \nu_j,
\end{equation}
with $\ell \in \N$ and each $\nu_j \colon (0,T) \to \scr M_L(\R^N)$, satisfying the following assumptions.

\begin{asss}[On the L\'{e}vy measure] \label{ass_nu}
For each $j \in \{1,\dots,\ell\}$, there are a $d_j$-dimensional linear subspace $V_j \subset \R^N$ and a radial Borel measure $\bar\nu_j$ such that $\bar\nu_j \resmes V_j  \in \scr M_L(V_j)$,
\begin{gather} \label{eq_boundsnu}
\nu_j \geq \bar\nu_j \resmes V_j,
\\ \label{eq_condnurad0}
\forall\, M \in \GLN, \quad M_\sharp \bar\nu_j \ll  \bar\nu_j, \quad \text{with} \quad M \mapsto \frac{\di M_\sharp \bar\nu_j}{\di \bar\nu_j} \in L^\infty_\loc(\GLN;L^\infty(\R^N,\bar\nu_j)),
\\ \label{eq_condnurad1}
\exists\, \underline c,\bar c >0, \ \alpha_j \in [0,2), \ \beta_j \in (\alpha_j - 1,\alpha_j] \ : \quad  \underline c \,r^{\alpha_j} \leq \int_{B_r \cap V_j} |\var|^2 \,\di\bar\nu_j \leq \bar c\,{r^{\beta_j}} \quad \forall\, r \leq 1.
\end{gather}
\end{asss}

Let us shortly comment about these assumptions.
\begin{itemize}[leftmargin=*]
\item The L\'evy measures $\nu_j$ may a priori be supported on the whole $\R^N$, but condition~\eqref{eq_boundsnu} tells that they are known to be coercive on a subspace $V_j$, with lower bound $\bar \nu_j$, which is a radial measure actually quantifying the diffusivity on $V_j$. Note that $\nu_j$ itself does not need to be radial, nor symmetric, indeed   \emph{any} $\nu_j$ that dominates $\bar\nu_j \resmes V_j$ is   admissible in our setting. For instance, given $\bar\nu_j$ as above, we can consider $\nu_j = h_j \bar\nu_j \resmes V_j$ for any bounded and uniformly positive $h_j \in \R^N \to \R_+$; on the other hand, in \Cref{ex_ntak} below, we also exhibit an anisotropic measure that fulfills Assumptions~\ref{ass_nu}, but whose density is not the product of a radial one and a uniformly positive one.

\item Condition~\rife{eq_condnurad0}, holding for both $M$ and $M^{-1}$,  is actually equivalent to
\[
\underline c_{M,j} \bar\nu_j \leq M_\sharp \bar\nu_j \leq \bar c_{M,j} \bar\nu_j,
\]
with $\underline c_{M,j}, \bar c_{M,j}>0$ locally bounded in $M$. In particular, when $M=\lambda I_N$ for some $\lambda > 0$, we see that this requires some amount of homogeneity to $\bar \nu_j$. When $\bar\nu_j$ has a density $\frk K(|\var|)$, it can be regarded as a \emph{relaxed} \emph{homogeneity} requirement on $\frk K$; indeed, in that case, it is also equivalent to $\frk K(\lambda r) \asymp_\lambda \frk K(r)$ locally uniformly in $\lambda$, i.e.
\begin{equation} \label{eq_wok'}
\inf_{(\lambda, r) \in [a,b] \times \R_+} \frac{\frk K(\lambda r)}{\frk K(r)} > 0 \qquad \forall\, a,b \in \R_+.
\end{equation}
For example, for continuous densities, such a property boils down to requiring that
\[
\lim_{r \dto 0} \frac{\frk K(\lambda r)}{\frk K(r)} > 0, \qquad \lim_{r \uto +\infty} \frac{\frk K(\lambda r)}{\frk K(r)} > 0.
\]

\item Condition \rife{eq_condnurad1} in turn ensures some scaling property of $\bar \nu_j$ on arbitrarily small balls, whose main example is obviously the fractional Laplacian kernel $|z|^{-d-\sigma}$ ($\sigma \in (0,2)$) when $V_j$ is $d$-dimensional, in which case $\alpha_j=\beta_j=2-\sigma$. However,  \eqref{eq_condnurad1} together with \eqref{eq_condnurad0} also implies a general control on every (possibly large) ball $B_R$:
\begin{equation} \label{eq_intnuR1}  
\forall\, R > 1 \quad \exists\, \underline C_{j,R}, \bar C_{j,R} > 0 \ : \quad \underline C_{j,R} \leq \int_{B_R \cap V_j} |\var|^2 \,\di\bar\nu_j \leq \bar C_{j,R}\,;
\end{equation}
this  will turn out to be useful in \Cref{sec_oscdecay}.
\end{itemize}

The basic idea of the above setting is that there are a number of diffusive operators, acting on possibly different subspaces $V_j\subset \R^N$, which may (or may not) include both second order  and jump operators. For the sake of a uniform notation, we will then set $V_0 \defeq \Ran  \bar A$ (given by \rife{eq_condAt}). It may well be the case that $V_0$ is just the sum of $V_j$'s, and the operator $\call L$ in \rife{eq_nleq} rewrites as
\be\label{somme_Levy}
\call L^Hu= \sum_{j=1}^\ell \cL_ju + H^{\ell+1} \cdot D_{x^{\ell+1}} u\,,\qquad x=(x^1, \dots, x^{\ell}, x^{\ell+1}),\ x^i\in \R^{d_i},
\ee
where $H^{\ell+1}\colon (0,T) \times \R^N \to \R^{d_{\ell+1}}$ is the drift part acting on the $x^{\ell+1}$ variable and $\cL_j\var = \tr(A_j D^2_{x^j}\var) + \call I_j[\var]$ are L\'evy operators acting on the variable $x^j$, with their own measure $\nu_j$; in this case, the matrix $A$ appearing in \rife{eq_nleq} would just be composed by the sum of the $\ell$ matrices $A_j$. 

The example \rife{somme_Levy} is certainly the simplest we have in mind, yet possibly different hypoelliptic conditions may link the set of diffusive variables $(x^1,\dots,x^\ell)$ with the deterministic variables $x^{\ell+1}$ (recall for instance the more general structure presented in \Cref{intro-thm1gen}). 
We will introduce later the needed assumptions to this purpose, and we will comment on them in the following section. 

The next examples serve as an illustration of possibly different operators that are included in the above general framework, including the interesting case that the subspaces $V_j$ have non trivial intersection.

\begin{exa}[A basic model] \label{ex_basicop} On account of  Assumption \ref{ass_nu}, the simplest model to keep in mind is when $A$ is degenerate and
\[
V_{j} = \{ x^{-j} = 0\}, \qquad \nu_{j} = \bar\nu_{j} \resmes V_{j}, \qquad \bar\nu_{j} = \frac{2-\sigma_j}{|\var|^{1+\sigma_j}} \scr H^1, \quad \sigma_j \in (0,2).
\] 
This  includes operators $\call L^H$ with $\call I$ being the sum of fractional Laplacians of different orders along different directions, such as, for example,
\[
\call L^H = - (-\Delta_x)^{\frac13} - (-\Delta_y)^{\frac34} + \Delta_z + y \cdot D_w,
\]
on $\R^4$ with typical coordinates $(x,y,z,w)$.  
\end{exa}

\begin{exa}[Intersecting diffusive directions] \label{ex_idd1}
The subspaces $V_j$ need not have trivial intersections. For example, on $\R^4$ with typical coordinates $(x,y,z,w)$, we can consider operators such as 
\[
\call L^H = - (-\Delta_{(x,y)})^{\frac12} - (-\Delta_{(y,z)})^{\frac34} + z \cdot D_w,
\]
which identifies the subspaces $V_1 = \R^2_{(x,y)}$ and $V_2 = \R^2_{(y,z)}$, according to Assumptions~\ref{ass_nu}.

As we are interested in smoothing estimates for  solutions to \eqref{eq_nleq}, it will be relevant, in the above case, to use the stronger diffusion $(-\Delta_{(y,z)})^{\frac34}$ in the $(y,z)$-direction and to consider the weaker $(-\Delta_{(x,y)})^{\frac12}$ ``just'' to deal with the smoothing effect in the $x$-direction --- while an improvement in regularity in the $w$-direction will be obtained  from the hypoelliptic  structure induced by the drift.  See \Cref{ex_idd2} for further discussion.
\end{exa}

Besides the model cases of the previous examples, and simple variations of them, Assumptions~\ref{ass_nu} allow to include more general L\'{e}vy measures, as 
we discuss below.

\begin{exa}[Two-scale kernels and non-homogeneous densities] \label{exa_interpfl}
Let $\nu = {\frk K}(|\var|) \scr H^d \resmes V$ with $d = \mathrm{dim}V$ and ${\frk K} \colon \R_+ \to \R_+$. Clearly, the kernel of a fractional Laplacian satisfies \eqref{eq_wok'}, and its associated measure satisfies \eqref{eq_condnurad1}; therefore, $\nu$ satisfies Assumptions~\ref{ass_nu} if $\frk K$ is continuous and
\begin{equation} \label{eq_exkasym}
\frk K \sim |\var|^{-d-\sigma} \quad \text{near $0$}, \qquad \frk K \sim |\var|^{-d-\sigma'} \quad \text{near $\infty$},
\end{equation}
with $\sigma,\sigma' \in (0,2)$. In fact, \eqref{eq_condnurad1} holds with the same exponents $\alpha = \beta = 2 - \sigma$ as in \Cref{ex_basicop}, since they only depend on the behavior of $\frk K$ on small balls about the origin. More generally, as condition \eqref{eq_condnurad1} allows $\alpha_j$ and $\beta_j$ to be different, we can also replace the homogeneous kernels in \eqref{eq_exkasym} with others like $\lvert\log(|\var|)\rvert^\beta |\var|^{-d-\sigma}$ with $\beta \in \R$ and $\sigma \in (0,2)$; note that these satisfy \eqref{eq_wok'} because
\[
\frac{\lvert\log(\lambda r)\rvert}{\lvert\log(r)\rvert} \to 1 \quad \text{as either $r \dto 0$ or $r \uto +\infty$}.
\]
\end{exa}

\begin{exa}[Non-trivial anisotropic kernels] \label{ex_ntak}
On $\R^N$, $N \geq 2$, with coordinates $(x^i)_{1 \leq i \leq N}$, consider a finite measure $\nu$ such that $\mathrm{supp}\, \nu = \tildee{\R^2} \defeq \{ x \in \R^N :\, x^i = 0 \ \forall\, i > 2\}$ and 
\[
\frac{\di \nu (x)}{\di x} = \frac{1}{(|x^1| + |x^2|^2)^\alpha} \1_{\tildee{\R^2} \cap B_1\!}(x) + \frac{1}{|x|^\alpha}\1_{\tildee{\R^2} \cap B_1\compl\!}(x), \qquad \alpha \in \Bigl( 2, \, \frac52 \Bigr).
\]
Clearly, $\di\nu(x) \geq |x|^{-2-\sigma} \1_{\tildee{\R^2}}(x)\, \di x$ for $\sigma = \alpha - 2 \in (0,\frac12)$, thus it satisfies Assumptions~\ref{ass_nu}, so we only need to show that $\nu$ is in fact L\'{e}vy. To this end, it suffices to note that by the change of variables $(x^1,x^2) = (y^2z,y)$ we have
\begin{equation*} \begin{split}
\int_{B_1 \cap \tildee{\R^2}} |\var|^2 \di \nu &= \int_0^1 y^{4-2\alpha} \int_{|z| < \frac{\sqrt{1-y^2}}{y^2}} \frac{y^2z^2 +1}{(|z|+1)^\alpha} \,\di z\di y \lesssim \int_0^1 y^{4-2\alpha} \,\di y \lesssim 1.
\end{split}
\end{equation*} 
\end{exa}

\begin{rmk}[Nonlocal \emph{as} local diffusions] \label{rmk_extloc}
 We  point out that, at the price of a less transparent presentation, we could have omitted   to explicitly write the local diffusion $\tr(AD^2\var)$, while   admitting instead that, for some $j \in \{1,\dots,\ell\}$,
\begin{equation} \label{eq_nltolocdef} \begin{split}
&\int_{\R^N} \Bigl( u(t,x+z) - u(t,x) - Du(t,x) \cdot z \1_{|z| \leq 1} \Bigr) \,\nu_j(t,\di z) \\
&\defeq \lim_{\sigma \uto 2}\, (2-\sigma) \int_{V_j} \Bigl( u(t,x+z) - u(t,x) - Du(t,x) \cdot z \1_{|z| \leq 1} \Bigr) \, \frac{\scr H^{d_j}(\di z)}{|z|^{d_j + \sigma}},
\end{split} \end{equation}
so that, when $u$ has bounded derivatives at $x$ up to the second order in the directions in $V_j$, the limit above equals $\frac{\omega_{d_j-1}}{2d_j} \tr(\call P_{V_j} D^2 u(t,x))$. 
in fact,   condition \eqref{eq_condnurad1} includes the scaling provided by local diffusions, when $\alpha_j = \beta_j = 0$.
\end{rmk}

\begin{rmk}[Integral averages of diffusions]
We could also consider further extensions to cover even more general diffusions. For instance, given $\mu \in \scr P([0,1])$, one can consider, in \eqref{eq_defcallI}, $\nu = \int_0^1 \nu_\lambda\, \di \mu(\lambda)$; note that this construction is reminiscent of that in \cite[Remark~5.7]{CDVARMA}. 
The proof of our main results would be similar to the case of $\mu$ discrete, with minor adaptations that we leave to the reader.
\end{rmk}

\subsection{Change of frame and time scaling}

As we explained in \Cref{model_LP}, in our approach the hypoellipticity is exploited by suitably changing the reference frame. Since time scaling is crucial to this purpose, we consider   a given time-dependent matrix  $\M \in C^1((0,T);\GLN)$ and   the new unknown function 
\begin{equation} \label{eq_defutrans}
u_\M(t,x) \defeq u(t,\M(t)x)\,.
\end{equation}
Whenever $u$ satisfies \rife{eq_nleq}, we readily obtain that $u_\M$ solves
\begin{equation} \label{eq_geneqcf}
(\de_t - \call L_\M^H)u_\M \defeq \de_t u_{\M} - \tr(A_\M D^2u_\M) - \call I_{\M}[u_{\M}] - H_{\M} \cdot Du_{\M} = 0,
\end{equation}
where $A_\M(t) \defeq \M(t)^{-1} A(t) (\M(t)^*)^{-1}$,
\begin{equation} \label{eq_deftildeH}
H_{\M}(t,x) \defeq \M(t)^{-1} \bigl( \dot{\M}(t)x + H(t,\M(t) x) \bigr),
\end{equation}
and
\begin{equation} \label{eq_deftildeI}
\call I_{\M}[v](t,x) \defeq \int_{\R^N} \Bigl( v(t,x + z) - v(t,x) - Dv(t,x) \cdot z \1_{|\M(t)z| \leq 1} \Bigr) \nu_{\M}(t)(\di z), \qquad \nu_\M(t) \defeq \M(t)^{-1}_{\,\sharp} \nu(t).
\end{equation}
Recall that $\nu= \sum_j \nu_j$, where each $\nu_j$ acts on  a subspace $V_j$ and satisfies Assumption \ref{ass_nu}. We now suppose  that $\M(t)^{-1} V_j$ is constant for each $t \in (0,T)$ and $j \in \{0,\dots,\ell\}$  (recall the $V_0 \defeq \Ran \bar A$) --- that is, there exist subspaces $V_j'$ such that 
\begin{equation} \label{eq_defSt}
\M(t) V_j' = V_j \qquad \forall\,t \in (0,T), \ j \in \{0,\dots,\ell\}.
\end{equation}
Then, by \eqref{eq_boundsnu} and \eqref{eq_condnurad0}, each $\nu_{\M,j}(t) \defeq \M(t)^{-1}_{\,\sharp} \nu_j(t)$ satisfies, for $t>0$,
\begin{equation} \label{eq_boundTnus}
s_{j}(t) \nu_{\M,j}(t) \geq \bar\nu_j \resmes V_j', \qquad 
\end{equation}
for some   positive constant $s_j(t)$,  which is optimally given by the $L^\infty$-norm of the Radon--Nikodym derivative of $\bar\nu_j \resmes V_j'$ with respect to $\nu_{\M,j}(t)$:
\begin{equation} \label{eq_stintdef}
s_j(t) \defeq \biggl\lVert \frac{\di(\bar\nu_j \resmes V_j')}{\di \nu_{\M,j}(t)} \biggr\rVert_{L^\infty(\nu_{\M,j}(t))}.
\end{equation}
Similarly, as for the local diffusion, from \eqref{eq_condAt} we have $s_0(t)  A_\M(t) \geq \call P_{V_0'}$, where $s_0(t)^{-1}$ is the least positive eigenvalue of $A_\M(t)$, and $\call P_{V_0'}$ denotes the orthogonal projection of $\R^N$ onto $V_0'   = \M(t)^{-1}(\Ran \bar A)$.


Requiring \rife{eq_defSt} --- i.e.\ that $\M(t)^{-1} V_j$ is constant in time --- may appear as a strong restriction on the matrix $\M(t)$; let us stress that the main purpose  of this condition is to allow us to work with constant-in-time projections, which is of course not needed a priori in our approach. In \Cref{model_LP}, for the purely second order case, we have shown how to handle a change of frame $\M(t)$ where this condition is not satisfied, namely by taking $\M(t)= \sqrt{\call C(t)}$ for the solution $\call C(t)$ of \rife{ode_C}. However, in the linear H\"ormander case, we also explained that the main scaling laws of $\M(t)$ simply appear through the diagonal matrix $\call S(t)$ in \rife{eq_defscalmat}, hence a simpler choice of $\M(t)$ can possibly be effective as well. In what follows, the reader should keep in mind that condition \rife{eq_defSt} is actually satisfied by our model choice, which is $\M(t):= \call S(t) \call T$, for some diagonal \emph{scaling matrix} $\call S(t)$ and some time-independent {\em twisting matrix} $\call T$.

 \begin{rmk}[Time scaling for homogeneous radial densities] \label{rmk_nu0hom} 
Suppose that $\bar\nu_j$ has a $\gamma$-homogeneous density (with respect to $\scr H^{d_j} \resmes V_j$) and $\M(t)$ acts on $V_j'$ as $t^{\delta} \call T$, for some $\call T \in \GLN$ such that $\call TV_j' = V_j$. Then,  the bound in \eqref{eq_boundTnus} holds with $s_j(t) \lesssim_{\call T} t^{-(\gamma + d_j)\delta}$. In particular,
\[
\begin{dcases}
\bar\nu_j = |\var|^{-d_j - \sigma} \scr H^{d_j} \resmes V_j \\
\M(t)|_{V_j'} = t^{\frac1\sigma} \call T
\end{dcases}
\quad\implies\quad
s_j(t) \leq c_{\call T} t,
\]
where we can take $\smash{c_{\call T} = \norm{\call T|_{V_j'}}^{-d_j-\sigma} \lvert\det(\call T|_{V_j'})\rvert}$. We also note that the above choice of $\delta=\sigma^{-1}$ is optimal when it comes to comparing diffusion and drift: indeed, $\smash{\M(t)^{-1} \dot{\M}(t)|_{V_j'} = t^{-1} I_N|_{V_j}}$, so $s_j(t) = O(t)$ is also needed to ensure that $\smash{s_j(t) \M(t)^{-1} \dot{\M}(t)|_{V_j'}}$ be bounded as $t\dto0$ --- cf.~the scalings \eqref{scaling_dance}, where $\gamma(t)^{-1}$ plays the role of $s_j(t)$.
\end{rmk}

\begin{rmk}[Boundedness of the time scaling]
In general, the function $s_j(t)$ defined in \eqref{eq_stintdef} is uniformly bounded on $(0,T)$ if $\M(t)$ is. Indeed, letting
\[
f(t,x) \defeq \frac{\di(\bar\nu_j \resmes V_j')}{\di \nu_{j}(t)} , \qquad g(M,x) \defeq \frac{\di M_\sharp(\bar\nu_j \resmes V_j')}{\di(\bar\nu_j \resmes V_j')},
\]
we have
\[
\frac{\di(\bar\nu_j \resmes V_j')}{\di \nu_{\M,j}(t)} =  \frac{\di(\bar\nu_j \resmes V_j')}{\di\M(t)^{-1}_{\,\sharp}(\bar\nu_j \resmes V_j')} \frac{\di\M(t)^{-1}_{\,\sharp}(\bar\nu_j \resmes V_j')}{\di \nu_{\M,j}(t)} = g(\M(t),\M(t)\var)f(t,\M(t)\var) \quad\ \text{$\nu_{\M,j}(t)$-a.e.},
\]
so \eqref{eq_boundsnu} and \eqref{eq_condnurad0} yield
\[ \begin{split}
s_j(t)
\leq \norm{g(\M(t),\var)}_{L^\infty(\bar\nu_j)} \norm{f(t,\var)}_{L^\infty(\nu_j(t))} \leq \norm{g(\M(t),\var)}_{L^\infty(\bar\nu_j)}.
\end{split}
\] 
\end{rmk}

\subsection{Bounded-to-Lipschitz smoothing}
 
One of our main results is the short-time estimate of \Cref{thm_mainsmnl} below, which provides with the bounded-to-Lipschitz decay for solutions of \rife{eq_nleq}. In order to have it in a possibly general setting, the statement is proved in the class of viscosity solutions, which is the natural framework of weak solutions for the maximum principle to hold --- see also \Cref{rmk_onuLip} about the stability of vanishing viscosity limits.
We refer to \cite{CIL} for a classical background of viscosity solutions, and to \cite{BarImbHP} for L\'evy-type operators (cf.~\cite[Section 2.1]{PorAdv} for a specific context similar as \rife{eq_nleq}). 

In the following, we will make use of possibly singular norms, and suitable regularizations of them. Slightly more generally, for $\epsilon \in [0,1)$, we will represent them by functions $\scr N_\epsilon \colon \R^N \to \R_{\geq0}$, which are not actually required to be norms, but still satisfy the following:
\be\label{ass_Neps}
\text{$\scr N_\epsilon \in C^2(\R^N)$ \ \ for $\epsilon \in (0,1)$ \quad with \quad $|D{\scr N}_\epsilon| \leq   C_{\scr N}$} 
\ee
\begin{equation} \label{eq_condcallN'}
c_{\scr N}|z| \leq  \scr N_0(z) \leq C_{\scr N}|z|, \qquad   {\scr N}_0(z) - C_{\scr N}\epsilon \leq{ \scr N}_\epsilon(z) \leq {\scr N}_0(z) \qquad \forall\, z \in \R^N
\end{equation}
for some positive constants $c_{\scr N}$ and $C_{\scr N}$ independent of $\epsilon$.
Hence, morally, $\scr N_0$ is the singular norm, and $\epsilon > 0$ is the regularization parameter; note that, for example, $\scr N_\epsilon(z) = \sqrt{|z|^2 + \epsilon} - \epsilon$ satisfies both \eqref{ass_Neps} and \eqref{eq_condcallN'}.

We also formalize the following definition (cf.~\eqref{eq_defoscintro}), already used in \cite{ForPor}.

\begin{defn}[Oscillation semi-norms] \label{def_oscsmn}
Given $v \colon \R^N \to \R$, $\phi \colon \R^N \to (1,+\infty)$, and $\theta \in [0,1]$, we define the \emph{oscillation $\theta$-semi-norm of $v$ weighted by $\phi$}
\be\label{osci-seminorm}
[v]_{\phi,\theta} \defeq \sup_{x,y \in \R^N} \frac{|v(x) - v(y)|}{\bigl(\phi(x)+\phi(y)\bigr)(|x-y|^\theta \wedge 1)}.
\ee
\end{defn}

\begin{prop}[Bounded-to-Lipschitz smoothing] \label{thm_mainsmnl}
Let $u \in C^0([0,T] \times \R^N) \cap L^\infty((0,T);\Lip(\R^N))$\footnote{With a slight abuse of notation, $u \in L^\infty((0,T);\Lip(\R^N))$ means that $[u(t,\var)]_{0,1}$ is finite and (essentially) bounded in $t$.} be a viscosity solution to \eqref{eq_nleq}, where $H(t,\var)$ is Lipschitz continuous, $A$ satisfies \eqref{eq_condAt}, and $\call I$ is given by \eqref{eq_defcallI} with $\nu_j$ satisfying Assumptions~\ref{ass_nu}.

Suppose there exist $\M \in C^0((0,T];\GLN) \cap C^1((0,T);\GLN)$ satisfying \eqref{eq_defSt}--\rife{eq_boundTnus} and $\M(0^+) = 0$, a subset $I \subseteq \{0,\dots,\ell\}$, and a collection of (non-trivial) mutually orthogonal subspaces $\{W_j\}_{j \in I}$ such that the following hold:
\begin{itemize}[left=\parindent]
\item we have the decomposition
\begin{equation} \label{eq_sbspdec}
 \call V \defeq \sum_{0 \leq j \leq \ell} V_j' = \bigoplus_{j \in I} W_{j}, \qquad \text{with} \qquad W_{j} \subseteq V_{j}' \quad \forall j \in I;
\end{equation}
\item given $ H_\M$ defined as in \rife{eq_deftildeH}, there is $K_H > 0$ such that
\begin{equation} \label{eq_condnormSH}
\max_{j \in I} s_j(t) \, \| DH_\M(t, \cdot)\|_\infty 
 \leq K_H    \qquad \forall \, t \in (0,T);
\end{equation}
\item there are $\gamma, \gamma' > 0$ (independent of $\epsilon$) such that
\begin{equation} \label{eq_condcallN''}
DH_\M(t,x) z \cdot \call P_{\call V}^\perp D\scr N_\epsilon(z) \leq - \gamma  \bigl( \norm{DH_\M(t, \cdot)}_\infty \vee 1\bigr) \scr N_0(z)  + \gamma'\epsilon
\quad\ \forall\, z \perp \call V, \ (t,x) \in (0,T) \times \R^N
\end{equation}
for some $\scr N_\epsilon$ satisfying \rife{ass_Neps}--\rife{eq_condcallN'}. 
\end{itemize}
Finally, let us be given a function $\phi \in C^2((0,T) \times \R^N;(1,+\infty))$ such that
\begin{gather} \label{lyapt1}
(\de_t - \call L^H)\phi \ \ \text{is uniformly bounded below on} \ (0,T) \times \R^N, 
\\ \label{eq_pshar}
\forall\,R>0 \quad \exists\, C_R \geq 1 \ : \quad C_R^{-1} \sup_{x+B_R} |D\phi(t,\var)| \leq \phi(t,x) \leq C_R \inf_{x + B_R} \phi(t,\var) \quad\ \forall\, (t,x) \in (0,T) \times \R^N,
\end{gather}
and suppose there exists a function $\Phi \in C^2((0,T) \times \R^N;(1,+\infty))$ such that
\be\label{ultima_Phi}
(\de_t - \call L^H)\Phi \ \ \text{is uniformly bounded below on} \ (0,T) \times \R^N, \quad \lim_{|x|\to+\infty} \sup_{t \in (0,T)} \frac{|u(t,x)| \vee 1}{\Phi(t,x)}=0.
\ee
Then, we have
\begin{equation} \label{eq_gthm_usmoothing}
|u(t,x) - u(t,y)| \leq C_T \sup_{t\in(0,T)} [u(t,\var)]_{\phi(t,\var),0} \bigl(\phi(t,x)+\phi(t,y)\bigr) \bigl\lvert \M(t)^{-1} (x-y) \bigr\rvert \quad\ \forall\,t \in (0,T], \ x,y \in \R^N,
\end{equation}
where $C_T$  depends on $T$, $d$, $N$, $|I|$, $\|\M\|_\infty$, $c_{\scr N}$, $C_{\scr N}$, $K_H$, $\gamma$, $\gamma'$, $\underline c$, $\bar c$, and $\phi$.
\end{prop}

We  invite the reader to consider first the model case that $\phi=1$ in the above statement: this is the classical bounded-to-Lipschitz result; more general functions $\phi$, having the role of Lyapunov functions, may be required later to match this result with the long-time estimates. By contrast, the function $\Phi$ plays no role other than fixing a class of solutions where the proof of the estimate can be justified (indeed, the constant $C_T$ in \rife{eq_gthm_usmoothing} does not depend on $\Phi$); in addition, such a $\Phi$ can be readily constructed in most cases --- cf. \cref{rmk_onuLip}.

The proof of \Cref{thm_mainsmnl} will be given in \Cref{proofs}; the ideas are essentially those explained in \Cref{model_LP}, namely we apply a coupling argument, via maximum principle, using a reflection coupling on the diffusive spaces $W_{j}$, while synchronous coupling is applied on $\call V^\perp$. 
A few clarifications are given below on the role of $\call V$, $W_j$ and the consequent choice of norm $\scr N$ that we may use. 
The main explicit applications, detailed in Section \ref{sec:ex}, will then lead to the statements of the Introduction.

\begin{rmk}[On the abstract assumptions]\label{spiega_ipotesi}
Here is a quick orientation to understand the above statement.

A priori, there are a certain number of  diffusive operators in   \rife{eq_nleq} which act, respectively, on subspaces $V_j$; after the change of frame given by $\M(t)$, those diffusions act on $V_j'= \call R(t)^{-1}V_j$. Since there may be some overlapping between the subspaces $V_j$ (cf. \Cref{ex_idd1}), what really matters is to decompose $ \call V$ into mutually orthogonal $W_j$ which  identify all the (mutually independent) diffusive variables; this explains \rife{eq_sbspdec}.  

Assumption \rife{eq_condnormSH}  is a condition on the time-scaling induced by the moving frame given by $\call R(t)$; it ensures that the new drift $H_{\M}$ rescales in accordance with the natural scaling of the diffusive parts, given by the $s_j(t)$ from \rife{eq_boundTnus}. The reader should keep in mind what happens in the classical H\"ormander case, explained in \Cref{model_LP} --- cf.~in particular \eqref{scaling_dance}. 

Finally, assumption \rife{eq_condcallN''} is a condition on the twisting effect induced by the moving frame given by $\call R(t)$; it ensures that the new drift $H_{\M}$ is {strictly} dissipative on $\call V^\perp$ --- i.e., whenever all diffusions vanish. 
\end{rmk}

Let us stress that the choice of $W_j$ has indeed some degree of freedom, depending on the true interaction between diffusion and drift as well as on the optimal time scaling needed to observe the diffusivity on each variable. The following examples should serve as an illustration of this fact.

\begin{exa}[Splitting and norms in case of intersecting diffusive directions] \label{ex_idd2}
If we consider  the operator $\call L^H$ of \Cref{ex_idd1}, namely
$$
\call L^H = - (-\Delta_{(x,y)})^{\frac12} - (-\Delta_{(y,z)})^{\frac34} + z \cdot D_w\,,
$$
we have $V_1 = \R^2_{(x,y)}$ and $V_2 = \R^2_{(y,z)}$ which have non trivial intersection. Here, since it is convenient to associate the second diffusion to the variable $y$, we will  choose $ W_1= \R_x$ and then 
consider the splitting of $\call V = \R^3_{(x,y,z)}$ given by
\[
\begin{tikzcd}[column sep=-3pt]
\call V &=&\smash{V_{1}} &+ &V_2 
\\[-3pt]
&=&\R_x \arrow[u,symbol={\subseteq}] &\oplus &\R^2_{(y,z)} \arrow[u,symbol={=}] 
\end{tikzcd}
\]
and consequently $\R^4 = \call V \oplus \R_w$ with norm of the form $|x| + \sqrt{y^2+z^2} + |w|$. Then we will   show (cf.\ the proof of \Cref{intro-thm1gen}) that a suitable choice of  $\M(t)$ exists which fixes $V_1, V_2$ and satisfies \rife{eq_condnormSH}--\rife{eq_condcallN''} (for this norm). 

Similar considerations apply  when there are $i \in \{0,\dots,\ell\}$ and $J \subset \{0,\dots,\ell\} \setminus \{i\}$ such that $V_{i} \subseteq \sum_{j \in J} V_{j}$; for example,
\[
\call L^H = - (-\Delta_{(x,y)})^{\frac12} - (-\Delta_{y})^{\frac\sigma2} - (-\Delta_{(y,z)})^{\frac34} + z \cdot D_w.
\]
In this case, we may want to use the same decomposition as before, if $\sigma \leq \frac32$; such a choice ignores the second diffusion and, according to \Cref{intro-thm1gen}, yields the scaling $\call S(t) = (t^{2}, t^{\frac23}, t^{\frac23}, t^{\frac53})$. On the other hand, it $\sigma > \frac32$, the diffusion $(-\Delta_y)^{\frac\sigma2}$ is stronger than $(-\Delta_{(y,z)})^{\frac34}$ in the $y$-direction, thus suggesting to use the finer splitting $\R^4 = \R_x \oplus \R_y \oplus \R_z \oplus \R_w$ --- and thus a norm like $|x|+|y|+|z|+|w|$ --- in order to obtain better estimates; indeed, in this case, we have $\call S(t) = (t^{2}, t^{\frac1\sigma}, t^{\frac23}, t^{\frac53})$.
\end{exa}

\begin{exa}[Interplay between diffusion and drift]
If one has the stronger condition that \eqref{eq_condcallN''} holds for all $z$ orthogonal to $\sum_{j \in J} V_j'$, for some subset $J \subset \{0,\dots,\ell\}$, then one can take $I \subseteq J$ and a decomposition $\{W_j\}_{j \in I}$ satisfying \eqref{eq_sbspdec} with $\sum_{j \in J} V_j'$ in place of $\call V$. This reflects the fact that we need to perform a coupling by reflection only to deal with those directions for which the drift does not provide (up to some suitable twisting) the needed coercivity; on the other hand, using the coercivity of the drift instead of some diffusions may or may not be advantageous to have a better short-time scaling, depending on the strength of the diffusions in play.

For example, let us consider, on $\R^3$ with coordinates $(x,y,z)$,
\[
\call L^H = \Delta_x - (-\Delta_y)^{\frac\sigma2} + x \cdot D_y + y \cdot D_z,
\]
to which we can apply \Cref{intro-thm1gen} with $\phi = 1$. Given the structure of the drift, we may ignore the second diffusion and propagate regularity by using only the first Laplacian, thus obtaining the standard H\"ormander scaling $\call S(t) = (t^{\frac12}, t^{\frac32}, t^{\frac52})$. On the other hand, we may also consider that the regularity in the $z$-variable is obtained through the fractional Laplacian acting on the $y$-variable, thus obtaining $\call S(t) = (t^{\frac12}, t^{\frac1\sigma }, t^{\frac1\sigma +1})$, which is a better short-time scaling than the previous one if $\sigma > \frac23$.
\end{exa}

\begin{rmk}[Time-dependent decomposition in subspaces] \label{rmk_genSigma'}
Furthermore, it will be clear from the proof that we may allow some dependence of $I$ on time. 
For example, if $V_j = V$ for all $j \in \{0,\dots,\ell\}$, then we may choose $I = I(t) = \{j(t)\}$ for any $j \colon (0,T) \to \{0,\dots,\ell\}$. 
Such a flexibility allows to interpret the decomposition $\{W_j\}_{j \in I}$ by favoring certain measures $\nu_j$, in order to optimize the time rate eventually provided by $s_j(t)$, distinguishing for instance the short and long time decay. See \Cref{sharp-Kolm} for an application of this observation.
\end{rmk}

Let us also point out that the simpler scenario of $\sum_{0 \leq j \leq \ell} V_j = \R^N$ is included as well. In this case, condition \eqref{eq_condcallN''} is void, as we do not need the drift to provide coercivity to compensate any missing diffusive direction.
For instance, consider the case that $\M(t) = \call S(t)\call T$ with each $V_j$ invariant under $\call S \in C^1((0,T);\Sym(N))$ and $\call T \in \GLN$ such that an admissible decomposition $\{W_j\}_{j \in I}$ does exist; here, the twist $\call T$ plays no role other than ensuring that we can take the mutually orthogonal $W_j$'s, and the choice of $\call S$ is solely constrained by \eqref{eq_condnormSH}.

\begin{rmk}[Solutions need not be Lipschitz] \label{rmk_onuLip}  
In \Cref{thm_mainsmnl}, we assume that $u$ is a Lipschitz continuous solution such that $u(t,x)=o(\Phi(t,x))$ as $|x|\to \infty$ (uniformly in $t\in(0,T)$), for some supersolution $\Phi$. In fact, if we assume that $\frac{H(t,x)}{1+ |x|} \in L^\infty((0,T)\times\R^N)$, and the upper control $\nu_j \lesssim |\var|^{-d_j-\sigma_j} \scr H^{d_j} \resmes V_j$ for some $\sigma_j \in (0,2)$, one can explicitly exhibit such a $\Phi$, with some power growth at infinity, following e.g.\ the proofs of \cite[Lemma~3.1 and Proposition~6.1]{PorAdv}. Even for general L\'evy measures, without any prescribed decay at infinity, one can construct an unbounded supersolution, at least if $\nu$ is independent of $t$, following \cite[Lemma~4.9]{Espen+al}. This would allow, for instance, to use the above result for bounded Lipschitz solutions; however, note that the estimate \eqref{eq_gthm_usmoothing} is independent of the Lipschitz constant of $u$, as well as of the function $\Phi$:  
indeed, it only depends on the oscillation semi-norm of $u$, weighted by $\phi$.
As a consequence, with usual stability arguments, we can deduce the same result for general viscosity solutions, which a priori can be unbounded and not necessarily Lipschitz continuous. The strategy is standard: by considering first the case of bounded and Lipschitz initial data, one proves that solutions have stable  Lipschitz bounds at positive times, which guarantees the local compactness to get at the general case. 
This can be seen by following for instance the arguments of \cite[Proposition~6.1 and Lemma~6.2]{PorAdv} and \cite[Proposition~6.1]{ForPor}. On the other hand, by comparison,\footnote{Indeed, \eqref{lyapt1} implies that $(\de_t-\call L^H)\phi \geq -k_0$ for some $k_0>0$, so one readily checks that $(t,x) \mapsto e^{k_0t}\norm{\phi^{-1}u(0,\var)}_\infty\phi(t,x)$ is a supersolution of \eqref{eq_nleq}.}
\begin{equation} \label{eq_ucompC3}
|u| \lesssim \norm{\phi^{-1}u(0,\var)}\phi \quad \text{on} \ (0,1) \times \R^N,
\end{equation}
yielding $\sup_{t \in (0,T)} [u(t,\var)]_{\phi(t,\var),0} \lesssim \norm{\phi^{-1} u(0,\var)}_\infty$. Hence, our estimate can take the simplified form 
\begin{equation} \label{eq_gthm_usmoothing'}
|u(t,x) - u(t,y)| \leq C_T \norm{\phi^{-1} u(0,\var)}_\infty \bigl(\phi(t,x)+\phi(t,y)\bigr) \bigl\lvert \M(t)^{-1} (x-y) \bigr\rvert \qquad \forall\,t \in (0,T], \ x,y \in \R^N
\end{equation}
\emph{for any continuous viscosity solution $u$}  such that $u\in  \phi L^\infty((0,T)\times\R^N)$.
\end{rmk}

\subsection{Exponential decay of oscillation semi-norms} \label{sec_oscdecay}

Given the same setting used above, we now state the long-time decay estimate for the solution of \rife{eq_nleq}. To this purpose, we assume the existence of a Lyapunov function according to the following definition.

\begin{defn}[Stationary Lyapunov function] \label{ass_lyap}
We say that $\phi\in C^2(\R^N;(1,+\infty))$ is a \emph{stationary Lyapunov function} if
\begin{equation}  \label{eq_condLyapdec}
\lim_{|x| \to +\infty} \phi(x)=+\infty,
\qquad\quad 
\varpi_0 \defeq \liminf_{|x| \to +\infty} \inf_{t\in(0,T)} \frac{(-\call L^H\phi)(t,x)}{\phi(x)} > 0.
\end{equation}
\end{defn}

Additionally, we will also be requiring $\phi$ to satisfy the more technical condition~\eqref{eq_pshar} (where $\phi(t,x)=\phi(x)$).
\vskip0.5em
Next is our second main result, stating the global-in-time decay of $[u(t,\var)]_{\phi,\theta}$; the setting is the same as in \Cref{thm_mainsmnl}, except that we now use a stationary change of frame.

\begin{prop}[Decay of oscillation semi-norms] \label{prop_oscdec}
Let $u \in C^0([0,T] \times \R^N) \cap L^\infty((0,T);\Lip(\R^N))$ be a viscosity solution to \eqref{eq_nleq}, where $H(t,\var)$ is Lipschitz continuous, $A$ satisfies \eqref{eq_condAt}, and $\call I$ is given by \eqref{eq_defcallI} with $\nu_j$ satisfying Assumptions~\ref{ass_nu}.

Assume that $[u(0,\var)]_{\phi,\theta}$ is finite for some $\theta \in (0,1]$, and that there exists a stationary Lyapunov function $\phi$  satisfying \eqref{eq_pshar}. 
Suppose that there exists $\M \in \GLN$ and a collection of orthogonal subspaces $\{W_j\}_{j \in I}$, satisfying \rife{eq_sbspdec}, such that both \eqref{eq_condnormSH} and \eqref{eq_condcallN''} hold (for some $\scr N_\epsilon$ satisfying \rife{ass_Neps}--\rife{eq_condcallN'}).\footnote{Note that this means that we are requiring the assumptions of \Cref{thm_mainsmnl} to hold for some $\M$ \emph{independent of $t$}.}

Then we have
\begin{equation} \label{eq_gthm_uosc}
[u(t,\var)]_{\phi,\theta} \leq C_\theta e^{-\varpi t} [u(0,\var)]_{\phi,\theta} \qquad \forall\,t \in (0,T],
\end{equation}
for some $\varpi > 0$ and $C_\theta$ that  depend on $d$, $|I|$, $\norm\M$, $c_{\scr N}$, $C_{\scr N}$,  $K_H$, $\gamma, \gamma'$, $\underbar c, \bar c$, $\phi$, and $\theta$.
\end{prop}

The combination of \Cref{prop_oscdec,thm_mainsmnl} --- that is, of the smoothing estimate \eqref{eq_gthm_usmoothing} and the decay estimate \eqref{eq_gthm_uosc} --- also provides an exponential long-time decay of the weighted oscillation of $u(0,\var)$,\footnote{That is, of the oscillation $0$-semi-norm of $u(0,\var)$ weighted by $\phi$, according to \Cref{def_oscsmn}.}. This is known from \cite[Theorem~5.1]{ForPor}, but for the sake of completeness we provide a statement and the short proof here as well, in our more general setting. 

\begin{cor}[Decay of oscillation for viscosity solutions] \label{cor_imprdec}
Suppose that the assumptions of Proposition \ref{thm_mainsmnl} and Proposition \ref{prop_oscdec} hold.
Then we have
\begin{equation} \label{eq_decoscimp}
|u(t,x) - u(t,y)|  \leq C  e^{-\varpi t} \norm{\phi^{-1} u(0,\var)}_\infty \bigl(\phi(x)+\phi(y)\bigr) \qquad \forall\,t \in (0,T], \ x,y \in \R^N,
\end{equation}
with the same $\varpi$ as in \eqref{eq_gthm_uosc} and some constant $C$ independent of $u$, $t$, and $T$.
\end{cor}

\begin{proof}
Estimate \eqref{eq_gthm_usmoothing'} with $t=1$ gives
\begin{equation} \label{eq_estu1}
|u(1,x) - u(1,y)| \lesssim \norm{\phi^{-1}u(0,\var)}_\infty \bigl(\phi(x)+\phi(y)\bigr) |x-y|,
\end{equation}
and due to \eqref{eq_ucompC3} we can replace $|x-y|$ in \eqref{eq_estu1} with $|x-y| \wedge 1$. Recalling \Cref{def_oscsmn}, it follows that $[u(1,\var)]_{\phi,1} \lesssim \norm{\phi^{-1}u(0,\var)}_\infty$, thus from \Cref{prop_oscdec} we deduce that \eqref{eq_decoscimp} holds for all $t \geq 1$; on the other hand, \eqref{eq_decoscimp} trivially holds for $t \in (0,1)$ as well, by \eqref{eq_ucompC3}.
\end{proof}

\section{ Main proofs} \label{proofs}

In this section we give the proof of \Cref{thm_mainsmnl,prop_oscdec}. We start by collecting a few  technical results based on maximum principle and coupling arguments, which play a crucial role.

\subsection{Maximum principle \& coupling arguments: preliminary toolbox}
For the sake of a shorter notation, for $f \colon (0,T) \times \R^N \to \R$, $p \in \R^N$, $\nu' \colon (0,T) \to \scr M_L(\R^N)$, and $M \colon (0,T) \to \GLN$, we write 
\begin{equation} \label{eq_deftildeIgen}
{\call I}_{\nu',M}[f](t,x,p) \defeq \int_{\R^N} \Bigl( f(t,x+z) - f(t,x) - p \cdot z \1_{|M(t)z| \leq 1} \Bigr) \, \nu'(t)(\di z)\,;
\end{equation}
note that, in particular, \rife{eq_deftildeI} implies
\[
{\call I}_\M[v](t,x) = \call I_{\nu_{\M},\M}[v](x,Dv(t,x)).
\]
We recall first a few well-known facts in viscosity solutions' theory (see  \cite{CIL}, and \cite{BarImbHP} for nonlocal operators). Whenever $v$ is a continuous viscosity solution of the equation
\be\label{ghost}
\partial_t v = \tr({\call A}(t)D^2v) + \call I_{\nu',M}[v]+ H(t,x) \cdot Dv 
\ee
for some ${\call A}$, $\nu'$, $M$, $H$, and whenever, for some smooth function $\zeta(t,x,y)$,  it happens that $v(t,x)-v(t,y)- \zeta(t,x,y)$ has a local maximum at $(\bar t, \bar x,\bar y)$, 
then, for every $n$ sufficiently large, there exist matrices  $X_n,Y_n \in \Sym^N$, and a constant $c_N$ only depending on $N$, satisfying  
\begin{equation} \label{eq_clvtineq}
  -  (n+ c_N\| D^2 \zeta(\bar t, \bar x, \bar y) \|) I_N \leq 
\begin{pmatrix}
X_n & 0 \\
0 & - Y_n
\end{pmatrix}
\leq D^2 \zeta(\bar t, \bar x, \bar y) + \frac1n(D^2\zeta(\bar t, \bar x, \bar y))^2
\end{equation}
and such that 
\begin{equation} \label{eq_pr1visth}\begin{split}
&\de_t \zeta (\bar t, \bar x, \bar y) - \tr(  \call A(\bar t) (X_n - Y_n)) 
\\
&- I_{\nu',M}[v](\bar t,\bar x,D_x\zeta (\bar t, \bar x, \bar y)) + I_{\nu',M}[v](\bar t, \bar y,-D_y\zeta (\bar t, \bar x, \bar y))
\\
&- H (\bar t,\bar x) \cdot D_x \zeta (\bar t, \bar x, \bar y) - H (\bar t,\bar y) \cdot D_y \zeta (\bar t, \bar x, \bar y) \leq 0\,.
\end{split}
\end{equation}
The matrices $X_n,Y_n$ are, roughly speaking, approximate values of the second order derivatives of $v$ at, respectively, $(\bar t, \bar x)$ and $(\bar t, \bar y)$. We refer e.g.\ to \cite[Theorem~2.3]{PorAdv} for a more precise statement of the previous property in a similar context as needed here. Obviously, 
by standard rescaling of test functions, if $u$ is a viscosity solution of \rife{eq_nleq}, then $u_\M$ (defined in \rife{eq_defutrans}) is a viscosity solution of \rife{eq_geneqcf}. Thus, the above considerations apply to $u$ as well as to $u_\M$. 

We now derive two key inequalities at the local maxima occurring for suitable test functions, obtained via coupling methods. The first concerns local diffusions while the second deal with nonlocal ones. In fact, the former is rather standard in the Ishii--Lions method \cite{IL}; see also \cite[Prop. 2.4]{PoPr}.

\begin{lem}[Local coupling] \label{lem_loccbr}
Let $v$ be a viscosity solution of \rife{ghost},  where $\call A \colon (0,T) \to \Sym_{\geq 0}(N) \setminus \{0\}$ and  $\lambda(t) > 0$ denotes the least positive eigenvalue of $\call A(t)$.  
Given functions $\Psi \colon \R^N \to \R$ and $\Phi,\phi \in C^{1,2}_\loc((0,T)\times \R^N)$, let 
\begin{equation} \label{eq_defzeta}
\zeta(t, x,y) \defeq (\phi(t,x) + \phi(t,y))\Psi(x-y) + \Phi(t,x) + \Phi(t,y), \qquad x,y \in \R^N, \,t\in (0,T)
\end{equation}
and suppose that
\begin{equation} \label{eq_argmaxcond}
(\bar t, \bar x, \bar y) \in \arg\max\{(t,x,y) \mapsto v(t,x) - v(t,y) - \zeta(t,x,y)\}.
\end{equation}
Assume that $\Psi$ is $C^2$ in a neighborhood of $\bar x - \bar y$.  Then, for any  sequence $\{(X_n,Y_n)\}_{n \in \N} \subseteq \Sym(N)^2$ satisfying \rife{eq_clvtineq},  for any subspace $W \subseteq V  = \Ran\, \call A(\bar t)$, and for any $R \in [0,2]$, one has 
\begin{equation} \label{eq_lemrcl} \begin{split}
&\limsup_{n \to +\infty} \tr(\call A(\bar t)(X_n - Y_n)) 
\\
&\leq \bigl( \tr(\call A(\bar t)D^2\phi(\bar t, \bar x)) + \tr(\call A(\bar t)D^2\phi(\bar t, \bar y)) \bigr) \Psi(\bar x- \bar y) + \tr(\call A(\bar t)D^2\Phi(\bar t, \bar x)) + \tr(\call A(\bar t)D^2\Phi(\bar t,\bar y))
\\
& \quad + R \lambda(\bar t) \Bigl( \bigl(\phi(\bar t, \bar x) + \phi(\bar t, \bar y) \bigr) D^2\Psi(\bar x - \bar y) \xi \cdot \xi + ( D\phi(\bar t,\var)|^{\bar x}_{\bar y} \cdot \xi ) ( D\Psi(\bar x - \bar y) \cdot \xi ) \Bigr) \1_{\call P_W(\bar x - \bar y) \neq 0},
\end{split}
\end{equation}
where $\call P_W$ denotes the orthogonal projection  onto $W$  and $\xi \defeq \wwidehat{\call P_W(\bar x-\bar y)}$.
\end{lem}

\begin{proof}
Given the orthogonal projection $\call P_W$, for $W$  a subspace of $V= \Ran\, \call A(\bar t)$, and 
setting $\xi \defeq \wwidehat{\call P_W(\bar x-\bar y)}$, we note that the following matrix is positive semidefinite, due to the assumptions on $\call A$, $\lambda_*$ and $R$:
\begin{equation} \label{eq_deffrkC}
\frk C \defeq \begin{pmatrix}
\call A(\bar t) & \call A(\bar t) - R\lambda(\bar t) \1_{\call P_W(\bar x - \bar y) \neq 0} \xi^{\otimes 2} \\
 \call A(\bar t) - R\lambda(\bar t) \1_{\call P_W(\bar x - \bar y) \neq 0} \xi^{\otimes 2} & \call A(\bar t)
\end{pmatrix}.
\end{equation}
Multiplying \eqref{eq_clvtineq} by $\frk C$ and taking the traces, we have
\[
\limsup_{n \to +\infty} \tr(\call A(\bar t)(X_n - Y_n)) \leq \tr\bigl(\call A(\bar t)(D^2_{ x  x} + D^2_{ y  y} + 2D^2_{ x  y})\zeta(\bar t, \bar x, \bar y)\bigr) - 2R\lambda(\bar t) \tr\bigl(\xi^{\otimes2}D^2_{ x  y}\zeta(\bar t, \bar x,\bar y)\bigr) \1_{\call P_W(\bar x - \bar y) \neq 0}.
\]
The conclusion now follows from direct computations of $D^2 \zeta$.
\end{proof}

By allowing $R \in [0,2]$, we are in fact including in the estimate \rife{eq_lemrcl} various different couplings, such as the \emph{coupling by reflection}, corresponding to $R = 2$, the \emph{coupling by projection}, when $R=1$, and \emph{synchronous coupling}, in the degenerate case of $R = 0$. Let us note that, however, any choice of $R\in(0,2)$ would be equivalent for our purposes; hence, with some abuse of terminology, we will only mention later the \emph{reflection} ($R \neq 0$) and \emph{synchronous} ($R = 0$) couplings.

In fact, it is also worth noting that the coupling induced by \eqref{eq_deffrkC} (even with $R = 2$ and $W = V$) is not exactly the probabilistic coupling by reflection in the sense of, e.g., \cite{ChenLi}; indeed, that would be, more naturally,
\[
\sigma \begin{pmatrix} I  & I-2(\bar x-\bar y)^{\otimes2} \\ I-2(\bar x-\bar y)^{\otimes2} & I \end{pmatrix}\sigma^*
\qquad \text{with} \quad \sigma\sigma^* = \call A.
\]
Nevertheless, we find the one provided by \eqref{eq_deffrkC} to be more convenient  inside an Ishii--Lions argument making use of distance-like test functions --- one main reason being that $\sigma(\bar x-\bar y)^{\otimes2}\sigma^*$ is not a projection, while $\xi \otimes \xi$ is.

We move now to the nonlocal counterpart of the above result. 

\begin{lem}[Nonlocal coupling] \label{lem_dvnl}
Let $v$, $\Psi$, $\Phi$ and $\phi$ be as in \Cref{lem_loccbr}. Given a subspace $V \subseteq \R^N$ and $\bar\nu \in \scr M_L(V)$ radial, suppose that there is $\lambda \colon (0,T) \to \R_+$ such that $\nu \geq \lambda \bar\nu$. Let $\zeta$ be as in \eqref{eq_defzeta} and suppose that $(\bar x,\bar y)$ satisfy \eqref{eq_argmaxcond}. 
Given $R\geq 0$, assume that $\Psi$ is $C^2$ in $\bar x-\bar y+B_R$. 
Then, for any subspace $W \subseteq V$, letting $\xi \defeq \wwidehat{\call P_W(\bar x-\bar y)}$, one has
\begin{align} \label{eq_lem1ineq1}
& {\call I}_{\nu,M}[v](\bar t,\bar x,D_x\zeta(\bar t,\bar x,\bar y)) -  {\call I}_{\nu,M}[v](\bar t,\bar y,-D_y\zeta(\bar t,\bar x,\bar y))
\notag \\
&\leq {\call I}_{\nu,M}[\Phi](\bar t,\bar x,D\Phi(\bar t,\bar x)) + {\call I}_{\nu,M}[\Phi](\bar t,\bar y,D\Phi(\bar t,\bar y)) + \bigl( {\call I}_{\nu,M}[\phi](\bar t,\bar x,D\phi(\bar t,\bar x)) + {\call I}_{\nu,M}[\phi](\bar t,\bar y,D\phi(\bar t,\bar y)) \bigr) \Psi(x-y)
\notag \\
&\quad + {2\lambda(\bar t)} \biggl( \int_{\call P_VB_R} \bar\phi(\bar t,z) |\xi \cdot z|^2 \int_0^1 (1-\tau) D^2\Psi\bigl(\bar x-\bar y+2\tau(\xi\cdot z )\xi\bigr) \xi \cdot \xi\,\di\tau\,\bar\nu(\di z)
\notag \\
&\quad  +  \Bigl( \int_{\call P_V B_R} |\var|^{2}\,\di \bar\nu \Bigr) C_{D\phi(\bar t,\var),R} |D\Psi(\bar x-\bar y) \cdot \xi| \biggr) \1_{\call P_W(\bar x-\bar y) \neq 0} ,
\end{align}
where
\begin{equation} \label{eq_constnlcoup}
\bar\phi(\bar t,z) \defeq 2\bigl(\phi(\bar t,\bar x+z) + \phi(\bar t,\bar y+z-2(\xi\cdot z)\xi)\bigr), \quad C_{D\phi,R} \defeq \norm{D\phi(\bar t,\var)}_{\infty;\bar x+B_R} + \norm{D\phi(\bar t,\var)}_{\infty;\bar y+B_R},
\end{equation}
and the operator ${\call I}$ is defined by \eqref{eq_deftildeIgen}.
\end{lem}

\begin{proof}
We follow the idea in the proof of \cite[Lemma~3.4]{PorAdv}; see also \cite[Lemma~A.2]{EJP}. As $(\bar t,\bar x, \bar y)$ is a maximum point of $v(t,x) - v(t,y) - \zeta(t,x,y)$, we have (omitting the dependence on $t$, and the bar over $x$ and $y$)
\[
v|^{x+z}_x - v|^{y+z'}_y \leq \zeta|^{(x+z,y+z')}_{(x,y)} \qquad \forall\, z,z' \in \R^N,
\]
which can be rewritten as
\begin{equation} \label{eq_maxineq}
\begin{split}
v|^{x+z}_x - v|^{y+z'}_y & \leq  \Phi|^{x+z}_x + \Phi|^{y+z'}_y  + (\phi|^{x+z}_x + \phi|^{y+z'}_y)\Psi(x-y)
\\
& \quad  \quad  + (\phi(x+z) + \phi(y+z')) \Psi|^{x-y+z-z'}_{x-y} \qquad \forall\, z,z' \in \R^N.
\end{split}
\end{equation}
Using that
\begin{equation*}
\begin{gathered}
D_x \zeta(x,y) = D\Phi(x) + D\phi(x)\Psi(x-y) + \bigl(\phi(x) + \phi(y)\bigr)D\Psi(x-y),
\\
D_y \zeta(x,y) = D\Phi(y) + D\phi(y)\Psi(x-y) - \bigl(\phi(x) + \phi(y)\bigr)D\Psi(x-y),
\end{gathered}
\end{equation*}
we have, for all $z,z' \in \R^N$,
\begin{equation} \label{eq_maxineq'}
\begin{split}
&v|^{x+z}_x - D_x\zeta(x,y) \cdot z \chi(z) -  v|^{y+z'}_y -D_y\zeta(x,y) \cdot z' \chi(z') 
\\
&\leq \Phi|^{x+z}_x - D\Phi(x) \cdot z \chi(z) + \Phi|^{y+z'}_y - D\Phi(y) \cdot z' \chi(z')
\\
&\quad + \bigl( \phi|^{x+z}_x - D\phi(x) \cdot z \chi(z) + \phi|^{y+z'}_y - D\phi(y) \cdot  z' \chi(z') \bigr)\Psi(x-y) 
\\
&\quad + \bigl(\phi(x+z) + \phi(y+z')\bigr) \Bigl( \Psi|^{x-y+z-z'}_{x-y} - D\Psi(x-y)\cdot \bigl(z \chi(z) -z' \chi(z')\bigr) \Bigr)
\\
&\quad + \bigl( \phi|^{x+z}_x + \phi|^{y+z'}_y \bigr) D\Psi(x-y)\cdot(z \chi(z) -z' \chi(z') ),
\end{split}
\end{equation}
for any arbitrary function $\chi$.
Therefore, given a L\'{e}vy measure $\mu$, we can integrate against $(I_N,I_N)_\sharp\mu$ the inequality \eqref{eq_maxineq'} with $\chi(z) = \1_{|Mz|\leq 1}$, to obtain (omitting the dependence on $M$)
\begin{equation} \label{eq_l1e1}
\begin{split}
&{\call I}_{\mu}[v](x,D_x\zeta(x,y)) - {\call I}_{\mu}[v](y,-D_y\zeta(x,y))
\\
& \leq {\call I}_{\mu}[\Phi](x,D\Phi(x)) + {\call I}_{\mu}[\Phi](y,D\Phi(y)) + \bigl( {\call I}_{\mu}[\phi](x,D\phi(x)) + {\call I}_{\mu}[\phi](y,D\phi(y)) \bigr) \Psi(x-y).
\end{split}
\end{equation}
On the other hand, given two L\'evy measure $\mu$ and $\mu'$ that are invariant under reflection through the origin (i.e., $(-I_N)_\sharp\mu = \mu$), then
\[
\int_{\R^N} z \1_{|Mz| > 1}\,\mu(\di z) = 0,
\]
hence, by integrating \eqref{eq_maxineq'} with $\chi = 1$, against a coupling $\pi$ of $\mu$ and $\mu'$, we have
\begin{equation} \label{eq_l1e0}
\begin{split}
&{\call I}_{\mu}[v](x,D_x\zeta(x,y)) - {\call I}_{\mu'}[v](y,-D_y\zeta(x,y))
\\
&\leq {\call I}_{\mu}[\Phi](x,D\Phi(x)) + {\call I}_{\mu'}[\Phi](y,D\Phi(y))
+ \bigl( {\call I}_{\mu}[\phi](x,D\phi(x)) + {\call I}_{\mu'}[\phi](y,D\phi(y)) \bigr) \Psi(x-y)
\\
&\quad + \int_{\R^N\times\R^N} \bigl(\phi(x+z) + \phi(y+z')\bigr) \bigl( \Psi|^{x-y+z-z'}_{x-y} - D\Psi(x-y)\cdot (z -z' ) \bigr) \pi(\di z,\di z')
\\
&\quad + \int_{\R^N \times \R^N} \bigl( \phi|^{x+z}_x + \phi|^{y+z'}_y \bigr) D\Psi(x-y)\cdot(z-z')\,\pi(\di z,\di z').
\end{split}
\end{equation}
Define now, for any ball $B_R \subset \R^N$,
\[
m \defeq \bar\nu \resmes \call P_VB_R, \qquad \tilde\nu \defeq \nu - \lambda m;
\]
note that $m$ and $\tilde\nu$ are both non-negative L\'{e}vy measures due to the standing assumptions, also $(-I_N)_\sharp m = m$. Now, given a subspace of $W \subseteq V$, if $\call P_W(\bar x-\bar y) \neq 0$ we let
\[
\xi \defeq \wwidehat{\call P_W(x-y)}, \qquad \call O \defeq I_N - 2\xi^{\otimes 2};
\]
in this case, $\call O$ is a reflection with respect to the line $\R\xi \subset V$, so we also have $\call O_\sharp m = m$ (recall that $\bar \nu$ is radial). Then, letting $\pi = (I_N,\call O)_\sharp m$ in \eqref{eq_l1e0}, we have
\begin{equation*}
\begin{split}
&{\call I}_{m}[v](x,D_x\zeta(x,y)) - {\call I}_{\call O_\sharp m}[v](y,-D_y\zeta(x,y))
\\
&\leq {\call I}_{m}[\Phi](x,D\Phi(x)) + {\call I}_{\call O_\sharp m}[\Phi](y,D\Phi(y)) + \bigl( {\call I}_{m}[\phi](x,D\phi(x)) + {\call I}_{\call O_\sharp m}[\phi](y,D\phi(y)) \bigr) \Psi(x-y)
\\
&\quad + \int_{\R^N} (\phi(x+z) + \phi(y+ \call Oz)) \bigl( \Psi|^{x-y+2(\xi\cdot z)\xi}_{x-y} - 2(\xi \cdot z) D\Psi(x-y) \cdot \xi \bigr) m(\di z) 
\\
&\quad + 2\int_{\R^N} (\phi|^{x+z}_x + \phi|^{y+\call Oz}_y) (D\Psi(x-y) \cdot \xi)(\xi \cdot z) \,m(\di z).
\end{split}
\end{equation*}
By Taylor's expansion with integral reminder, 
\begin{align*}
\Psi|^{x-y+2(\xi\cdot z)\xi}_{x-y} - 2(\xi \cdot z) D\Psi(x-y) \cdot \xi &= 4 \int_0^1 (1-\tau) |\xi \cdot z|^2 D^2\Psi\bigl(x-y+2\tau(\xi\cdot z )\xi\bigr)  \xi \cdot \xi \,\di\tau,
\\
\phi|^{x+z}_x + \phi|^{y+\call Oz}_y &= \int_0^1 \bigl( D\phi(x+\tau z) \cdot z + D\phi(y+\tau\call Oz) \cdot \call Oz \bigr)\,\di\tau;
\end{align*}
by using the latter identity, we can estimate
\begin{equation*} \begin{split}
&\int_{\R^N} (\phi|^{x+z}_x + \phi|^{y+\call Oz}_y) (D\Psi(x-y) \cdot \xi) (\xi \cdot z) \,m(\di z)
\\
&\leq \bigl(\norm{D\phi}_{\infty;x+ \call P_VB_R} + \norm{D\phi}_{\infty;y+ \call P_VB_R}\bigr) |D\Psi(x-y) \cdot \xi| \int_{\R^N} |(\xi \cdot z)z| \,m(\di z),
\end{split}\end{equation*}
where
\[
\int_{\R^N} |(\xi \cdot z)z| \,m(\di z) \leq \int_{\call P_V B_R} |\var|^{2}\,\di \bar\nu. 
\]
We have obtained that
\begin{equation*}
\begin{split}
&{\call I}_{m}[v](x,D_x\zeta(x,y)) - {\call I}_{m}[v](y,-D_y\zeta(x,y))
\\
&\leq {\call I}_{m}[\Phi](x,D\Phi(x)) + {\call I}_{m}[\Phi](y,D\Phi(y)) + \bigl( {\call I}_{m}[\phi](x,D\phi(x)) + {\call I}_{m}[\phi](y,D\phi(y)) \bigr) \Psi(x-y)
\\
&\quad + 4 \int_{\R^N} \bigl(\phi(x+z) + \phi(y+z-2(\xi\cdot z)\xi)\bigr) |\xi \cdot z|^2  \int_0^1 (1-\tau) D^2\Psi\bigl(x-y+2\tau(\xi\cdot z )\xi\bigr)  \xi \cdot \xi \,\di\tau\,m(\di z)
\\
&\quad + 2 \Bigl( \int_{\call P_V B_R} |\var|^{2}\,\di \bar\nu \Bigr) \bigl(\norm{D\phi}_{\infty;x+ \call P_VB_R} +\norm{D\phi}_{\infty;y+ \call P_VB_R}\bigr) |D\Psi(x-y) \cdot \xi|\,.
\end{split}
\end{equation*}
Summing now this inequality multiplied by $\lambda$ with \eqref{eq_l1e1} where we choose $\mu = \tilde \nu$, we obtain \eqref{eq_lem1ineq1}. The case when $\call P_W( x- y) = 0$ easily follows by taking $\call O = I_N$ or, equivalently, $R=0$.
\end{proof}

Note that a synchronous coupling is included in this nonlocal case as well, since we can take $R=0$.

\begin{rmk}[A natural splitting for a ``distance'' function] \label{rmk_structtest}
In order to exploit the structure of the term inside the integral in \eqref{eq_lem1ineq1} --- or its local analogue from \Cref{lem_loccbr} --- we are going to consider $\Psi$ of the form
\begin{equation} \label{eq_fromPsirmk}
\Psi = \psi \circ \rho, \qquad \rho(z) = \rho_*(\call P_Wz) + \rho_\perp(\call P_W^\perp z),
\end{equation}
for some $\psi \colon \R_{\geq 0} \to \R$, $\rho_*,\rho_\perp \colon \R^N \to \R_{\geq 0}$ of class $C^2$ in a neighborhood of $\rho(\bar x-\bar y)$, $\call P_W(\bar x - \bar y)$ and $\call P_W^\perp(\bar x - \bar y)$, respectively. In this case, letting $r \defeq |\call P_W(\bar x - \bar y)|$, by direct computations one has
\[
D\Psi(x-y) \cdot \xi = \psi'(\rho) D\rho_*(r\xi) \cdot \xi
\]
and
\[
D^2\Psi\bigl(\bar x-\bar y+2\tau(\xi\cdot z )\xi\bigr) \xi \cdot \xi = \psi'(\bar\rho_\tau)D^2\rho_*((r+2\tau \xi\cdot z)\xi) \xi \cdot \xi + \psi''(\bar\rho_\tau)| D\rho_*((r+2\tau\xi\cdot z)\xi) \cdot \xi |^2,
\]
where
\[
\bar\rho_\tau \defeq \rho_*((r+2\tau\xi\cdot z)\xi) + \rho_\perp(\call P_W^\perp(\bar x - \bar y)).
\]
\end{rmk}

Finally, we conclude this preliminary part with a technical  lemma about interior maxima of perturbed test functions, which is needed later to justify the use of the maximum principle for the augmented problem, as is classical for doubling variables methods.

\begin{lem}[Interior maxima of regularized penalized differences] \label{lem_exmax}
Let $\Phi,\phi \colon [0,T) \times \R^N \to (1,+\infty)$, let $\chi \in C^0([0,T);\R_{>0})$ satisfy $\chi(t) \to +\infty$ as $t \uto T$, and let $\psi \in C^0([0,\delta);\R_{\geq 0})$, for some $\delta \in \R \cup \{+\infty\}$, be a bounded and increasing function such that $\psi(0) = 0$ and, for some $\theta \in (0,1]$,
\begin{equation} \label{eq_psioverrt}
 c_{\psi,\theta,\delta} \defeq \liminf_{\rho \in (0,2\delta)} \frac{\psi(\varrho)}{\varrho^\theta \wedge 1} < +\infty.
\end{equation}
Let $v \in C^0([0,T]\times\R^N)   \cap L^\infty((0,T);C^{0,1}(\R^N))$.
For all $t \in (0,T)$, $x,y \in \R^N$, $\epsilon \in [0,1)$, and $\scr N_\epsilon$ satisfying \eqref{ass_Neps}--\eqref{eq_condcallN'},
define
\begin{gather*}
\zeta_\epsilon(t,x,y) \defeq K_0\bigl(\phi(t,x)+\phi(t,y)\bigr)\psi(\rho_\epsilon) + \Phi(t,x) + \Phi(t,y) + \chi(t), \qquad \rho_\epsilon \defeq \scr N_\epsilon(x-y),
\\
\call U_\epsilon(t,x,y) \defeq v(t,x) - v(t,y) - \zeta_\epsilon(t,x,y),
\end{gather*}
for some $K_0 >0$. 
Assume that $\call U_0$ has a positive supremum $\frk m_0$ on $(0,T) \times \Omega$, with
\begin{equation} \label{eq_defOmega}
\Omega = \{ (x,y) \in \R^N \times \R^N :\, \rho_0 < \delta \}.
\end{equation}
Then there is $\epsilon_0\in(0,1]$ such that for any $\epsilon \in [0,\epsilon_0)$ the function $\call U_\epsilon$ has a positive supremum $\frk m_\epsilon \geq \frk m_0$ on $(0,T) \times \Omega$ which is an interior maximum, in either one of the following two cases: 
\begin{equation} 
\begin{gathered}\label{eq_choiceK0NL'}
\lim_{t \dto 0}\, [v(t,\var)]_{\phi(t,\var),\theta} = 0, \qquad 
\lim_{|x| \to +\infty} \sup_{t \in [t_0,T]} \frac{|v(t,x)|}{\Phi(t,x) \vee \phi(t,x)} = 0 \quad \forall\, t_0 \in (0,T], \\
K_0 \geq \frac{\sup_{t\in(0,T)} [v(t,\var)]_{\phi(t,\var),0}}{\psi(\delta/2)}\,\1_{\delta \neq +\infty} \vee \frac{2 \sup_{t \in (0,T)} [v(t,\var)]_{\phi(t,\var),\theta}}{(c_{\scr N} \wedge 1)c_{\psi,\theta,\delta}}.
\end{gathered}
\end{equation}
or
\begin{equation} \label{eq_choiceK0NL''}
\lim_{|x| \to +\infty} \sup_{t \in [0,T]} \frac{|v(t,x)|}{\Phi(t,x) \vee \phi(t,x)} = 0, 
\qquad
K_0 \geq \frac{\sup_{t\in(0,T)} [v(t,\var)]_{\phi(t,\var),0}}{\psi(\delta/2)}\,\1_{\delta \neq +\infty} \vee \frac{2 [v(0,\var)]_{\phi(0,\var),\theta}}{(c_{\scr N} \wedge 1)c_{\psi,\theta,\delta}}
\end{equation}
\end{lem}

\begin{proof}
By   \eqref{eq_condcallN'} we have $\rho_\epsilon\leq \rho_0$, so  the monotonicity of $\psi$ implies that $\call U_\epsilon \geq \call U_0$ for any $\epsilon \in (0,1]$; then $\call U_\epsilon$ has a positive supremum as well. Consider now the scenario \eqref{eq_choiceK0NL'}. Note that, from \rife{eq_condcallN'}, we have $|x-y|\leq  c_{\scr N}^{-1} \rho_0$; in addition, we have  $\rho_0 \leq 2\rho_\epsilon$ whenever $\rho_0 \geq 2C_{\scr N}\epsilon$. Therefore, letting $\ell_v \defeq \sup_{t\in(0,T)} [v(t,\var)]_{0,1}$, we estimate
\begin{equation*} \begin{split}
v(t,x) - v(t,y) &\leq 
{[v(t,\var)]_{\phi(t,\var),\theta}} \bigl(\phi(t,x) + \phi(t,y)\bigr) (|x-y|^\theta \wedge 1) \1_{\rho_0 \geq 2C_{\scr N}\epsilon} +   \ell_v |x-y| \1_{\rho_0 < 2C_{\scr N}\epsilon} 
\\
&\leq 
(1 \vee (2c_{\scr N}^{-1})^\theta )  {[v(t,\var)]_{\phi(t,\var),\theta}} \bigl(\phi(t,x) + \phi(t,y)\bigr) (\rho_\epsilon^\theta \wedge 1) \1_{\rho_0 \geq 2C_{\scr N}\epsilon} +  2C_{\scr N}c_{\scr N}^{-1} \ell_v \epsilon  \1_{\rho_0 < 2C_{\scr N}\epsilon}  
\\ 
&\leq {(1 \vee (2c_{\scr N}^{-1})^\theta ) } c_{\psi,\theta,\delta}^{-1} {[v(t,\var)]_{\phi(t,\var),\theta}} \bigl(\phi(t,x) + \phi(t,y)\bigr) \psi(\rho_\epsilon)  + 2C_{\scr N}c_{\scr N}^{-1} \ell_v \epsilon
\end{split}\end{equation*}
where we also used \eqref{eq_psioverrt} in the last inequality.
Since $\frac12 (c_{\scr N}\wedge 1)\leq (\frac{c_{\scr N}}2)^\theta\wedge 1$ we deduce the implication
\[
[v(t,\var)]_{\phi(t,\var),\theta} \leq \frac12{K_0(c_{\scr N}\wedge 1)c_{\psi,\theta,\delta}} \quad \implies \quad \call U_\epsilon(t,\var) \leq 2C_{\scr N}c_{\scr N}^{-1} \ell_v \epsilon - \chi(t),
\]
so the first and third conditions in \eqref{eq_choiceK0NL'} ensure that there is $t_0>0$ (independent of $\epsilon$) such that $\call U_\epsilon\leq 0$ in $[0,t_0)$, for any $\epsilon>0$ small enough. Hence, using also 
the second condition in \eqref{eq_choiceK0NL'}, we deduce that each $\frk m_\epsilon$ is actually a maximum on $[0,T] \times \bar\Omega$, obtained for $t\geq t_0$. Since both $\call U(t,x,x)$ and $\call U_\epsilon(t,x,x)$ are positive, it cannot be attained for $\rho = 0$, and since $\chi(t) \to +\infty$ as $t \uto T$, it cannot be attained at $t = T$. On the other hand, if also $2C_{\scr N}  \epsilon < \delta$, then, by the monotonicity of $\psi$ and \eqref{eq_condcallN'}, we have $\psi(\rho_\epsilon) \geq \psi(\delta/2)$ when $\rho = \delta$, so no maximum can be attained for $\rho = \delta$ either, because of the lower bound on $K_0$. This concludes the proof that $\frk m_\epsilon$ is an interior maximum.

Similarly we conclude if the scenario \eqref{eq_choiceK0NL''} holds. In fact, following the same argument as before, the  second condition now implies that $\frk m_\epsilon$ cannot be attained at $t=0$, due to the lower bound on $K_0$. 
The rest of the proof works as in the former setting. 
\end{proof}

\subsection{Proof of the general estimates.}

We are now ready to give the proof of the bounded-to-Lipschitz smoothing estimate.

\begin{proof}[\bfseries Proof of Proposition \ref{thm_mainsmnl}]
Given $\vartheta \in (0,1)$ to be determined, let $\smash{\psi(\varrho) \defeq \varrho(1 - \frac1{2(1+\vartheta)}\varrho^\vartheta)}$; note that
\[
\psi(\varrho)\varrho^{-1}, \psi'(\varrho) \in (\tfrac12,1] \quad \text{if} \ \varrho < 1,
\]
so $\psi$ satisfies the assumptions of \Cref{lem_exmax} for $\delta \leq 1$ (to be determined) and $\theta = 1$. Let $\tilde\phi \defeq \phi_{\M}$ and $\tilde\Phi \defeq \Phi_{\M}$; note that, without loss of generality, we can suppose that $\phi$ and $\Phi$ are supersolutions of \eqref{eq_nleq}, and thus $\tilde\phi \defeq \phi_{\M}$ and $\tilde\Phi \defeq \Phi_{\M}$ are supersolutions of \eqref{eq_geneqcf}.\footnote{Indeed, if $(\de_t - \call L^H) \phi \geq -k_0$, it suffices to consider $e^{k_0t}\phi(t,x)$ in place of $\phi(t,x)$, and likewise for $\Phi$; the verification that the change of coordinates preserves the supersolution property is elementary.}

For $\epsilon \in [0,1)$, $\eta>0$, define
\[
\zeta_\epsilon(t,x,y) \defeq K_0\bigl(\tilde\phi(t,x)+\tilde\phi(t,y)\bigr)\psi(\rho_\epsilon) + \eta \Bigl( \tilde\Phi(t,x) + \tilde\Phi(t,y) + \frac1{T-t} \Bigr), 
\]
where
\[
\rho_\epsilon = \rho_\epsilon(x-y) \defeq \sum_{j \in I} |\call P_{W_j}(x-y)|_\epsilon + \kappa\scr N_\epsilon({\call P}_{\call V}^\perp(x-y)), \qquad |\var|_\epsilon \defeq \sqrt{|\var|^2 + \epsilon^2} - \epsilon,
\]
and $K_0,\kappa > 0$ are constants to be determined. As $|\var|_0 = |\var|$, we will also use the lighter notations $\rho = \rho_0$, $\scr N = \scr N_0$, and $\zeta = \zeta_0$.

Given $v \defeq u_{\M}$ and $\Omega$ as \eqref{eq_defOmega} (with $\delta \leq 1$), we choose $K_0$ as in \eqref{eq_choiceK0NL'}, and we want to prove that we have
\[
\call U(t,x,y) \defeq v(t,x) - v(t,y) - \zeta(t,x,y) \leq 0 \quad\ \forall\, (t,x,y) \in (0,T) \times \Omega.
\]
This is sufficient to conclude, since by letting $\eta \dto 0$ and using $\psi'\leq 1$ we get  
\[
v(t,x) - v(t,y) \leq K_0 \bigl(\tilde\phi(t,x)+\tilde\phi(t,y)\bigr)\rho \quad \text{if $\rho < 1$},
\]
which yields, due to \eqref{eq_condcallN'} and the value of $K_0$,  
\[
v(t,x) - v(t,y) \leq C 
\sup_{t\in(0,T)} [v(t,\var)]_{\tilde\phi(t,\var),0} \bigl(\tilde\phi(t,x)+\tilde\phi(t,y)\bigr) |x - y| \qquad  \forall\, t \in (0,T], \ x,y \in \R^N
\]
(where $C$ only depends on $\vartheta$, $\delta$, $|I|$, $\kappa$, and $C_{\scr N}$). This implies \eqref{eq_gthm_usmoothing} by noting that $\sup_{t \in (0,T)}[v(t,\var)]_{\tilde\phi(t,\var),0} = \sup_{t \in (0,T)} [u(t,\var)]_{\phi(t,\var),0}$.

Seeking for a contradiction, we suppose that $\call U$ has a positive supremum $\frk m$ on $(0,T) \times \Omega$. By \Cref{lem_exmax}, that is an interior maximum, and by choosing $\epsilon$ sufficiently small, $\call U_\epsilon(t,x,y) \defeq v(t,x) - v(t,y) - \zeta_\epsilon(t,x,y)$ also has a positive interior maximum $\frk m_\epsilon$. From this point on, functions will be understood to be valued at a point of $\arg\max_{(0,T) \times \Omega} \call U_\epsilon$.

By viscosity theory, as mentioned before (cf. \rife{eq_pr1visth}) we have
\begin{equation} \label{eq_pr1visth2}\begin{split}
&\de_t \zeta_\epsilon(t,x,y) - \tr(A_\M(t)(X_n - Y_n)) 
\\
&- {\call I}_{\nu_{\M}(t),\M(t)}[v(t,\var)](t,x,D_x\zeta_\epsilon(t,x,y)) + {\call I}_{\nu_{\M}(t),\M(t)}[v(t,\var)](t, y,-D_y\zeta_\epsilon(t,x, y))
\\
&- H_\M(t,x) \cdot D_x \zeta_\epsilon(t,x,y) - H_\M(t,y) \cdot D_y \zeta_\epsilon(t,x,y) \leq 0,
\end{split}
\end{equation}
for a sequence $(X_n,Y_n) \in \Sym(N)$ satisfying \eqref{eq_clvtineq} (with $\zeta_\epsilon(t,\var)$ in place of $\zeta$).
 We now apply \Cref{lem_loccbr} and \Cref{lem_dvnl}, where we choose $\call A = A_\M$, $\lambda = s_j^{-1}$, $\nu = \nu_j$, $V=V_j'$, $W= W_j \1_{j \in I}$, $s=s_j$, and $R = R_j \in [0,1]$  to be fixed; from \eqref{eq_pr1visth} we obtain 
\begin{align} \label{eq_contreqlns}
&\eta \Bigl( (\de_t - \call L_\M^H)\tilde\Phi(t,x) + (\de_t - \call L_\M^H)\tilde\Phi(t,y) + \frac{1}{(T-t)^2} \Bigr) + K_0\bigl( (\de_t - \call L_\M^H)\tilde\phi(t,x) + (\de_t - \call L_\M^H)\tilde\phi(t,y) \bigr) \psi(\rho_\epsilon)
\notag \\
&\leq K_0\bigl(\tilde\phi(t,x)+\tilde\phi(t,y)\bigr)\psi'(\rho_\epsilon) H_\M(t,\var)\bigr|^x_y \cdot D_x\rho_\epsilon
\notag \\
&\quad + \frac{R_0 K_0}{s_0(t)} \Bigl( \bigl(\tilde\phi(t,x) + \tilde\phi(t,y) \bigr) D^2(\psi \circ \rho_\epsilon)(x - y) \xi_0 \cdot \xi_0 + C_{D\tilde\phi(t,\var),0}  |D(\psi \circ \rho_\epsilon)( x- y) \cdot \xi_0| \Bigr) \1_{\call P_{W_0}(x - y) \neq 0 \wedge 0 \in I}
\notag \\
&\quad + \sum_{j \in I \setminus \{0\}} \frac{2K_0}{s_{j}(t)} \biggl(  \int_{\call P_{V_j'} B_{R_j}} \bar{\tilde\phi}(t,z) |\xi_j \cdot z|^2 \int_0^1 (1-\tau) D^2(\psi\circ\rho_\epsilon)( x- y+2\tau(\xi_j\cdot z)\xi_j ) \xi_j \cdot \xi_j\, \di\tau\,\bar\nu_{j}(\di z)
\notag \\
&\quad + \Bigl( \int_{\call P_{V_j'} B_{R_j}} |\var|^{2}\,\di \bar\nu_j \Bigr) C_{D\tilde\phi(t,\var),R_j}  |D(\psi \circ \rho_\epsilon)( x- y) \cdot \xi_j| \biggl) \1_{\call P_{W_j}(x-y) \neq 0},
\end{align}
where $\xi_j \defeq \wwidehat{\call P_{W_j}(x-y)}$,  while $\bar{\tilde\phi}$ and $C_{D\tilde\phi(t,\var),R_j}$ are defined as in \eqref{eq_constnlcoup}.

By \eqref{eq_pshar},  letting $\tilde C \defeq (1 \vee \sup_{(0,T)}\norm{\M})C_{\sup_{(0,T)}\norm{\M}}$, we have\footnote{It suffices to note that, for any $R >0$, one has $\M(t)(y+B_R) \subseteq \M(t)y + B_{\norm{\M(t)}R}$, so 
\[
\sup_{y + B_R} D\tilde\phi(t,\var) \leq \norm{\M(t)} C_{\norm{\M(t)}R} \tilde\phi(t,y), \qquad \inf_{z \in B_R} \tilde\phi\bigl(t,y+z-2(\xi\cdot z)\xi\bigr) = \inf_{y+B_R}\tilde\phi(t,\var) \geq C^{-1}_{\norm{\M(t)}R} \tilde\phi(t,y),
\]
and likewise for $\tilde\phi(t,x)$.}
\begin{equation} \label{eq_contrtLyap}
\tilde C^{-1} C_{D\tilde\phi(t,\var),R} \leq \tilde\phi(t,x)+\tilde\phi(t,y) \leq \tilde C \inf_{z \in B_{R}} \bar{\tilde\phi}(t,z) \qquad \forall\, R \leq 1;
\end{equation}
also, since $\bar\nu_j$ is radial, \eqref{eq_condnurad1} gives
\[
\int_{\call P_{V_j'} B_{R_j}} |\var|^{2}\,\di \bar\nu_j = \int_{B_{R_j} \cap V_j} |\var|^{2}\,\di \bar\nu_j \leq \bar cR_j^{\beta_j}
\]
Hence, from \rife{eq_contreqlns} we deduce
\begin{align} \label{eq_contreqlns'0}
&0 \leq \psi'(\rho_\epsilon) H_\M(t,\var)\bigr|^x_y \cdot D_x\rho_\epsilon
 \notag\\
 &\quad+ \frac{ R_0}{s_0(t)} \Bigl( D^2(\psi \circ \rho_\epsilon)(x - y) \xi_0 \cdot \xi_0 + \tilde C  |D(\psi \circ \rho_\epsilon)( x- y) \cdot \xi_0| \Bigr) \1_{\call P_{W_0}(x - y) \neq 0 \wedge 0 \in I}
\notag \\
&\quad + \sum_{j \in I \setminus \{0\}} \frac{2}{s_{j}(t)} \biggl(  \int_{\call P_{V_j'} B_{R_j}} \tilde C^{-1} |\xi_j \cdot z|^2 \int_0^1 (1-\tau) D^2(\psi\circ\rho_\epsilon)( x- y+2\tau(\xi_j\cdot z)\xi_j ) \xi_j \cdot \xi_j\, \di\tau\,\bar\nu_{j}(\di z)
\notag \\
&\quad + \bar c \tilde C R_j^{\beta_j} |D(\psi \circ \rho_\epsilon)( x- y) \cdot \xi_j| \biggl) \1_{\call P_{W_j}(x-y) \neq 0},
\end{align}
provided that
\begin{equation} \label{eq_pr1trneg}
D^2(\psi\circ\rho_\epsilon)( x- y+2\tau(\xi_j\cdot z)\xi_j ) \xi_j \cdot \xi_j \leq 0.
\end{equation}
To show \eqref{eq_pr1trneg}, we note that, letting
\[
r_j \defeq |\call P_{W_j}(x-y)|, \qquad \rho_{\epsilon,\tau} = \sum_{j \in I} |r_j+2\tau\xi_j \cdot z|_\epsilon + \scr N_\epsilon(\call P_{\call V}^\perp(x-y)),
\]
and following \Cref{rmk_structtest}, we have
\begin{equation*}
D^2(\psi \circ \rho_\epsilon)(x-y+2\tau(\xi_j \cdot z)\xi_j) \xi_j \cdot \xi_j =  \frac{\epsilon^2}{(|r_j+2\tau \xi_j \cdot z|^2 + \epsilon^2)^{\frac32}} \, \psi'(\rho_{\epsilon,\tau}) + \frac{|r_j+2\tau\xi_j\cdot z|^2}{|r_j+2\tau\xi_j\cdot z|^2 + \epsilon^2}\, \psi''(\rho_{\epsilon,\tau}).
\end{equation*}
Note that if $2|r_j+2\tau\xi\cdot z| \geq r_j \geq \sqrt{\epsilon}$, then we have, for $\epsilon$ small enough, 
\[
\frac{\epsilon^2}{(|r_j+2\tau \xi \cdot z|^2 + \epsilon^2)^{\frac32}} \leq \frac{\epsilon^2}{(\frac\epsilon{4}+ \epsilon^2)^{\frac32}} \leq  8\sqrt{\epsilon}\leq  16( r_j - \epsilon ) \leq 32|r_j+2\tau\xi_j \cdot z|_\epsilon \leq 32\rho_{\epsilon,\tau},
\]
and   the above inequality $\sqrt\epsilon \leq 4|r_j+2\tau\xi_j\cdot z|_\epsilon$ also gives, for $\epsilon$ small enough,
\[
\frac{|r_j+2\tau\xi_j\cdot z|^2}{|r_j+2\tau\xi_j\cdot z|^2 + \epsilon^2} \geq \Bigl( \frac{|r_j+2\tau\xi_j\cdot z|_\epsilon}{|r_j+2\tau\xi_j\cdot z|_\epsilon + \epsilon}\Bigr)^2 \geq \frac12.
\]
Therefore, recalling that $\psi$ is increasing and concave, we obtain that
\[
D^2(\psi \circ \rho_\epsilon)(x-y+2\tau(\xi_j \cdot z)\xi_j) \xi_j \cdot \xi_j \leq   32\psi'(\rho_{\epsilon,\tau})\rho_{\epsilon,\tau} + \frac12\,\psi''(\rho_{\epsilon,\tau}) \qquad \text{if \ $2|r_j+2\tau\xi_j\cdot z| \geq r_j \geq \sqrt{\epsilon}$,}
\]
where we may assume that the right-hand side is negative (so that \eqref{eq_pr1trneg} holds) provided that we choose $\delta$ small enough. On the other hand, following \Cref{rmk_structtest} we also have
\[
|D(\psi\circ\rho_\epsilon)(x-y) \cdot \xi_j | = \frac{r_j}{\sqrt{r_j^2+\epsilon^2}}\,\psi'(\rho_\epsilon) \leq \psi'(\rho_\epsilon).
\]

From the above considerations (which also cover the local diffusion, if there is any, since we can let $\tau = 0$) we deduce the following implication: if
\[
R_j \leq \frac{r_j}4 \,\1_{r_j \geq \sqrt\epsilon} \quad \forall\, j \in I,
\]
then
\begin{equation} \label{eq_contreqlns'} 
\sum_{j \in I} \frac{\scr I_j  + \tilde C' R_j^{\beta_j} \psi'(\rho_\epsilon)}{s_{j}(t)}  
+ \tilde C \psi'(\rho_\epsilon) H_\M(t,\var)\bigr|^x_y \cdot D_x\rho_\epsilon 
 \geq 0,
\end{equation}
where $\tilde C'$ only depends on $\bar c$ and $\tilde C$, while, for the sake of a more compact notation, we set $\beta_0 \defeq 0$ (with the convention that $0^0 \defeq 0$) and 
\begin{equation} \label{eq_prdefI_j}
\scr I_j \defeq \begin{dcases}
 \int_0^1 (1-\tau) \int_{\call P_{V_j'} B_{R_j}} |\xi_j \cdot z|^2 f(\rho_{\epsilon,\tau}) \bar\nu_{j}(\di z)\,\di\tau & j \neq 0
 \\
f(\rho_\epsilon) & j = 0
\end{dcases}
\end{equation}
with
\begin{equation} \label{eq_pmtfdef}
f(\varrho) \defeq 64\psi'(\varrho)\varrho + \psi''(\varrho).
\end{equation}
Furthermore, we can also replace $\call P_{V_j'} B_{R_j}$ in \eqref{eq_prdefI_j} with any subset of it, because $f$ is non-positive by construction.

At this point, we distinguish two cases, where $c_0 < 1$ is a constant to be determined.

$\bullet$ \emph{Case 1: there is $h \in I$ such that $r_h \geq \sqrt{\epsilon}$ and $r_{h} \geq c_0\rho$.} First, suppose that $h \neq 0$. By the previous observations, we can consider \eqref{eq_contreqlns'} with $R_{j} = \frac{r_h}4 \1_{j = h}$, and
\[
B_{\frac{r_h}4}^- \defeq \{ z \in B_{\frac{r_h}4} :\ \xi_{h} \cdot z \leq 0\}
\]
in place of $B_{R_{j}}$ in \eqref{eq_prdefI_j}.
Using \eqref{eq_deftildeH}, \eqref{eq_condnormSH} and that $|\var|_\epsilon$ and $\scr N_\epsilon$ have bounded gradients, 
\begin{equation} \label{eq_estdriftdp} \begin{split}
s_{h}(t) \bigl|H_\M(t,\var)\bigr|^x_y \cdot D_x \rho_\epsilon\bigr| &\leq K_H(1+C_{\scr N}) |x-y| \leq K (1+\kappa C_{\scr N}) (\kappa c_{\scr N})^{-1}(\rho_\epsilon + 2\epsilon).
\end{split}
\end{equation}
Consequently, we can exploit the bound $\sqrt\epsilon \leq 4(r-\epsilon) \leq 4\rho_\epsilon$ to obtain
\begin{equation} \label{eq_esttildeH}
s_{h}(t) \bigl|H_\M(t,\var)\bigr|^x_y \cdot D_x \rho_\epsilon\bigr| \leq \tilde{K}\rho_\epsilon,
\end{equation}
where $\tilde{K}$ only depends on $K_H$, $C_{\scr N}$, $c_{\scr N}$, and $\kappa$. 
On the other hand, by the radial symmetry of $\bar\nu_h$ and \eqref{eq_condnurad1},
\be\label{coerc-h}
\begin{split}
\int_{{\call P}_{V_h'} B_{\frac{r_h}4}^-} |\xi_{h} \cdot z|^2 \,\bar\nu_{h}(\di z) = \frac1{2d_{h}} \int_{{\call P}_{V_h'} B_{\frac{r_h}4}} |\var|^2 \,\di\bar\nu_{h} \geq \frac{\underline c}{2d_{h}} \Bigl(\frac{r_h}4 \Bigr)^{\alpha_{h}} \geq \frac{\underline c}{2d_{h}} \Bigl( \frac{c_0\rho_\epsilon}{4} \Bigr)^{\alpha_{h}},
\end{split}
\ee
so from \eqref{eq_esttildeH} we obtain
\begin{equation} \label{eq_p1ltn1}
s_{h}(t) \bigl|H_\M(t,\var)\bigr|^x_y \cdot D_x \rho_\epsilon\bigr| \leq \frac{2^{1+2\alpha_h}d_h\tilde K}{\underline cc_0^{\alpha_h}} \rho_\epsilon^{1-\alpha_h} \int_{{\call P}_{V_h'} B_{\frac{r_h}4}^-} |\xi_{h} \cdot z|^2 \,\bar\nu_{h}(\di z),
\end{equation}
while, since $\beta_h \leq \alpha_h$ and $r_h \geq c_0 \rho\geq c_0 \rho_\epsilon$,
\begin{equation} \label{eq_p1ltn2}
\Bigl(\frac{r_h}4\Bigr)^{\beta_h} 
\leq \frac{2^{1+2(\alpha_h-\beta_h)}d_h}{\underline cc_0^{\alpha_h-\beta_h}} \rho_\epsilon^{\beta_h - \alpha_h} \int_{{\call P}_{V_h'} B_{\frac{r_h}4}^-} |\xi_{h} \cdot z|^2 \,\bar\nu_{h}(\di z) .
\end{equation}
Combining \eqref{eq_prdefI_j} (where $R_{j} = \frac{r_h}4 \1_{j = h}$), \eqref{eq_p1ltn1} and \eqref{eq_p1ltn2}, we have
\begin{equation} \label{eq_pr1ad1}
\int_0^1 (1-\tau) \int_{\call P_{V_h'} B_{\frac{r_h}4}^-} |\xi_h \cdot z|^2 \Bigl( f(\rho_{\epsilon,\tau}) + \tilde K' \psi'(\rho_\epsilon) \rho_\epsilon^{\alpha'} \Bigr) \bar\nu_{h}(\di z)\,\di\tau \geq 0,
\end{equation}
where $\alpha' \defeq \min_{j \in I} ( ( \beta_{j} - \alpha_{j} ) \wedge ( 1 - \alpha_{j} )^- ) \in (-1,0]$, and $\tilde K'>0$ is a constant depending only on $N$, $\tilde K$, $\underline c$, $\bar c$, $c_0$, and $\tilde C$. Furthermore, since $\psi$ is concave and  $\rho_{\epsilon,\tau} \in (0,\rho_\epsilon)$ in $\call P_{V_h'} B_{\frac{r_h}4}^-$,
 \eqref{eq_pr1ad1} implies
\[
\int_0^1 (1-\tau) \int_{\call P_{V_h'} B_{\frac{r_h}4}^-} |\xi_h \cdot z|^2 \Bigl( f(\rho_{\epsilon,\tau}) + \tilde K' \psi'(\rho_{\epsilon,\tau}) \rho_{\epsilon,\tau}^{\alpha'} \Bigr) \bar\nu_{h}(\di z)\,\di\tau \geq 0,
\]
so we see that we obtain a contradiction by choosing $\vartheta < 1+\alpha'$ (in the definition of $\psi$) and $\delta$ sufficiently small (depending only on $\tilde{\call K}'$, $\alpha'$, and $\vartheta$). The case if $h=0$ is similar and easier, so we omit the details; we only note that it suffices to take $R_j = \1_{j = 0}$ and use \eqref{eq_esttildeH} alongside the definition of $\scr I_0$.

$\bullet$ \emph{Case 2: for all $j \in I$, $r_j < c_0\rho$ or $r_j < \sqrt\epsilon$.} Take $R_j=0$ for all $j$. Letting ${\call P}_{\call V}^\perp(x-y) = R\xi^\perp$, with $R \geq 0$ and $\xi^\perp \in \bb S^{N-1}$, we have
\[
D_x\rho_\epsilon = \sum_{j \in I} \frac{r_j}{\sqrt{r_j^2+\epsilon^2}}\,\xi_j + \kappa {\call P}_{\call V}^\perp D\scr N_\epsilon(R\xi^\perp),
\]
so, using \rife{ass_Neps}--\rife{eq_condcallN'} along with \eqref{eq_condcallN''}, we obtain
\[
\begin{split}
&-  H_\M(t,\var)\bigr|^x_y \cdot D_x\rho_\epsilon 
\\
&= - \int_0^1  DH_\M(t,y+\tau(x-y)) \Bigl(\, \sum_{j \in I} r_j\xi_j+R\xi^\perp\Bigr) \cdot \Bigl(\, \sum_{j \in I}\frac{r_j}{\sqrt{r_j^2+\epsilon^2}}\,\xi_j + \kappa {\call P}_{\call V}^\perp D\scr N_\epsilon(R\xi^\perp) \Bigr)\,\di\tau
\\
&\geq  \kappa \gamma \scr N(R\xi^\perp)\, \bigl( \|D H_\M(t,\var)\|_\infty \vee 1\bigr) - \kappa \gamma'\epsilon -  \|D H_\M(t,\var)\|_\infty \Bigl( (|I|+\kappa C_{\scr N}) \sum_{j \in I} r_j + |I|R \Bigr) 
\\
&\geq \bigl( \|D H_\M(t,\var)\|_\infty \vee 1\bigr) \biggl( \Bigl( \kappa \gamma - \frac{|I|}{c_{\scr N}} \Bigr) \scr N(R\xi^\perp) -  (|I|+\kappa C_{\scr N}) \sum_{j \in I} r_j \biggr) 
- \kappa \gamma' \epsilon
\end{split}
\]
Hence, as $\scr N(R\xi^\perp) = \rho - \sum_{j \in I} r_j \geq (1-|I|c_0)\rho - |I|\sqrt\epsilon$, it suffices to choose $\kappa$ large enough (depending only on  $c_{\scr N}$, $\gamma$, and $|I|$) and $c_0$ small enough (depending only on  $C_{\scr N}$, $|I|$, and $\kappa$) to obtain that
\begin{equation} \label{eq_estcontrcase2}
-H_\M(t,\var)\bigr|^x_y \cdot D_x\rho_\epsilon \geq c\rho - c'\sqrt\epsilon
\end{equation}
for some $c,c'>0$ that only depend on the constants above. Therefore, recalling \eqref{eq_contreqlns'} we obtain $c\rho \leq c'\sqrt\epsilon$, which contradicts that $\call U$ has a positive interior maximum as we let $\epsilon \dto 0$.

For the sake of completeness, we also mention that if $\call V^\perp = \{0\}$, this second case reduces to just having $\rho \lesssim \sqrt\epsilon$ by assumption, and we again obtain the previous contradiction.
\end{proof}

We now come to the proof of the  long-time exponential decay of the oscillation.

\begin{proof}[\textbf{Proof of Proposition \ref{prop_oscdec}}] 
This proof is a combination of those of \Cref{thm_mainsmnl} and  \cite[Theorem~3.2]{ForPor}, with some adaptations.
Suppose that $\theta \in \bigl(0,1 \wedge (2-\max_{j \in I}\alpha_j)\bigr)$, and let $\psi(\varrho) \defeq 1 - e^{-C_1\varrho^\theta}$, for $C_1>1$ to be determined.   With $v \defeq u_\M$ and $\tilde\phi \defeq \phi_\M$, set $\call U_\epsilon \defeq v(t,x) - v(t,y) - \zeta_\epsilon(t,x,y)$ where
\[
\zeta_\epsilon(t,x,y) \defeq K_0 e^{-\varpi t} \bigl(\tilde\phi(x)+\tilde\phi(y)+L\bigr)\psi(\rho_\epsilon) + \frac\eta{T-t}, 
\]
$L>0$ is to be determined, $K_0$ satisfies \eqref{eq_choiceK0NL''} with $\delta = +\infty$, and $\rho_\epsilon$ is as in the proof of \Cref{thm_mainsmnl}. Supposing that $\call U_0$ has a positive supremum in $(0,T) \times \R^{2N}$, by \Cref{lem_exmax}, $\call U_\epsilon$ has a positive interior maximum for any sufficiently small $\epsilon>0$. 
With the same notation and arguments as in the proof of \Cref{thm_mainsmnl}, we obtain 
\begin{align} \label{eq_nov}
&   \bigl( (\de_t - \call L_{\call R}^H)\tilde\phi(t,x) + (\de_t - \call L_{\call R}^H)\tilde\phi(t,y) - \varpi (\tilde\phi(x)+\tilde\phi(y)+L)\bigr) \psi(\rho_\epsilon) \notag
\\ 
& \leq  \bigl(\tilde\phi(t,x)+\tilde\phi(t,y)+L\bigr)\psi'(\rho_\epsilon) H_\M(t,\var)\bigr|^x_y \cdot D_x\rho_\epsilon
\notag \\
&\quad + \frac{R_0}{s_0(t)} \Bigl( \bigl(\tilde\phi(t,x) + \tilde\phi(t,y) \bigr) D^2(\psi \circ \rho_\epsilon)(x - y) \xi_0 \cdot \xi_0 + C_{D\tilde\phi(t,\var),0}  |D(\psi \circ \rho_\epsilon)( x- y) \cdot \xi_0| \Bigr) \1_{\call P_{W_0}(x - y) \neq 0 \wedge 0 \in I}
\notag \\
&\quad + \sum_{j \in I \setminus \{0\}} \frac{2}{s_{j}(t)} \biggl(  \int_{\call P_{V_j'} B_{R_j}} \bar{\tilde\phi}(t,z) |\xi_j \cdot z|^2 \int_0^1 (1-\tau) D^2(\psi\circ\rho_\epsilon)( x- y+2\tau(\xi_j\cdot z)\xi_j ) \xi_j \cdot \xi_j\, \di\tau\,\bar\nu_{j}(\di z)
\notag \\
&\quad + \Bigl( \int_{\call P_{V_j'} B_{R_j}} |\var|^{2}\,\di \bar\nu_j \Bigr) C_{D\tilde\phi(t,\var),R_j}  |D(\psi \circ \rho_\epsilon)( x- y) \cdot \xi_j| \biggl) \1_{\call P_{W_j}(x-y) \neq 0},
\end{align}
where $\xi_j \defeq \wwidehat{\call P_{W_j}(x-y)}$.

By assumption, there is $R>1$ (depending on $\phi$) such that $\call L_{\call R}^H\tilde\phi(t,x) \geq \frac{\varpi_0}2 \tilde \phi(x)$ whenever $|x| > R$.

Suppose that, at the maximum of $\call U_\epsilon$, $\rho_\epsilon > R' \defeq 2C_{\scr N}R$ (by taking $R$ larger if needed, we can assume that $R' \geq 1$); then by \eqref{eq_condcallN'} and the triangle inequality we can assume without loss of generality that (at the maximum) $|x| > R$. In this case, using 
\eqref{eq_nov}
with $R_j = 0$ for all $j$, we have 
\begin{equation} \label{eq_HLosc1}
\Bigl( \frac{\varpi_0}2 \tilde\phi(x) - \call L^H_\M \tilde\phi(t,y) - \varpi (\tilde\phi(x)+\tilde\phi(y)) \Bigr)\psi(\rho_\epsilon) - \tilde{K} \bigl(\tilde\phi(x) + \tilde\phi(y) + L\bigr) \psi'(\rho_\epsilon) \rho_\epsilon - L\varpi\psi(\rho_\epsilon) \leq 0,
\end{equation}
where $\tilde K$ is such that $|H_\M(t,\var)|^x_y \cdot D_x \rho_\epsilon| \leq \tilde{K}\rho_\epsilon$, which is given by the Lipschitz bound on $H$. 
By the definition of $\psi$ we have $\psi(\rho_\epsilon) \geq \frac12$ (because $R' \wedge C_1 > 1$) and $\varrho \mapsto \psi'(\varrho)\varrho$ monotonically decreases to $0$ as $\smash{\varrho \in (C_1^{-1/\theta},+\infty) \supset (R',+\infty)}$; hence, $\psi'(\rho_\epsilon)\rho_\epsilon < 2\psi(\rho_\epsilon)\psi'(R')R'$, and by taking $R'$ larger (that is, $R$ larger) if needed, we can suppose that 
\begin{equation} \label{eq_choiceR'}
\psi'(R')R' < \frac{\varpi_0}{8\tilde{K}}.
\end{equation}
From \eqref{eq_HLosc1} we obtain
\begin{equation} \label{eq_HLosc2}
\begin{split}
L\bigl(2\tilde{K}\psi'(R')R'+\varpi) &\geq \Bigl(\frac{\varpi_0}4 - \varpi\Bigr) \tilde\phi(x) - \call L^H_\M \tilde\phi(t,y) - \Bigl( \varpi + \frac{\varpi_0}4 \Bigr) \tilde\phi(y) \geq \frac{\varpi_0}8 \phi(x) - \tilde k_0
\end{split}
\end{equation}
where the last inequality holds if $\varpi \leq \frac{\varpi_0}8$, and we also used that, due to \eqref{eq_condLyapdec}, $\tilde k_0 \defeq \sup (\call L\phi + (\varpi + \frac{\varpi_0}4 ) \phi )$ is finite because $\varpi + \frac{\varpi_0}4 < \varpi_0$. By taking $R'$ larger if needed, we may assume that $\phi(x) \geq {16\tilde k_0}/\varpi_0$; therefore, we obtain a contradiction, provided that we choose
\begin{equation} \label{eq_choiceL}
L = \frac{\varpi_0 c_{R'}}{32\bigl(\tilde{K}\psi'(R')R'+\varpi)}, \qquad c_{R'} \defeq \inf_{|x| > R'} \phi(x).
\end{equation}

If, instead, $\rho_\epsilon \leq R'$, then we can proceed similarly as in the proof of \Cref{thm_mainsmnl}.

$\bullet$ \emph{Case 1: there is $h \in I$ such that $r_h \geq \sqrt{\epsilon}$ and $r_{h} \geq c_0\rho$ (with $\rho_\epsilon \leq R'$).} 

Since $\varpi<\varpi_0$, there exists $k_0 > 0$ such that $ (\de_t - \call L_{\call R}^H)\tilde\phi- \varpi \tilde \phi \geq -k_0$.   Hence,   we deduce
\begin{align} \label{eq_nov1'}
& \bigl( (\de_t - \call L_{\call R}^H)\tilde\phi(t,x) + (\de_t - \call L_{\call R}^H)\tilde\phi(t,y) - \varpi (\tilde\phi(x)+\tilde\phi(y)+L)\bigr) \psi(\rho_\epsilon) \notag\\
&-  \bigl(\tilde\phi(t,x)+\tilde\phi(t,y)+L\bigr)\psi'(\rho_\epsilon) H_\M(t,\var)\bigr|^x_y \cdot D_x\rho_\epsilon \notag\\
&\geq  -\bigl(  2k_0 + L \varpi\bigr) \psi(\rho_\epsilon) - \bigl(\tilde\phi(t,x)+\tilde\phi(t,y)+L\bigr) \tilde K \, \psi'(\rho_\epsilon)\rho_\epsilon, 
\end{align}
with $\tilde K$ as above.
Note that $\varpi \leq \frac{\varpi_0}8$ and \rife{eq_choiceR'}, together with \rife{eq_choiceL}, imply $L\geq \frac{c_{R'}}8$. Then, using the monotonicity of $\psi(\varrho)$ for $C_1\varrho^\theta>1$, we estimate, for $\varrho \leq R'$:
\begin{equation} \label{eq_estpsiosc} \begin{split}\psi(\varrho) &\leq \frac{e}\theta \psi'(\varrho)\varrho \1_{C_1\varrho^\theta \leq 1} + \psi(\varrho) L \frac{32\bigl(\tilde{K}\psi'(\varrho)\varrho+\varpi)}{\varpi_0 c_{R'}}\1_{C_1\varrho^\theta > 1} \leq \frac{\ell L}{c_{R'}}\bigl(\psi'(\varrho)\varrho + \varpi\psi(\varrho)\bigr) ,
\end{split}\end{equation}
where $\ell$ only depends on $\theta$, $\varpi_0$, and $\tilde{K}$. We can also suppose that $R$ is so large that
\begin{equation} \label{eq_choicephiR}
\frac{k_0\ell}{c_{R'}} < \frac c4,
\end{equation}
where $c$ is the constant appearing in \eqref{eq_estcontrcase2}. Now, we deduce from \eqref{eq_nov} and \eqref{eq_nov1'} that
\begin{align} \label{eq_nov1}
&   - L \varpi \Bigl(  1+ \frac c2 \Bigr) \psi(\rho_\epsilon)
   \notag   
   \\
   &\leq  \bigl(\tilde\phi(t,x)+\tilde\phi(t,y)+L\bigr) \Bigl(\tilde K + \frac c2 \Bigr)\, \psi'(\rho_\epsilon)\rho_\epsilon 
\notag \\
&\quad + \frac{R_0}{s_0(t)} \Bigl( \bigl(\tilde\phi(t,x) + \tilde\phi(t,y) \bigr) D^2(\psi \circ \rho_\epsilon)(x - y) \xi_0 \cdot \xi_0 + C_{D\tilde\phi(t,\var),0}  |D(\psi \circ \rho_\epsilon)( x- y) \cdot \xi_0| \Bigr) \1_{\call P_{W_0}(x - y) \neq 0 \wedge 0 \in I}
\notag \\
&\quad + \sum_{j \in I \setminus \{0\}} \frac{2}{s_{j}(t)} \biggl(  \int_{\call P_{V_j'} B_{R_j}} \bar{\tilde\phi}(t,z) |\xi_j \cdot z|^2 \int_0^1 (1-\tau) D^2(\psi\circ\rho_\epsilon)( x- y+2\tau(\xi_j\cdot z)\xi_j ) \xi_j \cdot \xi_j\, \di\tau\,\bar\nu_{j}(\di z)
\notag \\
&\quad + \Bigl( \int_{\call P_{V_j'} B_{R_j}} |\var|^{2}\,\di \bar\nu_j \Bigr) C_{D\tilde\phi(t,\var),R_j}  |D(\psi \circ \rho_\epsilon)( x- y) \cdot \xi_j| \biggl) \1_{\call P_{W_j}(x-y) \neq 0},
\end{align}
Following now the argument in Case 1 in the proof of \Cref{thm_mainsmnl}, we obtain  
\begin{equation} \label{eq_couplosc}
\int_0^1 (1-\tau) \int_{\call P_{V_h'} B_{\frac{r_h}4}^-} |\xi_h \cdot z|^2 \Bigl( f(\rho_{\epsilon,\tau}) + \tilde K' \psi'(\rho_{\epsilon,\tau}) \rho_{\epsilon,\tau}^{\alpha'} \Bigr) \bar\nu_{h}(\di z)\,\di\tau + \varpi \Bigl( 1 + \frac c2 \Bigr) \psi(\rho_\epsilon) \geq 0,
\end{equation}
Direct computations show that we can choose $C_1$ large enough (depending only on the constants in \eqref{eq_couplosc}) to have
\[
 f(\rho_{\epsilon,\tau}) + \tilde K' \psi'(\rho_{\epsilon,\tau}) \rho_{\epsilon,\tau}^{\alpha'} \leq - \tilde c \rho_{\epsilon,\tau}^{\theta-2} e^{-C_1\rho_{\epsilon,\tau}^\theta} \leq - \tilde c \rho_{\epsilon}^{\theta-2} e^{-C_1\rho_{\epsilon}^\theta},
\]
with $\tilde c$ depending on $\theta$ and $C_1$; therefore, recalling \eqref{coerc-h}, \eqref{eq_couplosc} implies that  
\begin{equation} \label{eq_lhsofw}
\frac{\tilde c \underline c c_0^{\alpha_h}}{4^{1+\alpha_h}d_{h}} e^{-C_1\rho_{\epsilon}^\theta} \rho_{\epsilon}^{\theta-2+\alpha_h} \leq \varpi\Bigl(1 + \frac c2 \Bigr).
\end{equation}
This gives a contradiction by choosing $\varpi$ sufficiently small, depending on the constants above and $R'$, provided that $\theta < 2-\max_{j \in I}\alpha_j$, because in this case the function on the left-hand side of \eqref{eq_lhsofw} is decreasing in $\rho_\epsilon$.

\noindent$\bullet$ \emph{Case 2: for all $j \in \{0,\dots,k\}$, $r_j < c_0\rho$ or $r_j < \sqrt\epsilon$ (with $\rho_\epsilon \leq R'$).} We use \eqref{eq_nov} with $R_j=0$ for all $j$. Estimate \eqref{eq_estcontrcase2} holds by the same computations as in \Cref{thm_mainsmnl}; then using also \eqref{eq_estpsiosc} and \eqref{eq_choicephiR} as before, we obtain
\[
\psi'(\rho_\epsilon) \bigl( c\rho_\epsilon - c'\sqrt\epsilon \bigr) \leq \frac c2\psi'(\rho_\epsilon)\rho_\epsilon + \varpi \Bigl( 1 + \frac c2 \Bigr)\psi(\rho_\epsilon).
\]
Choosing ${\varpi \leq \frac{c}{4+2c} \inf_{\varrho \in (0,R')} \frac{\psi'(\varrho)\varrho}{\psi(\varrho)}}$, we obtain $\frac c4 \rho_\epsilon \leq c'\sqrt\epsilon$, and the contradictions follows by letting $\epsilon \dto 0$.
This concludes the proof for $\theta \in \bigl(0,1 \wedge (2-\max_{j \in I}\alpha_j)\bigr)$. 

Consider now the case that $\theta \in \bigl[1 \wedge (2-\max_{\sigma \in\Sigma'}\alpha_\sigma), 1\bigr]$. As we proved that \eqref{eq_gthm_uosc} holds for some $\theta' < \theta$, it suffices to prove that there is $K_0 \leq  C_\theta' [u(0,\var)]_{\phi,\theta}$ (with $C_\theta'$ that only depends on $d$, $|I|$, $\norm\M$, $c_{\scr N}$, $C_{\scr N}$, $K$, $\gamma$, $\phi$, and $\theta$) such that
\[
|u(t,x) - u(t,y)| \leq K_0 e^{-\varpi t} \bigl(\phi(x) + \phi(y)\bigr) \rho^\theta \qquad \text{for } \rho < \delta,
\]
where $\delta$ can be fixed arbitrarily small. This can be done by choosing $\psi(\varrho) \defeq \varrho^\theta - \frac1{2(1+\theta')} \varrho^{1+\theta'} \1_{\theta = 1}$, and following the previous argument, with $L=0$ and $\delta$ in place of $R'$ (or, equivalently, the argument in the proof of \Cref{thm_mainsmnl} with $\phi(t,x)$ replaced by $e^{-\varpi t}\phi(x)$). Note that we obtain from  \Cref{lem_exmax} that we can choose  
\begin{equation*} \begin{split}
K_0 &= \frac{\sup_{t\in(0,T)} [u(t,\var)]_{\phi,0}}{\psi(\delta/2)}\, \vee \frac{2[u(0,\var)]_{\phi,\theta}}{(c_{\scr N} \wedge 1)c_{\psi,\theta,\delta}} \leq \frac{\sup_{t\in(0,T)} [u(t,\var)]_{\phi,\theta'}}{\psi(\delta/2)}\, \vee \frac{2[u(0,\var)]_{\phi,\theta}}{(c_{\scr N} \wedge 1)c_{\psi,\theta,\delta}}
\\
&\leq \frac{C_{\theta'} [u(0,\var)]_{\phi,\theta'}}{\psi(\delta/2)}\, \vee \frac{2[u(0,\var)]_{\phi,\theta}}{(c_{\scr N} \wedge 1)c_{\psi,\theta,\delta}} \leq \Bigl( \frac{C_{\theta'}}{\psi(\delta/2)}\, \vee \frac{2}{(c_{\scr N} \wedge 1)c_{\psi,\theta,\delta}} \Bigr) [u(0,\var)]_{\phi,\theta}.
\end{split}
\end{equation*}
This concludes the whole proof.
\end{proof}

\section{Main applications: H\"ormander hierarchies} \label{sec:ex}

In this section, we specialize our results to the models presented in the Introduction, which provide the main explicit applications of \Cref{thm_mainsmnl,prop_oscdec}. 
We start with the proof of Theorem \ref{intro-thm1}, which describes the cascade dynamical system 
\begin{equation}\label{cascade}
\begin{dcases}
\di X^0_t = H^0(t,X^0_t,\dots,X^r_t)\di t + \di Z_t
\\
\di X_t^1 = H^1(t,X^0_t,\dots,X^r_t)\di t\\
\di X_t^2 = H^2(t,X^1_t,\dots,X^r_t)\di t \\[2pt]
 \ \ \smash{\vdots} \\
\di X_t^r = H^r(t,X^{r-1}_t,X^r_t)\di t;
\end{dcases}
 \end{equation}
here the infinitesimal generator of $Z_t$ is what we called a standard diffusion operator acting on the variable $x_0$ (cf. \cref{def_opintro}).

\vskip1em

\begin{proof}[{\bfseries Proof of \Cref{intro-thm1}}]
According to the general setting in \Cref{sec:genset} and on account of \Cref{rmk_extloc}, we can suppose that $\ell =1$, $V_1=\R^{d_0} \times \{0\}$, and $\bar\nu_1 = \lim_{\sigma\uto \bar\sigma_\cL} (2-\sigma)|\var|^{-d_0-\sigma} \scr H^{d_0}$.
We first consider the case where $d_i = d$ for all $i$; here we can suppose that $J_i = I_d$ in \eqref{eq_HhormC1intro}. 

\emph{Part (1).} In order to apply \Cref{thm_mainsmnl},\footnote{Let us point out that \Cref{thm_mainsmnl} will apply   here to bounded solutions,  in which case an unbounded supersolution $\Phi$ satisfying \rife{ultima_Phi} always exists for globally Lipschitz drift functions $H$ and standard diffusive operator $\cL$; see \Cref{rmk_onuLip}.}  we need to show that there are $\{W_j\}_{j\in I}$, $\M \in C^0((0,T];\GLN) \cap C^1((0,T);\GLN)$ with $\M(0^+) = 0$, satisfying \eqref{eq_defSt}--\rife{eq_boundTnus}, and $\scr N_\epsilon$ satisfying \rife{ass_Neps}--\rife{eq_condcallN'}, such that \eqref{eq_sbspdec}, \eqref{eq_condnormSH}, and \eqref{eq_condcallN''} hold. 
We make the ansatz that $\M(t) = \call S(t)\call T$ for some $\call T \in \GLN$, and with $\call S(t)$ given by \eqref{eq_scalintro}, i.e.\ $\call S(t) = t^{\bar\sigma^{-1}} \mathrm{diag}(t^{i})_{0\leq i\leq r} \otimes I_d$, where we let $\bar\sigma = \bar\sigma_\cL$ for the sake of a lighter notation. For some $\mu> 1$ and $\alpha_i,\beta_i > 0$ to be determined, define
\begin{equation} \label{eq_defrotcoord}
\call T \defeq \begin{pmatrix}
1 	& -\mu \beta_1	&0 & \dots 	& 0 \\
0		& 1				& \!\!\!-\mu\beta_2 &\dots 	& 0 \\
\vdots 	& \vdots					& \ddots	& \ddots	 & \vdots \\
0 & 0 & \dots & 1 & - \mu \beta_r \\
0		& 0						& \dots	& 0	& 1
\end{pmatrix}
\otimes I_d
\end{equation}
and, for all $\epsilon \in (0,1]$,
\begin{equation} \label{eq_chNHorm}
\scr N_\epsilon(z) \defeq \sum_{0 \leq j \leq r} \alpha_j\mu^j|z^j|_\epsilon, \qquad z \in \R^N.
\end{equation}
where $|\var|_\epsilon = \sqrt{|\var|^2 + \epsilon^2} - \epsilon$. Clearly, $\scr N_\epsilon$  satisfies \rife{ass_Neps}--\rife{eq_condcallN'}, with constants only depending on $\alpha_j$, $\mu$, and $r$.
On the other hand, note that $\call T$ is the identity on $\R^{d_0} \times \{0\}$, so we have $I = \{1\}$ and $W_1 = V_1' = \R^{d_0} \times \{0\}$.   
Then, by the assumptions on $\cL$ (cf. \rife{nulb_intro}) and on account of Remark \ref{rmk_nu0hom},  $\M(t) = \call S(t)\call T$ satisfies \rife{eq_boundTnus} with $s_1(t) = t$.
In addition, we have
\[
\call S^{-1} \bigl( \dot{\call S} + DH\call S \bigr) = \frac1t \Bigl( \mathrm{diag}(i+\bar\sigma^{-1})_{0 \leq i \leq r} \otimes I_d + \bigl( t^{j-i+1} D_j H^i \bigr)_{0 \leq i,j \leq r} \Bigr) 
\]
so that, since $D_j H^i \neq 0$ only for $j\geq i-1$, we deduce that 
\[
s_1(t) \sup_{x \in \R^N} \lVert DH_\M(t,x)\rVert \leq K,
\]
with $K$ depending on $\bar\sigma$, $\norm{DH}_\infty$, $r$ and $T$, as well on the choice of  $\mu$, $\alpha_i$, and $\beta_i$; this shows that \eqref{eq_condnormSH} holds. 

We are left to check condition \rife{eq_condcallN''}. To this purpose,
we note that\footnote{We use the usual convention for products that $\prod_{k \in \emptyset} \beta_k \defeq 1$.}
\begin{gather*}
\call T^{-1} = \Bigl( \mu^{j-i} \prod_{k={i+1}}^{j} \beta_k  \1_{j\geq i} \Bigr)_{0 \leq i,j \leq r},
\\
D\scr N_\epsilon(z) = \call B_\epsilon(z)z, \qquad \call B_\epsilon(z) \defeq \mathrm{diag}\Bigl( \frac{\alpha_j\mu^j}{\sqrt{|z^j|^2 + \epsilon^2}}\Bigr)_{0 \leq j \leq r} \otimes I_d;
\end{gather*}
therefore, since ${\call P_{\call V}} = \call P_{V_1} = \mathrm{diag}(1,0,\dots,0) \otimes I_d$,
\[
 \call B_\epsilon(z)^* \call P_{\call V}^\perp \call T^{-1} = \Bigl( \frac{\alpha_i\mu^j}{\sqrt{|z^i|^2 + \epsilon^2}} \prod_{k=i+1}^{j} \beta_k   \1_{j\geq i \neq 0} \Bigr)_{0 \leq i,j \leq r} \otimes I_d,
\]
and then
\begin{gather*}
t \call B_\epsilon(z)^* \call P_{\call V}^\perp \M^{-1} \dot{\M} = \biggl( \mathrm{diag} \Bigl( \frac{(j+\bar\sigma^{-1})\alpha_j\mu^j}{\sqrt{|z^j|^2 + \epsilon^2}}  \1_{j \neq 0} \Bigr)_{0 \leq j \leq r} + \Bigl( \frac{\alpha_i\mu^j}{\sqrt{|z^i|^2 + \epsilon^2}} \prod_{k=i+1}^{j} \beta_k   \1_{j>i \neq 0} \Bigr)_{0 \leq i,j \leq r} \biggr) \otimes I_d ,
\\
t \call B_\epsilon(z)^* \call P_{\call V}^\perp \M^{-1} DH \M =  \Bigl(\sum_{h=i\vee1}^{r} \frac{\alpha_i\mu^h}{\sqrt{|z^i|^2 + \epsilon^2}} \prod_{k=i+1}^{h} \beta_k \bigl( t^{j-h+1} D_j H^h - \mu \beta_j t^{j-h} D_{j-1} H^h \bigr) \Bigr)_{0 \leq i,j \leq r}.
\end{gather*}
Since $|z^i|/\sqrt{|z^i|^2 +\epsilon^2} \leq 1$, we have, for all $z \in \ker \call P_{\call V}$,
\[
t \call B_\epsilon(z)^* \call P_{\call V}^\perp \M^{-1} \dot{\M} z \cdot z \leq \sum_{j=1}^r \mu^j \Bigl( (j+\bar\sigma^{-1})\alpha_j + \sum_{i=1}^{j-1} \alpha_i \prod_{k=i+1}^j\beta_k \Bigr) |z^j|;
\]
then, recalling \eqref{eq_HhormC1intro} and letting $\gamma_j>0$ such that $D_{j-1}H^j \geq \gamma_j I_d$, we have
\begin{equation*} \begin{split}
&t \call B_\epsilon(z)^* \call P_{\call V}^\perp \M^{-1} DH \M z \cdot z
\\
&= \sum_{j =1}^r \frac{\alpha_j\mu^{j}}{\sqrt{|z^j|^2 + \epsilon^2}} t D_j H^j z^j \cdot z^j  - \sum_{j = 1}^r \frac{\alpha_j\mu^{j+1}\beta_j}{\sqrt{|z^j|^2 + \epsilon^2}} D_{j-1} H^j z^j \cdot z^j + \sum_{j=1}^r \frac{\alpha_j\mu^{j+1}}{\sqrt{|z^j|^2 + \epsilon^2}} \beta_{j+1} D_j H^{j+1} z^j \cdot z^j 
\\
&\quad + \sum_{\substack{i,j=1\\j\neq i}}^r \sum_{h=i}^{j+1} \frac{\alpha_i\mu^{h}}{\sqrt{|z^i|^2 + \epsilon^2}} \prod_{k=i+1}^h \beta_k t^{j-h+1} D_j H^h z^i \cdot z^j - \sum_{\substack{i,j=1\\j\neq i}}^r \sum_{h=i}^{j} \frac{\alpha_i\mu^{h+1}}{\sqrt{|z^i|^2 + \epsilon^2}} \beta_j \prod_{k=i+1}^h \beta_k t^{j-h} D_{j-1} H^h z^i \cdot z^j
\\
&\leq C \sum_{j =1}^r \biggl( {\alpha_j\mu^{j}} ( 1 +\mu \beta_{j+1}) +  \sum_{i=1}^{j-1} \alpha_i \Bigl(\, \sum_{h=i}^{j+1} \mu^h \prod_{k=i+1}^h \beta_k + \sum_{h=i}^j  \mu^{h+1} \beta_j \prod_{k=i+1}^h \beta_k \Bigr)  + \alpha_{j+1} \mu^{j+1} \biggr) |z^j|
\\
&\quad - \sum_{j = 1}^r \frac{\alpha_j \beta_j\mu^{j+1} \gamma_j |z^j|^2}{\sqrt{|z^j|^2 + \epsilon^2}} ,
\end{split}\end{equation*}
where $C$ only depends on $\norm{D_jH^i}_\infty$ and $T$. Since 
\[
\frac{|z^i|^2}{\sqrt{|z^i|^2 + \epsilon^2}} = \frac{|z^i|_\epsilon(|z^i|_\epsilon + 2\epsilon)}{|z|_\epsilon + \epsilon} \geq |z^i|_\epsilon \geq |z^i| - \epsilon,
\]
we deduce that
\be\label{nov26}
t \call B_\epsilon(z)^* \call P_{\call V}^\perp \M^{-1} \bigl( \dot{\M} + DH \M \bigr) z \cdot z \leq -\sum_{j=1}^r c_j |z^j| + c\,\epsilon \qquad \forall \, z \in \ker {\call P_{\call V}}
\ee
with $c=\sum_{j=1}^r \gamma_j\alpha_j\beta_j \mu^{j+1}$ and 
\begin{equation*}
c_j = \mu^{j+1} \biggl( \gamma_j \alpha_j \beta_j - C \Bigl( \alpha_j \beta_{j+1} + (\beta_{j+1} + \beta_j ) \sum_{i=1}^{j-1} \alpha_i  \prod_{k=i+1}^j \beta_k + \alpha_{j+1} \Bigl) + \frac{C'}\mu \biggr),
\end{equation*}
where $C'$ only depends on $C$, $\alpha_i$, $\beta_i$, and $r$. Now we observe that
\begin{equation*}
\frac{c_j}{\mu^{j+1}} - \frac{C'}\mu = C\alpha_j \biggl( \frac{\gamma_j \beta_j}{C} - \frac{\alpha_{j+1}}{\alpha_j} - \beta_{j+1} - \frac{\beta_j + \beta_{j+1}}{\alpha_j}\sum_{i=1}^{j-1} \alpha_i  \prod_{k=i+1}^j \beta_k  \biggr);
\end{equation*}
therefore, if we recursively choose
\[
\beta_{j+1} = \beta_j^2, \qquad \alpha_j = \alpha_0^j \beta_j,
\]
we have
\[
\frac{c_j}{\mu^{j+1}} - \frac{C'}\mu = C\alpha_j \beta_j \biggl( \frac{\gamma_j}{C} - \alpha_0 - \beta_{j}^2 - (1 + \beta_{j}) \sum_{i=1}^{j-1} \frac1{\alpha_0^{j-i}} \prod_{k=i}^{j-1} \beta_k  \biggr) \geq \frac13 \alpha_j \beta_j \gamma_j,
\]
provided that we take $\alpha_0 = \frac{1}{2C}\min_j \gamma_j$, and $\beta_1$ small enough (depending on $\gamma_j$, $\alpha_0$, and $r$). At this point, it suffices to take $\mu$ large enough (depending on $\alpha_0$, $\beta_1$, $\gamma_j$, and $C'$) to have
\[
c_j \geq c'\alpha_j \mu^j,
\]
with $c'$ depending on $\mu$, $\beta_1$, and $\min_j \gamma_j$; due to \rife{nov26} and to the definition of $\call N_0$, this proves that \eqref{eq_condcallN''} holds with $\gamma =c'$.

\emph{Part (2)}. By the previous computations, it is enough to take $\M = \call S(1)\call T$ with $\call S(t)$ and $\call T$ as above, in order to apply \Cref{prop_oscdec}.

This concludes the proof if $d_i = d$ for all $i$; if this is not the case, it suffices to take
\begin{equation} \label{eq_Tgen}
\call T \defeq \begin{pmatrix}
I_{d_0} 	& -\mu \beta_1 J_1	&0 & \dots 	& 0 \\
0		& I_{d_1}				& -\mu \beta_2 J_2 &\dots 	& 0 \\
\vdots 	& \vdots					& \ddots	& \ddots	 & \vdots \\
0 & 0 & \dots & I_{d_{r-1}} & - \mu \beta_r J_r \\
0		& 0						& \dots	& 0	& I_{d_r}
\end{pmatrix}
\end{equation}
instead and proceed along the same lines as above.   
\end{proof}

To extend this result and obtain \Cref{intro-thm1gen}, let us note that the structure assumptions therein essentially produce multiple, pairwise disjoint, basic H\"{o}rmander hierarchies, where by \emph{basic} we mean given by  \eqref{eq_HhormC1intro}. The underlying philosophy can be visually represented as follows: if we can separate the variables and obtain disjoint chains
\begin{equation} \label{eq_chainvis}
\begin{gathered}
x^1 \rightarrow x^{j^1_1} \rightarrow \dots \rightarrow x^{j^1_{r_1}} \\
\vdots \\
x^q \rightarrow x^{j^q_1} \rightarrow \dots \rightarrow x^{j^q_{r_q}}
\end{gathered}
\end{equation}
(with $r_1 + \dots + r_q + q = N$), where $x^1,\dots,x^q$ are diffusive directions, and such that we can read the structure \eqref{eq_HhormC1intro} when considering only the coordinates in each chain, then we will be able to deal with each chain ``separately''.

For example, consider $\call L = \call L_x - \call L_y$ in $\R^{N_1} \times \R^{N_2}$ with typical coordinates $(x,y)$, where $\call L_x$ and $\call L_y$ are operators respectively independent of $x$ and $y$, and satisfying the basic H\"{o}rmander hierarchy \eqref{eq_HhormC1intro} --- when restricted to $\R^{N_1}$ and $\R^{N_2}$, respectively; in this case the idea is that the use of an appropriate norm $\scr N$ essentially decouples the equation by allowing to deal with the contributions of $\call L_x$ and $\call L_y$ separately, when it comes to the analysis required to prove \Cref{thm_mainsmnl} (or \Cref{intro-thm1}).

As far as the time scaling is concerned, given that of \eqref{eq_scalintro}, in the scenario of multiple hierarchies as above, we expect to have a short-time scaling $(t^{s_0},\dots,t^{s_N})$ with corresponding exponents 
\begin{equation} \label{eq_scalexgenh00}
s_{j^k_h} = \bar\sigma_{\cL_k}^{-1} + h;
\end{equation}
that is, the chains of exponents relative to the coordinates in \eqref{eq_chainvis} will be
\begin{equation*}
\bar\sigma_{\cL_k}^{-1} \rightarrow \bar\sigma_{\cL_k}^{-1} + 1 \rightarrow \dots \rightarrow \bar\sigma_{\cL_k}^{-1}+ r_k, \qquad k \in \{1,\dots,q\}. 
\end{equation*}

\begin{proof}[{\bfseries Proof of \Cref{intro-thm1gen}}]
The corresponding general setting is $\ell = q$, $V_j= \{ x \in \R^N :\ x^{-(j,0)} = 0 \}$, and $\bar\nu_j = \lim_{\sigma\uto \bar\sigma_{\cL_j}} (2-\sigma)|\var|^{-d_{j,0}-\sigma} \scr H^{d_{j,0}}$.
Taking into account the above discussion, we take $\call T$ and $\scr N_\epsilon$ such that their restrictions to each of the subspaces $ \{ x \in \R^N :\ x^{-k} = 0 \}$, $k \in \{1,\dots,q\}$, have the structure in \eqref{eq_defrotcoord} and \eqref{eq_chNHorm}, respectively. Then, the computations in the proof of \Cref{intro-thm1} imply that the assumptions of \Cref{thm_mainsmnl,prop_oscdec} are readily verified, since one has to deal with direct sums of maps that verify them, with $s_k(t) \lesssim t$ and $\norm{DH^k_{\M}} \lesssim t^{-1}$ for all $k$. 
\end{proof}

\begin{rmk}\label{inter_catene}
It is possible to further weaken the assumptions of \Cref{intro-thm1gen}, by asking for instance that each $H^{k,0}$ and $H^{k,1}$ be independent of $x^{-k,-0}$ (instead of $x^{-k}$); that is, it may also depend on $x^{i,0}$ for any $i \in \{0,\dots,q\}$. In this case, $DH$ has no longer a direct sum structure that allows to replicate the proof of \Cref{intro-thm1} by separately dealing with each coordinate chain. In particular, letting
\[
\call S(t) = \bigoplus_{k=1}^q t^{\sigma_{k}^{-1}} \mathrm{diag}(t^iI_{d_{k,i}})_{0 \leq i \leq r_k}
\]
for some $\sigma_k \in (0,2]$, one has (for each $i,j,h,k$ with $h,k \in \{1,\dots,q\}$, $j \in \{0,\dots, r_k\}$, $i \in \{0,\dots, r_h\}$)
\[
s_k(t) \lesssim t^{\bar\sigma_{\cL_k}/\sigma_k}, \qquad (\call S^{-1} DH \call S)_{(h,i),(k,j)} = t^{j+ \sigma_k - i - \sigma_h} D_{x^{k,j}}H^{h,i};
\] 
therefore,
\[
\call S^{-1} \bigl( \dot{\call S} + DH\call S \bigr) \lesssim \frac1{t^{1+\max_{h,k}(\sigma_h-\sigma_k)}},
\]
so we see that the best choice that we are allowed to make is $\sigma_k = \min_h \bar\sigma_{\cL_h}$ for all $k$. This is clearly strongly related to the needed interplay between diffusion and drift to obtain Lipschitz smoothing in all directions. On the other hand, the verification of continue to hold as we have not added any dependence on $x^{i,j}$ with $j \neq 0$. We propose the next examples in order to illustrate and clarify the above constructions and the present remark.
\end{rmk} 

\begin{exa}[Choice of scaling, twist and norm with multiple H\"{o}rmander hierarchies] \label{ex_genhorm1}
In $(0,T) \times \R^5$, with typical space coordinates $(x^0,\dots,x^4)$, let
\be\label{ex-nov26}
\call L^H = \Delta_{x^0} - (-\Delta_{x^1})^{\frac\sigma2} + x^1D_{x^2} + x^0D_{x^3} + x^3 D_{x^4}.
\ee
Then we can identify two coordinate chains $x^0 \to x^3 \to x^4$ and $x^1 \to x^2$. We can rescale by
\begin{equation} \label{eq_scalexgenh1}
\call S(t) = \mathrm{diag}(t^{\frac12},t^{\frac1\sigma},t^{\frac1\sigma+1},t^{\frac32},t^{\frac52})
\end{equation}
and obtain
\[
\call L_{\call S}^H = \frac1t\bigl( \Delta_{x^0} - (-\Delta_{x^1})^{\frac\sigma2} \underbrace{{} + x^0 D_{x^0} + x^1 D_{x^1} + (x^1+x^2)D_{x^2} + (x^0+x^3)D_{x^3} + (x^3+x^4) D_{x^4}}_{\displaystyle-H_{\call S} \cdot D_x} \bigr)\,.
\]
Then we can use $\M(t) = \call S(t) \call T$, where the twist matrix is given by 
\[
\call T \defeq \begin{pmatrix}
1 & 0 & 0 & -\mu \beta_3 & 0 \\
0 & 1 & -\mu & 0 & 0 \\
0 & 0 & 1 & 0 & 0 \\
0 & 0 & 0 & 1 & -\mu \\
0 & 0 & 0 & 0 & 1
\end{pmatrix}
\]
and  consider the norm (neglecting the $\epsilon$-regularization for the sake of simplicity)
\[
\scr N(x) = |x^0| + |x^1| + \mu |x^2| + \alpha_3 \mu |x^3| + \alpha_4 \mu^2 |x^4|.
\]
Therefore, for $z = (0,0,z^2,z^3,z^4)$ we have $D\scr N(z) \in \{0\}\times\{0\} \times \R^3$ and 
\[
\call T^{-1} DH_{\call S} \call T z \cdot D\scr N(z) \leq -\mu^2(1+O(\mu^{-1}))|z^2| - \mu^2 \bigl( \alpha_3(\beta_3-1)-\alpha_4 + O(\mu^{-1})\bigr) |z^3| - \mu^3(\alpha_4 - 1+O(\mu^{-1})) |z^4|,
\]
so \eqref{eq_condcallN''} is verified, for instance,  upon choosing $\beta_3 = 2$, $\alpha_3 = 3$, $\alpha_4 = 2$, and $\mu$ large enough.

We notice that the same construction would apply if, in   \rife{ex-nov26}, we replace $\Delta_{x^1}$ with $\Delta_{(x^0,x^1)}$; indeed, \Cref{intro-thm1gen} does not exclude the possibility that some variable is affected by more than one of the diffusion operators. Of course, in this case we would continue to  use the scaling $t^{\frac12}$ for the variable $x^0$ thanks to the second order term in \rife{ex-nov26}.
\end{exa}

\begin{exa}[Choice of scaling with multiple H\"{o}rmander hierarchies, revisited]
Take
\[
\call L^H = \Delta_{x^0} - (-\Delta_{x^1})^{\frac\sigma2} + (x^0+x^1)D_{x^2} + x^0D_{x^3} + x^3 D_{x^4}.
\]
Note that here both $x^0$ and $x^1$ are transferring diffusivity to $x^2$, so we are in the situation described by \cref{inter_catene}. In that case, the scaling 
 \eqref{eq_scalexgenh1} would not work anymore; indeed, it would produce
\[
\call L^H_{\call S} = \frac1t \biggl( \Delta_{x^0} - (-\Delta_{x^1})^{\frac\sigma2} + \biggl(\, \frac1{t^{\frac1\sigma - \frac12}}\,x^0 + x^1+x^2\biggr)D_{x^2} + (x^0+x^3)D_{x^3} + (x^3+x^4) D_{x^4} \biggr),
\]
that does not satisfy \eqref{eq_condnormSH} as $\sigma < 2$. On the other hand, taking $\call S'(t) = \mathrm{diag}(t^{\frac1\sigma}I_2,t^{\frac1\sigma+1} I_2, t^{\frac1\sigma+2})$ gives
\[
\call L^H_{\call S'} = \frac1t \biggl( \frac1{t^{\frac2\sigma-1}} \Delta_{x^0} - (-\Delta_{x^1})^{\frac\sigma2} + (x^0 + x^1+x^2)D_{x^2} + (x^0+x^3)D_{x^3} + (x^3+x^4) D_{x^4} \biggr),
\]
for which \eqref{eq_condnormSH} holds. The rest of the argument would then be similar to that in \Cref{ex_genhorm1}.
\end{exa}

We conclude with the proof the the remaining statement of the Introduction, namely \Cref{sharp-Kolm}, which is concerned with the case that the drift fields $H^i$ in \rife{cascade} only depend on $x^{i-1}$, for $i\geq 1$. This shows that \Cref{thm_mainsmnl} can  provide, at the same time, sharp rates for both short-time bounded-to-Lipschitz smoothing estimates and long-time decay of derivatives.

\begin{proof}[{\bfseries Proof of \Cref{sharp-Kolm}}]
We can consider the general setting in \Cref{sec:genset}, with $\ell = 2$, $V_1 = V_2 = \R^{d_0} \times \{0\}$, $\bar\nu_1 = \lim_{\sigma\uto \bar\sigma_\cL} (2-\sigma)|\var|^{-d_0-\sigma} \scr H^{d_0}$, and $\bar\nu_2 = \lim_{\sigma\uto \underline\sigma_\cL} (2-\sigma)|\var|^{-d_0-\sigma} \scr H^{d_0}$.
We choose $\M(t) = \call S(t)\call T$, with $\call T$ as in \eqref{eq_Tgen} and 
\[
\call S(t) = t^{\tilde\sigma(t)^{-1}} \mathrm{diag}(t^{i}I_{d_i})_{0\leq i\leq r}
\]
where $\tilde\sigma \in C^1_b(\R_{\geq0};(0,2])$ is decreasing and satisfies
\[
\tilde\sigma([0,1]) = \{\bar\sigma_\cL\}, \qquad \tilde\sigma([2,+\infty)) = \{\underline\sigma_\cL\}.
\]
Following the proof of \Cref{intro-thm1}, and recalling \Cref{rmk_genSigma'}, we can take $\{W_j\}_{j \in I} = \{ V_j \}_{j \in I}$ with $I = I(t)$ such that $I((0,1]) = \{\{1\}\}$ and $I((1,+\infty)) = \{\{2\}\}$; we obtain that
\[
\max_{j \in I(t)} s_j(t) = s_1(t) \1_{t \in (0,1]} + s_2(t) \1_{t > 1} \lesssim t^{\tilde\sigma(t)} \quad \text{uniformly in $T$},
\]
hence $\max_{j \in I(t)} s_j(t) \norm{\call S(t)^{-1} \bigl( \dot{\call S}(t) + DH(t,\var)\call S(t) \bigr)}$ is bounded uniformly in $T$ as well. 
Hence, in this case, by arguing as in the proof of \Cref{intro-thm1}, we see that all hypotheses of \Cref{thm_mainsmnl} hold with constants independent of $T$. On the other hand, from the proof of \Cref{thm_mainsmnl}, we see that the only explicit dependence of $C_T$ on $T$ is through the constant in \eqref{eq_contrtLyap}; 
therefore, as $u$ is bounded and thus we can take $\phi = 1$, part (1) of \Cref{intro-thm1} holds with $C_T$ actually independent of $T$.
\end{proof}

Compared to \Cref{model_LP}, when we consider second order equations under H\"ormander conditions, the structural assumptions of \Cref{sharp-Kolm} are verified when $E_{ij} = 0$ for all $i,j$ in the decomposition in \Cref{rmk_cfm}; in other words, \Cref{sharp-Kolm} can be applied to the principal part of the H\"{o}rmander operator in \eqref{Hor-eq}; see \cite{LaPo}.

\smallskip

Finally, we conclude by showing how the setting we developed before can apply to models of heat conduction through anharmonic chain of oscillators, of the type investigated, e.g., in \cite{Eckmann+al,angeliki,Bellet}.

\begin{exa}[A class of chains of oscillators]\label{chains}
With reference to \cite{angeliki}, we consider the following system of SDEs (up to constants which we normalize to $1$ without loss of generality):
\[ \begin{dcases}
\di q^i_t = p^i_t \,\di t \qquad (1 \leq i \leq N) \\
\di p^i_t = - D_{q^i}\call H\,\di t \qquad (2 \leq i \leq N-1)\\
\di p^1_t = - (D_{q^1}\call H + p^1)\,\di t + \di W^1_t \\
\di p^N_t = - (D_{q^N}\call H + p^N)\,\di t + \di W^N_t,
\end{dcases}\]
where $q^i \in \R^d$ are the displacement of atoms with respect to their equilibrium positions, $p^i \in \R^d$ are their momenta, $W^1$ and $W^N$ are independent Brownian motions, and the Hamiltonian $\call H$ is given by
\[
\call H(p,q) = \sum_{i=1}^N \Bigl( \frac{|p^i|^2}2 + V(q^i) \Bigr) + \sum_{i=1}^{N-1} U(q^{i+1}-q^i).
\]
for some pinning potential $V$ and interaction potential $U$. Note that this amounts to having the equation
\begin{equation} \label{ang-pde}
\de_t u - \Delta_{p^1} u - \Delta_{p^N} u + p^1 \cdot D_{p^1} u + p^N \cdot D_{p^N}u - \sum_{i=1}^N \bigl( p^i \cdot D_{q^i} u - D_{q^i} \call H \cdot D_{p^i} u \bigr) = 0.
\end{equation}
More precisely, a \emph{small perturbation of the harmonic chain}, with $d=1$, is studied in \cite{angeliki}, where $V$ and $U$ are quadratic up to errors whose Hessians vanish as $N$ grows, so in particular $U$ is convex for $N$ large.

Such a chain of oscillators is a case where noise affects both the first and the last variables of the chain, so, differently from \eqref{cascade}, the diffusion(s) can propagate regularity from the first coordinate down to the last one, but also vice-versa. However, this is only apparently different, as we can pair the variables as
\[ 
X^i = \begin{dcases}
\bigl(p^{\frac i2 +1}, p^{N - \frac i2}\bigr) & \text{$i$ even} \\
\bigl(q^{\lfloor \frac{i}2 \rfloor +1}, q^{N - \lfloor \frac i2 \rfloor}\bigr) & \text{$i$ odd}
\end{dcases}
\] 
to get the same structure as \eqref{intro-chosc}, with in particular
\[
D_{x^{i-1}} H^i = \begin{dcases}
\mathrm{diag}\bigl(U''(q^{\frac i2 +1} - q^{\frac i2}), U''(q^{N-\frac i2+1} - q^{N-\frac i2}) \bigr) & \text{$i\geq2$ even} \\
I_2 & \text{$i$ odd}.
\end{dcases}
\]
Note that this construction works if $N$ is even (in which case $i \in \{0,\dots,N-1\}$), while for odd $N$ we need to let $X^{N-1} = p^{\lceil N/2 \rceil}$ and $X^N = q^{\lceil N/2 \rceil}$, so that $D_{x^{N-2}} H^{N-1} = \bigl(U''(q^{\lceil N/2 \rceil} - q^{\lfloor N/2 \rfloor}), U''(q^{\lceil N/2 \rceil +1} - q^{\lceil N/2 \rceil}) \bigr) \in \R^2$ and $D_{x^{N-1}} H^N = 1$.
In any case, we see that the structure assumptions of \Cref{intro-thm1} are satisfied (with in particular $J_i = I_2$ for all $i$ if $N$ is even, while $J_N = 1$ and, e.g., $J_{N-1} = \begin{psmallmatrix}1\\0\end{psmallmatrix}$ if $N$ is odd), so such a model falls into the setting covered by our result.

Interestingly, we have multiple ways to visualize how the hypoelliptic structure propagates the noise through the variables: one is that we have just described; another one consists in observing that the original structure with the variables $(p,q)$ indeed satisfies the assumptions of \Cref{intro-thm1gen}, if one chooses two disjoint chains of coordinates $(p^1,q^1,p^2,q^2,\dots,p^n,q^n)$ and $(p^N,q^N,p^{N-1},q^{N-1},\dots,p^{n+1},q^{n+1})$, where $n = \lceil \frac N2 \rceil$. This is an optimal splitting as the scaling provided by $\call S(t)$ is the same as that obtained by applying \Cref{intro-thm1} with respect to the previous paired coordinates. On the other hand, we see that we might choose any other value of $n$, up to the extreme case $n=N$, which corresponds to ignoring one of the two diffusions; however, choices of $n \neq \lceil \frac N2 \rceil, \lfloor \frac N2 \rfloor$ worsen the short-time scaling associated to certain coordinates, as they do not fully exploit the symmetry of the system.

It is also worth noting that when $N$ is odd, the two different choices of $n=\lceil \frac N2 \rceil$ or $n = \lfloor \frac N2 \rfloor$ directly correspond to the two possible choices of $J_{N-1} =  \begin{psmallmatrix}1\\0\end{psmallmatrix}$ or $J_{N-1} =  \begin{psmallmatrix}0\\1\end{psmallmatrix}$. Indeed, the non-degeneracy condition in \eqref{eq_HhormC1intro} essentially tells that noise propagates from $x^{i-1}$ to $x^i$ because (when $d_i < d_{i-1}$) we can embed $\R^{d_i}$ into $\R^{d_{i-1}}$ in such a way that $D_{x^{i-1}} H^i|_{\R^{d_i}}$ defines a positive definite quadratic form; under this viewpoint, the two different choices of $J_{N-1}$ above are equivalent to saying that we propagate regularity to $X^{N-1} = p^{\lceil \frac N2 \rceil}$ either from the first component of $X^{N-1}$, hence from $q^{\lceil \frac N2 \rceil - 1}$, or from the second one, hence from $q^{\lceil \frac N2 \rceil + 1}$.

Finally, we point out that an argument along the lines of that proposed at the end of Section~\ref{sec_BernM}, yielding exponential decay of $\norm{Du(t,\var)}_\infty$, also applies to the setting in \cite{angeliki}, since (up to small perturbations) the drift in \eqref{ang-pde} has constant Jacobian $B$ which is asymptotically stable. 
\end{exa}

\section{Long-time contractivity of  hypoelliptic Fokker--Planck equations}

We introduce the following Fokker--Planck equation:
\begin{equation} \label{eq_fpg}
(\de_t - \call L^*)m = 0, \qquad \call L^*m \defeq \tr(AD^2) + \call I^*[m] - \mathrm{div}(Hm) = 0,
\end{equation}
with $H \colon (0,T) \times \R^N \to \R^N$, $A$ as in \eqref{eq_condAt}, and where
\begin{equation} \label{eq_defcallIdual}
\call I^*[m](t,x) \defeq \int_{\R^N} \Bigl( m(t,x+z) - m(t,x) - Dm(t,x) \cdot z \1_{|z| \leq 1} \Bigr) \, (-I_N)_\sharp \nu(t,\di z),  \qquad \nu = \sum_{1 \leq j \leq \ell} \nu_j,
\end{equation}
with $\nu_j$ satisfying Assumptions~\ref{ass_nu}.
We consider the following notion of solution to \eqref{eq_fpg}, which stems from the duality with \eqref{eq_nleq}.

\begin{defn}[Duality solutions of Fokker--Planck equations] \label{def_dualsol}
Let $m_0 \in \scr M(\R^N)$, $m \in C^0([0,T];\scr M(\R^N))$.\footnote{The space $\scr M(\R^N)$ is endowed with the weak-$*$ topology.} We say that $m$ is a \emph{duality solution} to \eqref{eq_fpg} if, for all $\tau \in (0,T)$, $g \in C_b^0(\R^N)$, $f \in C^{0,1}_b((0,\tau) \times \R^N)$, and $u \in C^0_b([0,\tau] \times \R^N)$ viscosity solution to
\begin{equation} \label{eq_dualsys}
\begin{dcases}
(-\de_t - \call L)u = f & \text{in} \ (0,\tau) \times \R^N \\
u(\tau,\var) = g & \text{in} \ \R^N,
\end{dcases}
\end{equation}
the following duality equality holds:
\begin{equation} \label{eq_dueq}
\int_{\R^N} g \,\di m(\tau) + \int_0^\tau \int_{\R^N} f(t,\var) \,\di m(t)\di t = \int_{\R^N} u(0,\var) \,\di m_0.
\end{equation}
\end{defn}

\begin{rmk}[Duality vs.\ weak solutions] \label{rmk_dualweak}
Duality solutions satisfy in particular
\begin{equation} \label{eq_wsol}
\int_0^T \int_{\R^N} \bigl( \de_t + \call L)\phi \,\di m(t)\di t + \int_{\R^N} \phi(0,\var) \,\di m_0 = 0
\end{equation}
for all $\phi \in C^{1,3}_c([0,T) \times \R^N)$, hence by a standard approximation argument they are proved to be weak solutions. On the other hand, under fairly general assumptions, one proves existence and uniqueness of duality solutions --- see, e.g., \cite[Theorem~2.6]{PorAdv}; we point out that the key ingredient is the availability of an estimate like \eqref{eq_gthm_usmoothing} for viscosity solutions to $(\de_t - \call L)u = f$, which, by following the proof of \Cref{thm_mainsmnl}, we can get if $f\in C^{0,1}_b$, as required in \rife{eq_dualsys}.
\end{rmk}

The main result we state in this section comes as a direct consequence of \Cref{cor_imprdec}, by exploiting the natural duality between solutions to \eqref{eq_nleq} and \eqref{eq_fpg}.

\begin{prop}[Contractivity for duality solutions] \label{prop_contrfp}
Let $m^1,m^2$ be duality solutions to \eqref{eq_fpg}, where $\call I^*$ is given by \eqref{eq_defcallIdual} with $\nu_j$ satisfying Assumptions~\ref{ass_nu}. Suppose that the assumptions of \Cref{thm_mainsmnl,prop_oscdec} hold, and let us be given a smooth stationary Lyapunov function  $\phi$ satisfying \rife{eq_pshar},   
and such  that $\norm{\phi m^i(0)}_\TV < +\infty$ for $i=1,2$. Then
\begin{equation} \label{eq_TVdec}
\norm{\phi(m^1(t) - m^2(t))}_\TV \leq K e^{-\varpi t} \norm{\phi(m^1(0) - m^2(0))}_\TV \qquad \forall \, t \in (0,T],
\end{equation}
where $K$ and $\varpi$ are positive constants independent of $m^i$, $t$, and $T$.
\end{prop}

\begin{proof}
By using \eqref{eq_wsol} and the second property in \eqref{eq_condLyapdec}, we can infer that
\begin{equation} \label{eq_phiinmt}
\int_0^T \int_{\R^N}  \phi \,\di m^i(t)\di t \, \lesssim \, T  + \int_{\R^N} \phi\,\di m^i(0)
\end{equation}
(with implied constant only depending on $\phi$), thus $\phi \in L^1(\R^N,m^i(t))$ for all $t \in (0,T)$, $i = 1,2$. By \eqref{eq_dueq} with $f=0$, $m=m^1-m^2$, $m_0 = m^1(0) - m^2(0)$, and $\tau = t$,
\begin{equation} \label{eq_dueq'}
\int_{\R^N} g \,\di m(t) = \int_{\R^N} u(0,\var) \,\di m_0,
\end{equation}
whenever $g \in C^0_b(\R^N)$ and $u$ is a viscosity solution to \eqref{eq_dualsys} (with $f=0$). In fact, the estimate \eqref{eq_gthm_usmoothing} ensures the stability of viscosity solutions under approximation via bounded functions of $g$ in $\phi L^\infty(\R^N)$ (by standard arguments for which we refer, e.g., to \cite[Corollary~4.6]{ForPor}); hence \eqref{eq_dueq'} holds for any $g \in \phi L^\infty(\R^N)$ and their respective $u$, because we also proved that $\phi \in L^1(\R^N,|m|(t))$.
Furthermore, as $m_0$ as zero total mass by construction, we can replace $u(0,\var)$ with $u(0,\var) + c$ for any constant $c$. It follows that
\begin{equation*} 
\int_{\R^N} g \,\di m(t) \leq \inf_{c \in \R} \norm{\phi m_0}_\TV \norm{\phi^{-1}(u(0,\var)+c)}_\infty
= \norm{\phi m_0}_\TV [u(0,\var)]_{\phi,0}
\lesssim e^{-\varpi t} \norm{\phi m_0}_\TV \norm{\phi^{-1}g}_\infty,
\end{equation*}
where the last bound is due to \eqref{eq_decoscimp} (recalling that the equation for $u$ is now backwards in time), while we refer to \cite[Lemma~4.1]{PorAdv} for the proof of the identity $[u(0,\var)]_{\phi,0} = \inf_{c \in \R} \norm{\phi^{-1}(u(0,\var)+c)}_\infty$. By the arbitrariness of $g \in \phi L^\infty(\R^N)$, we obtain \eqref{eq_TVdec}. 
\end{proof}

It is evident that \Cref{prop_contrfp} can be equivalently formulated as an exponential decay (to zero) property for solutions to the Fokker--Planck equation with zero-mass initial data, as follows.

\begin{prop}[Decay for zero-mass duality solutions] \label{prop_decafp}
Given the assumptions of \Cref{prop_contrfp}, let $m$ be duality solution to \eqref{eq_fpg} with $m(0)(\R^N) = 0$
. Then
\begin{equation*}
\norm{\phi m(t)}_\TV \leq K e^{-\varpi t} \norm{\phi m(0)}_\TV \qquad \forall \, t \in (0,T],
\end{equation*}
where $K$ and $\varpi$ are positive constants independent of $m$, $t$, and $T$.
\end{prop}

On the other hand, an interesting application of the above contractivity concerns the existence of, and convergence to, stationary measures for \eqref{eq_fpg} (i.e., time-independent duality solutions), as noted in \cite[Section~5.1]{PorAdv}. It makes sense to talk about long-time convergence of duality solutions, because when $H$ and $\nu$ are time-independent then short-time existence of duality solution implies that they extend to $[0,+\infty)$; hence, in the upcoming statement, duality solutions to \eqref{eq_fpg} are understood to be maps $m \colon [0,+\infty) \to \scr M(\R^N)$ that satisfy \Cref{def_dualsol} for all $T>0$.

\begin{prop}[Convergence to stationary measures] \label{prop_convstatmeasg}
Given the assumptions of \Cref{prop_contrfp}, suppose in addition that $H$ and $\nu$ are time-independent. There exists a unique $\bar m \in \call P(\R^N)$ solving $\call L^* \bar m = 0$; also, for any duality solution $m$ to \eqref{eq_fpg} with $m_0 \in \call P(\R^N)$ and $\norm{\phi m(0)}_\TV < +\infty$, one has
\begin{equation} \label{eq_statconv}
\norm{\phi(m(t) - \bar m)}_\TV \leq K e^{-\varpi t} \norm{\phi m(0)}_\TV \qquad \forall\, t > 0,
\end{equation}
where $K$ and $\varpi$ are positive constants independent of $m$, $t$, and $T$.
\end{prop}

\begin{proof}
Given $m_0 \in \call P(\R^N)$, existence and uniqueness of a duality solution $m$ to \eqref{eq_fpg} with time-independent $H$ and $\nu$ follow by combining the arguments in \cite[Appendix]{PorAdv} with \cite[Lemma~4.9]{Espen+al}, which implies the existence of an unbounded sublinear supersolution to \eqref{eq_nleq} belonging to $L^\infty((0,T);L^1(\R^N;m_0))$. By \eqref{eq_dueq} with $g = 1$, $f=0$, and $u=1$, we know that $m(t) \in \call P(\R^N)$ for all $t>0$. Hence, by \Cref{prop_decafp}, for any $\tau, t > 0$ we have
\begin{equation} \label{eq_decprstat}
\norm{\phi(m(\tau+t) - m(t))}_\TV \lesssim e^{-\varpi t} \norm{\phi(m(\tau) - m(0))}_\TV.
\end{equation}
This implies that $m(k) = m(0) + \sum_{j=0}^{k-1} (m(j+1)-m(j))$
converges in $(\call P(\R^N),\norm{\phi \var}_\TV)$ as $k \uto \infty$; we call $\bar m$ the limit measure. As $m$ is in particular a duality solution on each time interval $[k-1,k]$, we know that, for any $g \in C^0_b(\R^N) \cap \phi L^\infty(\R^N)$ such that $\call L g \in C^{0,1}_b(\R^N) \cap \phi L^\infty(\R^N)$,
\[
\int_{\R^N} g \,\di(m(k)-m(k-1)) = \int_{k-1}^k \int_{\R^N} \call L g \,\di m(t) \di t;
\]
therefore,
\[
\int_{\R^N} \call L g \,\di m(k) = \int_{k-1}^k \int_{\R^N} \call L g \,\di(m(k)-m(t)) \di t + \int_{\R^N} g \,\di(m(k)-m(k-1)).
\]
By using again \eqref{eq_decprstat}, we see that
\[
\biggl\lvert \int_{\R^N} \call L g \,\di m(k) \biggr\rvert \lesssim \int_{k-1}^k e^{-\varpi t} \norm{\phi(m(k-t)-m(0))}_\TV \,\di t + e^{-\varpi k} \norm{\phi(m(1)-m(0))}_\TV \lesssim e^{-\varpi k},
\]
where the integral over $[k-1,k]$ is controlled by recalling \eqref{eq_phiinmt}. Hence, $\int_{\R^N} \call L g\,\di \bar m = 0$, yielding that $\bar m$ is a stationary duality solution of \eqref{eq_fpg}. The estimate \eqref{eq_statconv} is obtained from \eqref{eq_decprstat}, by letting $\tau \uto \infty$ and noting that $\norm{\phi( \bar m - m(0))}_\TV \lesssim_\varpi \norm{\phi( m(1) - m(0))}_\TV  \lesssim_{\varpi,\phi} \norm{\phi m(0)}_\TV$.
\end{proof}

As an immediate corollary, one has \Cref{intro-thm2}, presented in the Introduction; indeed, in \Cref{sec:ex} we showed that the assumptions of \Cref{intro-thm1gen} allow to verify those of \Cref{thm_mainsmnl,prop_oscdec}, and thus \Cref{prop_convstatmeasg} applies.

\end{document}